\documentclass[bibliography=totoc,DIV=12]{scrartcl}

\usepackage{amsmath}
 \numberwithin{equation}{section}
 \allowdisplaybreaks[1]
\usepackage{amssymb}
\usepackage{amsthm}

\theoremstyle{definition}

\theoremstyle{remark}

\usepackage[main=english,german]{babel}
\usepackage[noisbn,noissn,fixlanguage]{babelbib}%
 \setbtxfallbacklanguage{english}
\usepackage{cite}%
\usepackage[T1]{fontenc}
\usepackage{hyperref}
\usepackage{lmodern}
\usepackage{mathdots}
\usepackage{mathscinet}
 
\usepackage{mathtools}
 \mathtoolsset{mathic=true}
\usepackage{microtype}
\usepackage{pmat}
\usepackage{url}

\usepackage{178arxiv}

\title{Schur analysis of matricial Hausdorff moment sequences}
\author{Bernd Fritzsche \and Bernd Kirstein \and Conrad M\"adler}
\date{}

\begin{document}
\maketitle

\begin{abstract}
 We develop the algebraic instance of an algorithmic approach to the matricial Hausdorff moment problem on a compact interval \(\ab\) of the real axis.
 Our considerations are along the lines of the classical Schur algorithm and the treatment of the Hamburger moment problem on the real axis by Nevanlinna.
 More precisely, a transformation of matrix sequences is constructed, which transforms Hausdorff moment sequences into Hausdorff moment sequences reduced by 1 in length.
 It is shown that this transformation corresponds essentially to the left shift of the associated sequences of canonical moments.
 As an application, we show that a matricial version of the arcsine distribution can be used to characterize a certain centrality property of \tnnH{} measures on \(\ab\).
\end{abstract}

\begin{description}
 \item[Keywords:] matricial Hausdorff moment problem, Schur--Nevanlinna type algorithm, matricial canonical moments, \tnnH{} block Hankel matrices, moments of matrix-valued measures, arcsine distribution
 \item[2010 MSC:] 44A60
\end{description}

\section{Introduction}
 This paper continues our former investigations on matricial Hausdorff moment sequences, \tie{}, moment sequences of matrix measures on a compact interval \(\ab\).
 In~\zita{MR3775449}, we treated the extension problem for truncated Hausdorff moment sequences in a purely algebraic way, obtaining deeper insights into the structure of such sequences.
 This leads us in~\zita{MR3979701} to a parametrization of all matricial Hausdorff moment sequences, or equivalently of the whole moment space generated by matrix measures on \(\ab\), in terms of matricial canonical moments in the general case of not necessarily invertible associated \tbHms{}.
 To that end, we slightly modified the approach of Dette/Studden~\zita{MR1883272} by using Moore--Penrose inverses of matrices on the one hand and by symmetrizing their construction of matrix canonical moments on the other hand.
 The one-to-one correspondence between the matricial moment space and the space of sequences of matricial canonical moments, established in~\zitaa{MR3979701}{\cthm{6.30}{}}, suggests in a next step to develop a Schur--Nevanlinna type algorithm in the context of the matricial power moment problem on \(\ab\).
 The construction of that particular instance of this algorithm, which is connected to matricial Hausdorff moment sequences, is the content of the present paper.
 
 The classical Schur algorithm developed by Schur~\zitas{MR1580948,MR1580958} consists in the successive decomposition of a holomorphic function mapping the open unit disc \(\mathbb{D}\defeq\setaca{z\in\C}{\abs{z}<1}\) of the complex plane \(\C\) into the closed unit disc \(\overline{\mathbb{D}}\) (so-called Schur functions), using the H.~A.~Schwarz lemma.
 This gives rise to a one-to-one correspondence between the class of all Schur functions, or equivalently the class of corresponding sequences of Taylor coefficients at the origin \(z_0=0\) (so-called Schur sequences), and a certain class of sequences of complex numbers from \(\overline{\mathbb{D}}\).
 This concept can also be applied in the context of finite measures \(\sigma\) on the real axis \(\R\), where the class of functions in question consists of functions \(S\) defined on the open upper half-plane \(\Pi_+\defeq\setaca{z\in\C}{\Im z>0}\) of \(\C\) by \(S(z)\defeq\int_\R\rk{t-z}^\inv\sigma\rk{\dif t}\) (so-called Stieltjes transforms).
 These are holomorphic functions mapping \(\Pi_+\) into its closure \(\overline{\Pi_+}\).
 In the case of existing power moments of \(\sigma\) of any non-negative integer order, the associated Stieltjes transform \(S\) admits an asymptotic expansion at \(z_0=\infc\) in sectorial regions \(\setaca{z\in\C}{\delta<\arg z<\pi-\delta}\), \(0<\delta<\pi/2\), and an expansion as finite or infinite formal continued fraction.
 The coefficients occurring in this continued fraction can be identified with the real and imaginary part of points in \(\overline{\Pi_+}\) and are in a one-to-one correspondence to the sequence of power moments of \(\sigma\).
 The relevant investigations go back to Nevanlinna~\zita{zbMATH02604576}.
 For finite measures on the right half-line \([0,\infp)\), similar results are due to Stieltjes~\zitas{MR1623484,MR1607517,MR1344726,MR1344720}.
 What concerns matricial generalizations of this approach, we mention Kovalishina~\zita{zbMATH03875641}, Chen/Hu~\zita{MR1624548}, and Fritzsche~et~al.~\zita{MR3380267,MR3014199} for the real axis \(\R\) and Chen/Hu~\zita{MR1807884}, Hu/Chen~\zita{MR2038751}, and Fritzsche~et~al.~\zitas{MR3765778,MR3611479,MR3611471} for the half-line \([0,\infp)\).
 
 It seems to the authors of the present paper that there was not as yet any attempt to apply the above explained concept for finite measures on a compact interval \(\ab\), even in the scalar case.
 Following our procedure for \(\R\) in~\zitas{MR3014199,MR3380267} and for \(\rhl\) in~\zitas{MR3611479,MR3611471}, we are going to close this gap, first focusing in this paper on the algebraic aspect, \tie{}, on matricial Hausdorff moment sequences.
 Our strategy is mainly based on using our former studies in~\zitas{MR3775449,MR3979701} about the intrinsic structure of \tFnnd{} sequences.
 
 The paper is organized as follows:
 In \rsec{S.MP}, we recall the underlying setting of matricial power moment problems on the real axis.
 To construct certain transformations for sequences of complex matrices, we give in \rsec{S0858} generalizations of the \tCP{} and the \trF{} corresponding to multiplication and (generalized) inversion, \tresp{}, of formal power series with matrix coefficients.
 In \rsec{S1216}, we derive some identities reflecting the interplay between several \tbHms{} formed by a sequence of complex \tpqa{matrices} at the one side and by the corresponding \trF{} on the other side (see in particular \rthmsss{H.T1431}{H.P0836}{H.P0938}).
 \rsecsts{H.S}{F.S.abms} provide for later use several facts on the classes of moment sequences of matrix measures on \(\R\), on \(\rhl\) and \(\lhl\), and on \(\ab\), partly taken from former work.
 In particular, the so-called \tfc{s} \(\cia{0},\cia{1},\cia{2},\dotsc\) from~\zitaa{MR3979701}{\csec{6}}, which establish a one-to-one correspondence between matricial \tFnnd{} sequences (which proved to be exactly the moment sequences of matrix measures on \(\ab\)) and a certain subset of sequences from the closed matricial unit interval \(\matint{\Oqq}{\Iq}\) are recalled.
 Using the constructions from \rsec{S.MP}, a transformation, called \tFT{ation}, of sequences of complex matrices is given in \rsec{F.S.FT}, which reduces the length of a finite sequence by one.
 We obtain representations for the four relevant types of \tbHms{} associated to the \tFT{ed} sequence and the interval \(\ab\) in terms of corresponding \tbHms{} (or Schur complements) built from the input sequence.
 To that end, some lengthy calculations are needed, which constitute the main part of the here presented work.
 The above mentioned representations are then used in \rsec{S.F.kFT} to verify that \tFnnd{ness} is respected by the \tFT{ation}.
 The iterated \tFT{ation} of an \tFnnd{} sequence gives rise to block \(LDU\)~factorizations of the four associated types of \tbHms{}, revealing connections between the transformed sequences and certain Schur complements \(\fpu{0},\fpu{1},\fpu{2},\dotsc\) already investigated in~\zitaa{MR3979701}{\csec{6}}.
 In our main result \rthm{F.T.SalgIP}, it is shown that the \tFT{ation} of an \tFnnd{} sequence corresponds essentially to the left shift of its \tfcf{} \(\cia{0},\cia{1},\cia{2},\dotsc\), justifying an interpretation in the context of Schur analysis.
 Since a matrix measure \(\sigma\) on the compact interval \(\ab\) is completely determined by the sequence of all its power moments of non-negative integer order, in \rsec{F.CM} we can apply the \tFT{ation} for such measures via its moment sequences as well.
 In terms of \tFT{s} \(\sigma^\FTa{k}\), we characterize for a given matrix measure \(\sigma\) on \(\ab\) the situation of having finite support, which corresponds to \tabCD{} of its moment sequence.
 Similarly, the case of \tZm{ity}, which is defined by having an \tabZ{} moment sequence, is shown to be equivalent to the condition that one of the \tFT{s} \(\sigma^\FTa{k}\) of \(\sigma\) essentially coincides with the arcsine distribution.

\section{Matricial moment problems on Borel subsets of the real axis}\label{S.MP}
 In this section, we are going to formulate a class of matricial power moment problems.
 Before doing this we have to introduce some terminology.
 We denote by \(\Z\) the set of all integers.
 Let \(\N\defeq\setaca{n\in\Z}{n\geq1}\).
 Furthermore, we write \(\R\) for the set of all real numbers and \(\C\) for the set of all complex numbers.
 In the whole paper, \(p\) and \(q\) are arbitrarily fixed integers from \(\N\).
 We write \(\Cpq\) for the set of all complex \tpqa{matrices} and \(\Cp\) is short for \(\Coo{p}{1}\).
 When using \(m,n,r,s,\dotsc\) instead of \(p,q\) in this context, we always assume that these are integers from \(\N\).
 We write \(A^\ad\) for the conjugate transpose of a complex \tpqa{matrix} \(A\).
 Denote by \(\Cggq\defeq\setaca{M\in\Cqq}{v^\ad Mv\in[0,\infp)\text{ for all }v\in\Cq}\) the set of \tnnH{} matrices from \(\Cqq\).

 Let \((\mathcal{X},\mathfrak{X})\) be a measurable space.
 Each countably additive mapping whose domain is \(\mathfrak{X}\) and whose values belong to \(\Cggq\) is called a \tnnH{} \tqqa{measure} on \((\mathcal{X}, \mathfrak{X})\).
 For the integration theory with respect to \tnnH{} measures, we refer to Kats~\zita{MR0080280} and Rosenberg~\zita{MR0163346}.
 
 Let \(\BsAR\) (\tresp{}\ \(\BsAC \)) be the \(\sigma\)\nobreakdash-algebra of all Borel subsets of \(\R\) (\tresp{}\ \(\C \)).
 In the whole paper, \(\Omega\) stands for a non-empty set belonging to \(\BsAR\).
 Let \(\BsAO \) be the \(\sigma\)\nobreakdash-algebra of all Borel subsets of \(\Omega\) and let \(\MggqO \) be the set of all \tnnH{} \tqqa{measures} on \((\Omega,\BsAO )\).
 Observe that \(\Mggoa{1}{\Omega}\) coincides with the set of ordinary measures on \((\Omega,\BsAu{\Omega})\) with values in \([0,\infp)\).

 Let \(\NO\defeq\setaca{m\in\Z}{m\geq0}\).
 In the whole paper, \(\kappa\) is either an integer from \(\NO\) or \(\infi\).
 In the latter case, we have \(2\kappa\defeq\infi\) and \(2\kappa+1\defeq\infi\).
 Given \(\upsilon,\omega\in\R\cup\set{-\infty,\infp}\), denote \(\mn{\upsilon}{\omega}\defeq\setaca{k\in\Z}{\upsilon\leq k\leq\omega}\).
 Let \(\Mggoua{q}{\kappa}{\Omega}\) be the set of all \(\mu\in\MggqO \) such that for each \(j\in\mn{0}{\kappa}\) the power function \(x\mapsto x^j\) defined on \(\Omega\) is integrable with respect to \(\mu\).
 If \(\mu\in\MggquO{\kappa}\), then, for all \(j\in\mn{0}{\kappa}\), the matrix
\beql{I.G.mom}
 \mpm{\mu}{j}
 \defeq\int_\Omega x^j\mu\rk{\dif x}
\eeq
 is called (power) \emph{moment of \(\mu\)} of order \(j\).
 Obviously, we have \(\Mggoa{q}{\Omega}=\Mggoua{q}{0}{\Omega}\subseteq \MggquO{\ell}\subseteq \MggquO{\ell+1}\subseteq\MggquO{\infi}\) for every choice of \(\ell\in\NO \) and, furthermore, \(\mpm{\mu}{0}=\mu\rk{\Omega}\) for all \(\mu\in\Mggoa{q}{\Omega}\).
 If \(\Omega\) is bounded, then one can easily see that \(\MggqO =\MggquO{\infi}\).
 In particular, for every choice of \(\ug,\obg\in\R\) with \(\ug<\obg\), we have \(\MggqF=\MggquF{\infi}\).
 We now state the general form of the moment problem lying in the background of our considerations:

\begin{Problem}[\mprob{\Omega}{\kappa}{=}]
 Given a sequence \(\seqska \) of complex \tqqa{matrices}, parametrize the set \(\MggqOsg{\kappa}\) of all \(\sigma\in\Mggoua{q}{\kappa}{\Omega}\) satisfying \(\mpm{\sigma}{j}=\su{j}\) for all \(j\in\mn{0}{\kappa}\).
\end{Problem}

 In this paper, we mainly focus on the case that \(\Omega\) is a compact interval \(\ab\) of the real axis \(\R\).
 Each solution of Problem~\mprob{\ab}{m}{=} generates in a natural way solutions to each of the problems \mprob{\rhl}{m}{=}, \mprob{\lhl}{m}{=}, and \mprob{\R}{m}{=}.
 The last mentioned three matricial moment problems were intensively studied in our former work (see~\zitas{MR2735313,MR2570113,MR2805417,MR3014201,MR3133464,MR3611479,MR3014199}).
 In particular, we analyzed the inner structure of matricial moment sequences associated with each of the sets \(\rhl\), \(\lhl\), and \(\R\) and obtained parametrizations of the corresponding moment spaces.
 In order to prepare the investigations of the present paper, we recall some of this material in \rsecss{H.S}{K.S}.

\section{Invertibility of sequences of complex matrices and the concept of the \trF{}}\label{S0858}
 In this section, we continue our investigations on the concept of invertibility of finite or infinite sequences of complex matrices which we started in~\zita{MR3014197}.
 The explicit computation of the inverse sequence associated with an invertible sequence lead us in~\zita{MR3014197} to the notion of the \trFa{a} given sequence from \(\Cpq\).
 A closer look at the construction of the \trF{} shows us now that important properties of this objects which were proved in~\zita{MR3014197} for the subclass of invertible sequences hold in the most general case.
 
 In this paper, the Moore--Penrose inverse of a complex matrix plays an important role.
 For each matrix \(A\in\Cpq\), there exists a unique matrix \(X\in\Cqp\), satisfying the four equations
\begin{align}\label{mpi}
 AXA&=A,&
 XAX&=X,&
 \rk{AX}^\ad&=AX,&
&\text{and}&
 \rk{XA}^\ad&=XA
\end{align}
 (see \teg{}~\zitaa{MR1152328}{\cprop{1.1.1}{}}).
 This matrix \(X\) is called the \emph{Moore--Penrose inverse of \(A\)} and is denoted by \(A^\mpi\).
 Concerning the concept of Moore--Penrose inverse we refer to~\zita{MR0338013},~\zitaa{MR1105324}{\cch{1}}, and~\zitaa{MR1987382}{\cch{1}}.
 For our purposes, it is convenient to apply~\zitaa{MR1152328}{\csec{1.1}}.
 
 Now we give some terminology.
 Let \(\Opq\) be the zero matrix from \(\Cpq\).
 Sometimes, if the size is clear from the context, we will omit the indices and write \(\NM\).
 Let \(\Iq\defeq\matauuo{\Kronu{jk}}{j,k}{1}{q}\) be the identity matrix from \(\Cqq\), where \(\Kronu{jk}\) is the Kronecker delta.
 Sometimes, if the size is clear from the context, we will omit the index and write \(\EM\).
 If \(A\) is a square matrix, then denote by \(\det A\) the determinant of \(A\).
 To emphasize that a certain (block) matrix \(X\) is built from a sequence \(\seqska\), we sometimes write \(X^{\ok{s}}\) for \(X\).
 
\bnotal{M.N.S}
 Let \(\seqska\) be a sequence of complex \tpqa{matrices}.
 For all \(m\in\mn{0}{\kappa}\), then let the block Toeplitz matrices \(\SLuo{m}{s}\) and \(\SUuo{m}{s}\) be given by
\[
 \SLuo{m}{s}
 \defeq
 \bMat
  s_0    &\NM        &\NM      & \hdots    &\NM   \\
  s_1    & s_0      &\NM      & \hdots    &\NM   \\
  s_2    & s_1      & s_0    & \hdots    &\NM   \\
  \vdots & \vdots   & \vdots & \ddots    & \vdots  \\
  s_m    & s_{m - 1} & s_{m - 2}& \hdots    & s_0
 \eMat\text{ and }
 \SUuo{m}{s}
 \defeq
 \bMat
   s_0    & s_1      & s_2    & \hdots    & s_m   \\
   \NM             & s_0      & s_1    & \hdots    & s_{m-1}   \\
   \NM             & \NM               & s_0    & \hdots    & s_{m-2}  \\
   \vdots        & \vdots          & \vdots        & \ddots    & \vdots  \\
   \NM             & \NM               & \NM             & \hdots    & s_0
 \eMat.
\]
 Whenever it is clear which sequence is meant, we will simply write \(\SLu{m}\) (\tresp{}\ \(\SUu{m}\)) instead of \(\SLuo{m}{s}\) (\tresp{}\ \(\SUuo{m}{s}\)).
\enota

 If \(m\in\NO\), then a sequence \(\seqs{m}\) of complex \tpqa{matrices} is called \emph{\tinv{}} if there is a sequence \(\seqa{r}{m}\) of complex \tqpa{matrices} such that the Moore--Penrose inverse \(\rk{\SLuo{m}{s}}^\mpi\) of the matrix \(\SLuo{m}{s}\) coincides with the block Toeplitz matrix \(\SLuo{m}{r}\).

\breml{R2-3}
 If \(m\in\NO\) and \(\seqs{m}\) is an \tinv{} sequence of complex \tpqa{matrices}, then there is a unique sequence \(\seqa{r}{m}\) of complex \tqpa{matrices} such that
\beql{G2-2}
 \rk{\SLuo{m}{s}}^\mpi
 =\SLuo{m}{r}.
\eeq
 Recalling~\zitaa{MR3014197}{\ccor{A.11}{}}, we see that this sequence, \(\seqa{r}{m}\), satisfies \(\rk{\SLuo{\ell}{s}}^\mpi=\SLuo{\ell}{r}\) for each \(\ell\in\mn{0}{m}\).
\erem 

 \rrem{R2-3} shows that it is reasonable that a 
  sequence \(\seqsinf\) of complex \tpqa{matrices} is called \emph{\tinv{}} if there is a sequence, \(\seqa{r}{\infi}\), of complex \tqpa{matrices} such that \eqref{G2-2} holds true for each non-negative integer \(m\).
 Given a \(\kappa\in\NOinf\), we will use \(\Ipqka\) to denote the set of all \tinv{} sequences, \(\seqska\), of complex \tpqa{matrices}.
 If \(\kappa\in\NOinf\) and if \(\seqska\in\Ipqka\), then the unique sequence \(\seqa{r}{\kappa}\) with \eqref{G2-2} for each \(m\in\mn{0}{\kappa}\), is called the \emph{\tiFa{\(\seqska\)}} and will be denoted by \(\seqa{s^\invF}{\kappa}\).

\breml{R2-7}
 Let \(\kappa\in\NOinf\)  and \(\seqska\in\Ipqka\).
 From \rrem{R2-3} we see that \(\su{0}^\invF=\su{0}^\mpi\) and that \(\seqs{m}\in\Ipqu{m}\) for all \(m\in\mn{0}{\kappa}\).
 Furthermore, for all \(m\in\mn{0}{\kappa}\), the \tiFa{\(\seqs{m}\)} is exactly \(\seqa{s^\invF}{m}\).
\erem

\breml{R2-8}
 Let \(\kappa\in\NOinf\) and \(\seqska\in\Ipqka\).
 Using \rrem{A.R.A++*}, we see that \(\seqa{s^\invF}{\kappa}\in\Iqpka\)  and, furthermore, that \(\seqska\) is the \tiFa{\(\seqa{s^\invF}{\kappa}\)}.
\erem

 The main focus in~\zitaa{MR3014197}{\csec{4}} was firstly on coming to an explicit description of the set \(\Ipqka\) of all \tinv{} sequences \(\seqska\) in \(\Cpq\) and secondly on finding an effective method for constructing the \tiFa{a} sequence \(\seqska\in\Ipqka\).
 
 In the following, we write \(\ran{A}\defeq\setaca{Ax}{x\in\Cq}\) and \(\nul{A}\defeq\setaca{x\in\Cq}{Ax=\Ouu{p}{1}}\) for the column space and the null space of a complex \tpqa{matrix} \(A\), \tresp{}
 If \(\kappa\in\NOinf\) and if \(\seqska\) is a sequence in \(\Cpq\).
 We then say that \(\seqska\) is \emph{dominated by its first term} (or, simply, that it is \emph{\tftd{}}) when
\begin{align*}%
 \bigcup_{j=0}^\kappa\ran{s_j}&\subseteq\ran{s_0}&
&\text{and}&
 \nul{s_0}&\subseteq\bigcap_{j=0}^\kappa\nul{s_j}.
\end{align*}
 The set of all \tftd{} sequences \(\seqska\) in \(\Cpq\) will be denoted by \(\Dpqka\).
 We recall~\zitaa{MR3014197}{\cdefn{4.13}{24}}:
 
\bdefnl{D1419}
 Let \(\seqska \) be a sequence of complex \tpqa{matrices}.
 Then we call the sequence \(\seq{s_j^\rez}{j}{0}{\kappa}\) defined by \(s_0^\rez\defeq s_0^\mpi\) and, for all \(j\in\mn{1}{\kappa}\), recursively by
\(
 s_j^\rez
 \defeq-s_0^\mpi\sum_{\ell=0}^{j-1}s_{j-\ell}s_\ell^\rez
\)
 the \emph{\trFa{\(\seq{s_j}{j}{0}{\kappa}\)}}.
 For each (block) matrix \(X\) built from the sequence \(\seqska\), we denote by \(X^\rez\) the corresponding matrix built from the \trFa{\(\seqska\)} instead of \(\seqska\).
\edefn
 

\breml{M.R.reztr}
 Let \(\seqska\) be a sequence of complex \tpqa{matrices} with \trF{} \(\seqr{\kappa}\).
 It is readily checked that, for each \(m\in\mn{0}{\kappa}\), the \trFa{\(\seqs{m}\)} coincides with \(\seqr{m}\).
\erem

 A first observation associated with the \trF{} is the following.
 
\bpropnl{\zitaa{MR3014197}{\cprop{5.10(a)}{}}}{P1210}
 If \(\kappa\in\NOinf\) and \(\seqska\) is a sequence from \(\Cqq\) then \(\seqa{s^\rez}{\kappa}\in\Dqpka\).
\eprop

 The following result emphasizes the relevance of the set \(\Dpqka\).
 
\bthmnl{\zitaa{MR3014197}{\cthm{4.21}{26}}}{T4-21}
 If \(\kappa\in\NOinf\), then \(\Ipqka=\Dpqka\).
 Furthermore, if \(\seqska \in \Ipqka\), then the \tiFa{\(\seqska\)} coincides with the \trFa{\(\seqska\)}.
\ethm
 
 Against to the background of \rthm{T4-21} the investigation of the \trF{} sequence in~\zita{MR3014197} was mainly concentrated on the particular case of \tftd{} sequences.
 Our next considerations are aimed to verify that several results for \tftd{} sequences obtained in~\zita{MR3014197} remain true in the most general case.
 The key to this observation comes from a particular representation formula for the \trF{} from which we obtain more insight into the nature of this object.
 The application of \rdefn{D1419} immediately yields
\begin{align}\label{A.G.srez012}
 \su{0}^\rez&=\su{0}^\mpi,&
 \su{1}^\rez&=-\su{0}^\mpi\su{1}\su{0}^\mpi,&
 \su{2}^\rez&=-\su{0}^\mpi\su{2}\su{0}^\mpi+\su{0}^\mpi\su{1}\su{0}^\mpi\su{1}\su{0}^\mpi,
\end{align}
 and \(
 \su{3}^\rez
  =-\su{0}^\mpi\su{3}\su{0}^\mpi+\su{0}^\mpi\su{1}\su{0}^\mpi\su{2}\su{0}^\mpi+\su{0}^\mpi\su{2}\su{0}^\mpi\su{1}\su{0}^\mpi-\su{0}^\mpi\su{1}\su{0}^\mpi\su{1}\su{0}^\mpi\su{1}\su{0}^\mpi\).
 Guided by these identities we are led to an expression for the \trF{} which is not given via recursions. 
 For all \(\ell,m\in\N\) with \(\ell\leq m\), denote by \(\pauu{\ell}{m}\) the set of all ordered \(\ell\)\nobreakdash-tuples \((k_1,k_2,\dotsc,k_\ell)\) with \(k_1,k_2,\dotsc,k_\ell\in\N\) satisfying \(k_1+k_2+\dotsb+k_\ell=m\).
 
\bthml{P1414}
 Let \(\seqska \) be a sequence of complex \tpqa{matrices}. 
 Then the sequence \(\seqa{t}{\kappa}\) given by \(t_0\defeq s_0^\mpi\) and, for all \(j\in\mn{1}{\kappa}\), by
\[
 t_j
 \defeq\sum_{\ell=1}^j(-1)^\ell\sum_{(k_1,k_2,\dotsc,k_\ell)\in\pauu{\ell}{j}}\su{0}^\mpi\su{k_1}\su{0}^\mpi\su{k_2}\dotsm\su{0}^\mpi\su{k_\ell}\su{0}^\mpi
\]
 coincides with the \trFa{\(\seq{s_j}{j}{0}{\kappa}\)}.
\ethm
\bproof
 We proceed by mathematical induction.
 First observe that \(\su{0}^\rez=\su{0}^\mpi=t_0\).
 If \(\kappa\geq1\), then \(\su{1}^\rez=-s_0^\mpi\su{1}\su{0}^\mpi=t_1\), according to \eqref{A.G.srez012} and \(\pauu{1}{1}=\set{(1)}\).
 Now assume \(\kappa\geq2\).
 Then there exists an \(m\in\mn{2}{\kappa}\) such that \(\su{j}^\rez=t_j\) for all \(j\in\mn{0}{m-1}\).
 For all  \(k\in\mn{1}{m-1}\) and all \(\ell\in\mn{k}{m-1}\) let
 \[
  \mathcal{I}_{k,\ell}
  \defeq\setaca*{(m-\ell,i_1,i_2,\dotsc,i_k)}{(i_1,i_2,\dotsc,i_k)\in\pauu{k}{\ell}}.
 \]
 For all \(k\in\mn{1}{m-1}\), then \(\pauu{k+1}{m}=\bigcup_{\ell=k}^{m-1}\mathcal{I}_{k,\ell}\), where the sets \(\mathcal{I}_{k,k},\mathcal{I}_{k,k+1},\dotsc,\mathcal{I}_{k,m-1}\) are pairwise disjoint.
 Consequently, we obtain
 \[
  \begin{split}
   -\sum_{\ell=1}^{m-1}s_{m-\ell}t_\ell
   &=-\sum_{\ell=1}^{m-1}s_{m-\ell}\rk*{\sum_{k=1}^{\ell}(-1)^k\sum_{(i_1,i_2,\dotsc,i_k)\in\pauu{k}{\ell}}\su{0}^\mpi\su{i_1}\su{0}^\mpi\su{i_2}\dotsm\su{0}^\mpi\su{i_k}\su{0}^\mpi}\\
   &=\sum_{\ell=1}^{m-1}\sum_{k=1}^{\ell}(-1)^{k+1}\sum_{(i_1,i_2,\dotsc,i_k)\in\pauu{k}{\ell}}s_{m-\ell}\su{0}^\mpi\su{i_1}\su{0}^\mpi\su{i_2}\dotsm\su{0}^\mpi\su{i_k}\su{0}^\mpi\\
   &=\sum_{k=1}^{m-1}\sum_{\ell=k}^{m-1}(-1)^{k+1}\sum_{(i_1,i_2,\dotsc,i_k)\in\pauu{k}{\ell}}s_{m-\ell}\su{0}^\mpi\su{i_1}\su{0}^\mpi\su{i_2}\dotsm\su{0}^\mpi\su{i_k}\su{0}^\mpi \\
   &=\sum_{k=1}^{m-1}(-1)^{k+1}\sum_{\ell=k}^{m-1}\sum_{(i_0,i_1,i_2,\dotsc,i_k)\in\mathcal{I}_{k,\ell}}\su{i_0}\su{0}^\mpi\su{i_1}\su{0}^\mpi\su{i_2}\dotsm\su{0}^\mpi\su{i_k}\su{0}^\mpi\\
    &=\sum_{k=1}^{m-1}(-1)^{k+1}\sum_{(i_0,i_1,\dotsc,i_{k})\in\pauu{k+1}{m}}\su{i_0}\su{0}^\mpi\su{i_1}\su{0}^\mpi\su{i_2}\dotsm\su{0}^\mpi s_{i_{k}}s_0^\mpi\\
   &=\sum_{\ell=2}^{m}(-1)^\ell\sum_{(n_1,n_2,\dotsc,n_\ell)\in\pauu{\ell}{m}}\su{n_1}\su{0}^\mpi\su{n_2}\su{0}^\mpi\su{n_3}\dotsm\su{0}^\mpi s_{n_\ell}s_0^\mpi.
  \end{split}
 \]
 Taking additionally into account \rdefn{D1419} and \(\pauu{1}{m}=\set{(m)}\), we infer then
\[
 \begin{split}
  \su{m}^\rez
  &=-s_0^\mpi\sum_{\ell=0}^{m-1}s_{m-\ell}s_\ell^\rez
  =-s_0^\mpi s_{m}s_0^\rez+s_0^\mpi\rk*{-\sum_{\ell=1}^{m-1}s_{m-\ell}s_\ell^\rez}\\
  &=-s^\mpi_0s_m s_0^\mpi+s_0^\mpi\rk*{\sum_{\ell=2}^{m}(-1)^\ell\sum_{(n_1,n_2,\dotsc,n_\ell)\in\pauu{\ell}{m}}\su{n_1}\su{0}^\mpi\su{n_2}\su{0}^\mpi\su{n_3}\dotsm\su{0}^\mpi s_{n_\ell}s_0^\mpi}\\
  &=(-1)^1\sum_{(n_1)\in\pauu{1}{m}}s^\mpi_0s_{n_1}s_0^\mpi+\sum_{\ell=2}^{m}(-1)^\ell\sum_{(n_1,n_2,\dotsc,n_\ell)\in\pauu{\ell}{m}}s_0^\mpi\su{n_1}\su{0}^\mpi\su{n_2}\su{0}^\mpi\su{n_3}\dotsm\su{0}^\mpi s_{n_\ell}s_0^\mpi\\
  &=\sum_{\ell=1}^m(-1)^\ell\sum_{(n_1,n_2,\dotsc,n_\ell)\in\pauu{\ell}{m}}s_0^\mpi\su{n_1}\su{0}^\mpi\su{n_2}\su{0}^\mpi\su{n_3}\dotsm\su{0}^\mpi s_{n_\ell}s_0^\mpi
  =t_{m}.\qedhere
 \end{split}
\]
\eproof

\bcorl{C1511}
 Let \(\seqska\) be a sequence of complex \tpqa{matrices}.
 For all \(k\in\mn{0}{\kappa}\) then \(\su{0}^\mpi\su{0}\su{k}^\rez=\su{k}^\rez\) and \(\su{k}^\rez\su{0}\su{0}^\mpi=\su{k}^\rez\) as well as \(\ran{\su{k}^\rez}\subseteq\ran{\su{0}^\mpi}\) and \(\nul{\su{0}^\mpi}\subseteq\nul{\su{k}^\rez}\).
\ecor
\bproof
 Apply \rthm{P1414} and \rrem{R.AA+B}.
\eproof

 It should be mentioned that \rcor{C1511} can also be obtained as an immediate consequence of \rprop{P1210} (see~\zitaa{MR3014197}{\cprop{5.10(b)}{}} and the corresponding proof there).

 We will use \rthm{P1414} to see that several results from~\zita{MR3014199} hold in full generality without any additional assumptions on the sequence \(\seqs{\kappa}\) made there.
 This will facilitate the anyhow lengthy calculations in the proofs of \rpropssss{F.P.FTH}{F.P.FTHa}{F.P.FTHb}{F.P.FTHab}, which are crucial for our main results.
 
 Beside its definition, the \trF{} can be also computed in a dual recursive way:

\bpropl{P0731}
 Suppose \(\kappa\geq1\) and let \(\seqska\) be a sequence from \(\Cpq\).
 For all \(m\in\mn{1}{\kappa}\), then
\(
 \su{m}^\rez
 =-\rk{\sum_{\ell=1}^ms^\rez_{m-\ell}s_\ell}s_0^\mpi.
\)
\eprop
\bproof
 We have \(\su{0}^\rez=\su{0}^\mpi\) and \(\su{1}^\rez=-s_0^\mpi\su{1}\su{0}^\mpi\), according to \eqref{A.G.srez012}, implying the asserted formula for \(m=1\).
 Now assume \(\kappa\geq2\) and consider an arbitrary \(m\in\mn{2}{\kappa}\).
 For all  \(k\in\mn{1}{m-1}\) and all \(\ell\in\mn{1}{m-k}\) let
 \[
  \mathcal{J}_{k,\ell}
  \defeq\setaca*{(j_1,j_2,\dotsc,j_k,\ell)}{(j_1,j_2,\dotsc,j_k)\in\pauu{k}{m-\ell}}.
 \]
 For all \(k\in\mn{1}{m-1}\), then \(\pauu{k+1}{m}=\bigcup_{\ell=1}^{m-k}\mathcal{J}_{k,\ell}\), where the sets \(\mathcal{J}_{k,1},\mathcal{J}_{k,2},\dotsc,\mathcal{J}_{k,m-k}\) are pairwise disjoint.
 Using \rthm{P1414}, we consequently obtain
 \[
  \begin{split}
   -\sum_{\ell=1}^{m-1}s^\rez_{m-\ell}s_\ell
   &=-\sum_{\ell=1}^{m-1}\rk*{\sum_{k=1}^{m-\ell}(-1)^k\sum_{(j_1,j_2,\dotsc,j_k)\in\pauu{k}{m-\ell}}\su{0}^\mpi\su{j_1}\su{0}^\mpi\su{j_2}\dotsm\su{0}^\mpi\su{j_k}\su{0}^\mpi}s_\ell\\
   &=\sum_{\ell=1}^{m-1}\sum_{k=1}^{m-\ell}(-1)^{k+1}\sum_{(j_1,j_2,\dotsc,j_k)\in\pauu{k}{m-\ell}}\su{0}^\mpi\su{j_1}\su{0}^\mpi\su{j_2}\dotsm\su{0}^\mpi\su{j_k}\su{0}^\mpi s_\ell\\
   &=\sum_{k=1}^{m-1}\sum_{\ell=1}^{m-k}(-1)^{k+1}\sum_{(j_1,j_2,\dotsc,j_k)\in\pauu{k}{m-\ell}}\su{0}^\mpi\su{j_1}\su{0}^\mpi\su{j_2}\dotsm\su{0}^\mpi\su{j_k}\su{0}^\mpi s_\ell\\
   &=\sum_{k=1}^{m-1}(-1)^{k+1}\sum_{\ell=1}^{m-k}\sum_{(j_1,j_2,\dotsc,j_k,j_{k+1})\in\mathcal{J}_{k,\ell}}\su{0}^\mpi\su{j_1}\su{0}^\mpi\su{j_2}\dotsm\su{0}^\mpi\su{j_k}\su{0}^\mpi s_{j_{k+1}}\\
    &=\sum_{k=1}^{m-1}(-1)^{k+1}\sum_{(j_1,j_2,\dotsc,j_{k+1})\in\pauu{k+1}{m}}\su{0}^\mpi\su{j_1}\su{0}^\mpi\su{j_2}\dotsm\su{0}^\mpi s_{j_{k+1}}\\
   &=\sum_{\ell=2}^{m}(-1)^\ell\sum_{(j_1,j_2,\dotsc,j_\ell)\in\pauu{\ell}{m}}\su{0}^\mpi\su{j_1}\su{0}^\mpi\su{j_2}\dotsm\su{0}^\mpi s_{j_\ell}.
  \end{split}
 \]
 Taking again into account \rthm{P1414} and additionally \(\pauu{1}{m}=\set{(m)}\), we infer then
 \[
  \begin{split}
   -\sum_{\ell=1}^ms^\rez_{m-\ell}s_\ell s_0^\mpi
   &=\rk*{-\sum_{\ell=1}^{m-1}s^\rez_{m-\ell}s_\ell}s_0^\mpi-s^\rez_0s_m s_0^\mpi\\
   &=\rk*{\sum_{\ell=2}^{m}(-1)^\ell\sum_{(j_1,j_2,\dotsc,j_\ell)\in\pauu{\ell}{m}}\su{0}^\mpi\su{j_1}\su{0}^\mpi\su{j_2}\dotsm\su{0}^\mpi s_{j_\ell}}s_0^\mpi-s^\mpi_0s_m s_0^\mpi\\
   &=\sum_{\ell=2}^{m}(-1)^\ell\sum_{(j_1,j_2,\dotsc,j_\ell)\in\pauu{\ell}{m}}\su{0}^\mpi\su{j_1}\su{0}^\mpi\su{j_2}\dotsm\su{0}^\mpi s_{j_\ell}s_0^\mpi+(-1)^1\sum_{(j_1)\in\pauu{1}{m}}s^\mpi_0s_{j_1}s_0^\mpi\\
   &=\sum_{\ell=1}^m(-1)^\ell\sum_{(j_1,j_2,\dotsc,j_\ell)\in\pauu{\ell}{m}}\su{0}^\mpi\su{j_1}\su{0}^\mpi\su{j_2}\dotsm\su{0}^\mpi\su{j_\ell}\su{0}^\mpi
   =\su{m}^\rez.\qedhere
  \end{split}
 \]
\eproof

 It should be mentioned that the result of \rprop{P0731}  was obtained in~\zitaa{MR3014197}{\cprop{5.24}{36}} for the particular case that \(\seqska\) belongs to \(\Dpqka\).
 The application of~\zitaa{MR3014197}{\cprop{5.24}{36}} lead us in~\zita{MR3014197} to several further results in the mentioned particular case.
 In view of \rthm{P1414} and \rprop{P0731}, we are now able to derive these results in a general case.
 
 Given an arbitrary \(n\in\N\) and arbitrary rectangular complex matrices \(A_1,A_2,\dotsc,A_n\), we write \(\col\seq{A_j}{j}{1}{n}=\col\rk{A_1,A_2,\dotsc,A_n}\) (\tresp{}, \(\row\seq{A_j}{j}{1}{n}\defeq\mat{A_1,A_2,\dotsc,A_n}\)) for the block column (\tresp{}, block row) built from the matrices \(A_1,A_2,\dotsc,A_n\) if their numbers of columns (\tresp{}, rows) are all equal.
 
 Let \(\seqska \) be a sequence of complex \tpqa{matrices}.
 For all \(\ell,m\in\NO\) with \(\ell\leq m\leq\kappa\), then let the block rows \(\yuu{\ell}{m}\) and \(\yhuu{\ell}{m}\) and the block columns \(\zuu{\ell}{m}\) and \(\zhuu{\ell}{m}\) be given by \(\yuu{\ell}{m}\defeq\col\seq{s_j}{j}{\ell}{m}\) and \(\yhuu{\ell}{m}\defeq\col\seq{s_{m-j}}{j}{0}{m-\ell}\) and by \(\zuu{\ell}{m}\defeq\row\seq{s_j}{j}{\ell}{m}\) and \(\zhuu{\ell}{m}\defeq\row\seq{s_{m-j}}{j}{0}{m-\ell}\).
 
 It is readily checked that the identities from \rdefn{D1419} and \rprop{P0731} can be subsumed:
 
\bcorl{ab.R1112}
 If \(\seqska\) is a sequence of complex \tpqa{matrices}, then \(\yuu{1}{m}^\rez=-\SLu{m-1}^\rez\yuu{1}{m}\su{0}^\mpi\) and \(\zuu{1}{m}^\rez=-\su{0}^\mpi\zuu{1}{m}\SUu{m-1}^\rez\) for all \(m\in\mn{1}{\kappa}\).
\ecor
\bproof
 The first identity follows from \rprop{P0731}, whereas the second identity is an immediate consequence of \rdefn{D1419}.
\eproof

 Note that \rcor{ab.R1112} was obtained for \tftd{} sequences in~\zitaa{MR3014197}{\ccor{4.23}{}} as a consequence of~\zitaa{MR3014197}{\cthm{4.21}{}}.

 \rthm{P1414} enables us an alternate approach to~\zitaa{MR3014197}{\cprop{5.16}{34}}:

\bpropl{H.R.rez*}
 If \(\seqska \) is a sequence of complex \tpqa{matrices} with \trF{} \(\seqr{\kappa}\), then \(\seq{r_j^\ad}{j}{0}{\kappa}\) coincides with the \trFa{\(\seq{\su{j}^\ad}{j}{0}{\kappa}\)}.
\eprop
\bproof
 In view of \rrem{A.R.A++*} we have \(\rk{\su{0}^\ad}^\mpi=\rk{\su{0}^\mpi}^\ad\).
 Combining this with the structure of the sets \(\pauu{\ell}{j}\) with \(j\in\N\) and \(\ell\in\mn{1}{j}\) the application of \rthm{P1414} yields the assertion.
\eproof

 Furthermore, in view of \rrem{A.R.l*A}, building the \trF{} is a homogeneous operation of degree~\(-1\) in the following sense:

\breml{H.R.rezhom}
 Let \(\lambda\in\C\) and let \(\seqska \) be a sequence of complex \tpqa{matrices} with \trF{} \(\seqr{\kappa}\).
 Then \(\seq{\lambda^\mpi r_j}{j}{0}{\kappa}\) coincides with the \trFa{\(\seq{\lambda s_j}{j}{0}{\kappa}\)}.
\erem

\bnotal{H.D.CP}
 Let \(\seqska\) be a sequence of complex \tpqa{matrices} and let \(\seqt{\kappa}\) be a sequence of complex \taaa{q}{r}{matrices}.
 As usual, then the notation \(\seq{(s\cp t)_j}{j}{0}{\kappa}\) stands for the \emph{\tCPa{\(\seqska\)}{\(\seqt{\kappa}\)}}, \tie{}, we have \((s\cp t)_j\defeq\sum_{\ell=0}^js_\ell t_{j-\ell}\) for each \(j\in\mn{0}{\kappa}\).
\enota

\breml{H.R.CPtr}
 Let \(\seqska\) be a sequence of complex \tpqa{matrices} and let \(\seqt{\kappa}\) be a sequence of complex \taaa{q}{r}{matrices}.
 Denote by \(\seq{w_j}{j}{0}{\kappa}\) the \tCPa{\(\seqska\)}{\(\seqt{\kappa}\)}.
 For each \(k\in\mn{0}{\kappa}\), then the matrix \(w_k\) is built from the matrices \(s_0,s_1,\dotsc,s_k\) and \(t_0,t_1,\dotsc,t_k\).
 In particular, for each \(m\in\mn{0}{\kappa}\), the \tCPa{\(\seqs{m}\)}{\(\seqt{m}\)} coincides with \(\seq{w_j}{j}{0}{m}\).
\erem

 It is well known and easily seen that the \tCP{} can be performed by multiplying block Toeplitz matrices of the forms described in \rnota{M.N.S}:
 
\breml{M.R.s*t}
 Let \(\seqska \) be a sequence of complex \tpqa{matrices}, let \(\seqt{\kappa}\) be a sequence of complex \taaa{q}{r}{matrices}, and let \(\seq{w_j}{j}{0}{\kappa}\) be a sequence of complex \taaa{p}{r}{matrices}.
 Then the following statements are equivalent:
 \baeqi{0}
  \il{M.R.s*t.i} \(\seq{w_j}{j}{0}{\kappa}\) is the \tCPa{\(\seqska\)}{\(\seqt{\kappa}\)}.
  \il{M.R.s*t.ii} \(\SLuo{m}{s}\SLuo{m}{t}=\SLuo{m}{w}\) for all \(m\in\mn{0}{\kappa}\).
  \il{M.R.s*t.iii} \(\SUuo{m}{s}\SUuo{m}{t}=\SUuo{m}{w}\) for all \(m\in\mn{0}{\kappa}\).
 \eaeqi
\erem

 We investigate now the \tCP{} of a sequence and its associated \trF{} in both possible orders.

\bpropl{H.L.rsst}
 Let \(\seqska\) be a sequence of complex \tpqa{matrices}, let \(\seqa{v}{\kappa}\) be the \tCPa{\(\seq{\su{j}^\rez}{j}{0}{\kappa}\)}{\(\seqska \)}, and let \(\seqa{w}{\kappa}\) be the \tCPa{\(\seqska \)}{\(\seq{\su{j}^\rez}{j}{0}{\kappa}\)}.
 For all \(j,k\in\mn{0}{\kappa}\), then \(v_j\su{k}^\rez=\Kronu{0j}\su{k}^\rez\) and \(\su{k}^\rez w_j=\Kronu{0j}\su{k}^\rez\).
\eprop
\bproof
 We have \(v_0=\su{0}^\mpi\su{0}\) and \(w_0=\su{0}\su{0}^\mpi\).
 Hence, \(v_0\su{0}^\rez=\rk{\su{0}^\mpi\su{0}}\su{0}^\mpi=\su{0}^\mpi=\su{0}^\rez\) and \(\su{0}^\rez w_0=\su{0}^\mpi\rk{\su{0}\su{0}^\mpi}=\su{0}^\mpi=\su{0}^\rez\).
 Now let \(j\in\mn{1}{\kappa}\).
 Because of \rprop{P0731} and \rdefn{D1419}, then
\[\begin{split}
  v_j
  =\sum_{k=0}^j\su{k}^\rez\su{j-k}
  =\sum_{\ell=0}^j\su{j-\ell}^\rez\su{\ell}
  &=\su{j}^\rez\su{0}+\sum_{\ell=1}^j\su{j-\ell}^\rez\su{\ell}\\
  &=\ek*{-\rk*{\sum_{\ell=1}^js^\rez_{j-\ell}s_\ell}s_0^\mpi}\su{0}+\sum_{\ell=1}^j\su{j-\ell}^\rez\su{\ell}
  =\sum_{\ell=1}^j\su{j-\ell}^\rez\su{\ell}(\Iq-\su{0}^\mpi\su{0})
\end{split}\]
 and
\[\begin{split}
  w_j
  =\sum_{k=0}^j\su{k}\su{j-k}^\rez
  =\sum_{\ell=0}^j\su{j-\ell}\su{\ell}^\rez
  &=\sum_{\ell=0}^{j-1}\su{j-\ell}\su{\ell}^\rez+\su{0}\su{j}^\rez\\
  &=\sum_{\ell=0}^{j-1}\su{j-\ell}\su{\ell}^\rez+\su{0}\rk*{-s_0^\mpi\sum_{\ell=0}^{j-1}s_{j-\ell}s_\ell^\rez}
  =(\Ip-\su{0}\su{0}^\mpi)\sum_{\ell=0}^{j-1}\su{j-\ell}\su{\ell}^\rez.
\end{split}\]
 Since \(\su{0}^\mpi\su{0}\su{k}^\rez=\su{k}^\rez\) and \(\su{k}^\rez\su{0}\su{0}^\mpi=\su{k}^\rez\) hold true for all \(k\in\mn{0}{\kappa}\) by \rcor{C1511}, the asserted identities follow.
\eproof

 The following result shows that the block Toeplitz matrices \(\SLu{m}^\rez\) and \(\SUu{m}^\rez\), built via \rnota{M.N.S} from the \trF{}, are in the sense of~\zitaa{MR1987382}{\cdefnp{1}{40}} indeed \tginv{2}s of \(\SLu{m}\) and \(\SUu{m}\), \tresp{:}

\bcorl{H.P.S12}
 If \(\seqska \) is a sequence of complex \tpqa{matrices}, then \(\SLu{m}^\rez\SLu{m}\SLu{m}^\rez=\SLu{m}^\rez\) and \(\SUu{m}^\rez\SUu{m}\SUu{m}^\rez=\SUu{m}^\rez\) for all \(m\in\mn{0}{\kappa}\).
\ecor
\bproof
 Consider an arbitrary \(m\in\mn{0}{\kappa}\).
 Denote by \(\seqa{v}{\kappa}\) the \tCPa{\(\seq{\su{j}^\rez}{j}{0}{\kappa}\)}{\(\seqska \)} and by \(\seqa{g}{\kappa}\) the \tCPa{\(\seqa{v}{\kappa}\)}{\(\seqa{s^\rez}{\kappa}\)}.
 Because of \rrem{M.R.s*t}, then \(\SLuo{m}{g}=\SLuo{m}{v}\SLu{m}^\rez\).
 We have \(v_0=\su{0}^\rez\su{0}=\su{0}^\mpi\su{0}\) and thus \(g_0=v_0\su{0}^\rez=\rk{\su{0}^\mpi\su{0}}\su{0}^\mpi=\su{0}^\mpi=\su{0}^\rez\).
 Using \rprop{H.L.rsst}, we obtain for each \(j\in\mn{1}{\kappa}\) furthermore \(g_j=\sum_{\ell=0}^jv_\ell\su{j-\ell}^\rez=\su{j}^\rez\).
 Therefore, the sequence \(\seqa{s^\rez}{\kappa}\) coincides with the \tCPa{\(\seqa{v}{\kappa}\)}{\(\seqa{s^\rez}{\kappa}\)}.
 From \rrem{M.R.s*t} we infer \(\SLuo{m}{v}=\SLu{m}^\rez\SLu{m}\).
 Consequently, \(\SLu{m}^\rez=\SLuo{m}{g}=\rk{\SLu{m}^\rez\SLu{m}}\SLu{m}^\rez=\SLu{m}^\rez\SLu{m}\SLu{m}^\rez\).
 Analogously, we obtain \(\SUu{m}^\rez=\SUu{m}^\rez\SUu{m}\SUu{m}^\rez\).
\eproof

\bpropnl{\zitaa{MR3014197}{\cprop{4.20}{26}}}{101.S216}
 If \(\seqska\in\Dpqu{\kappa}\), then \(\SLu{m}^\mpi=\SLu{m}^\rez\) and \(\SUu{m}^\mpi=\SUu{m}^\rez\) for all \(m\in\mn{0}{\kappa}\).
\eprop

 Given an arbitrary \(n\in\N\) and arbitrary rectangular complex matrices \(A_1,A_2,\dotsc,A_n\), we use \(\diag\seq{A_j}{j}{1}{n}\) or \(\diag\rk{A_1,A_2,\dotsc,A_n}\) to denote the corresponding block diagonal matrix.
 Furthermore, for arbitrarily given \(A\in\Cpq\) and \(m\in\NO\) we write
\beql{sdiag}
 \sdiag{A}{m}
 \defeq\diag\seq{A}{j}{0}{m}.
\eeq
 
\blemnl{\tcf{}~\zitaa{MR3014199}{\clem{4.15}{134}}}{103.M39}
 If \(\seqska \in\Dpqu{\kappa}\), then \(\SLu{m}\SLu{m}^\mpi=\sdiag{\su{0}\su{0}^\mpi}{m}\) and \(\SUu{m}\SUu{m}^\mpi=\sdiag{\su{0}\su{0}^\mpi}{m}\) and furthermore \(\SLu{m}^\mpi\SLu{m}=\sdiag{\su{0}^\mpi\su{0}}{m}\) and \(\SUu{m}^\mpi\SUu{m}=\sdiag{\su{0}^\mpi\su{0}}{m}\) for all \(m\in\mn{0}{\kappa}\).
\elem

 In view of \rrem{ab.R1052}, we infer from the first chain of equations in \rlem{103.M39} that the column spaces of the matrices \(\SLu{m}\) and \(\SUu{m}\) built from a sequence belonging to \(\Dpqu{\kappa}\) are exactly the \((m+1)\)\nobreakdash-fold direct product of the column space \(\ran{\su{0}}\), which implies block diagonal form of the transformation matrix associated to the corresponding orthogonal projection with respect to the standard basis.
 
 A closer look at~\zitaa{MR3014199}{\cdefn{4.16}{}} shows that several constructions done there even work for a slightly larger class then the set of all \tftd{} sequences.
 Let \(m\in\N\).
 A sequence \(\seqs{m}\) of complex \tpqa{matrices} for which \(\seqs{m-1}\in\Dpqu{m-1}\) is called a
 \emph{\tnftd{}} sequence.
 The set of all \tnftd{} sequences \(\seqs{m}\) will be denoted by \(\Dtpqu{m}\).
 We also set \(\Dtpqu{0}\defeq\Dpqu{0}\).
 Obviously, we have \(\Dpqu{m}\subseteq\Dtpqu{m}\) for all \(m\in\NO\).
 
\bleml{H.L.HsL}
 Let \(n\in\N\) and let \(\seqs{2n}\in\Dtpqu{2n}\).
 Then
 \beql{H.L.HsL.A}
 \Hu{n}
 =
 \bMat
  \Ip & \Ouu{p}{np}\\
  \yuu{1}{n}\su{0}^\mpi & \Iu{np}
 \eMat\rk{\zdiag{\su{0}}{\LLu{n}}}
 \bMat
  \Iq & \su{0}^\mpi \zuu{1}{n} \\
  \Ouu{nq}{q}&\Iu{nq}
 \eMat.
\eeq
\elem
\bproof
 From the assumption \(\seqs{2n}\in\Dtpqu{2n}\), we get \(\ran{\zuu{1}{n}}\subseteq\ran{\su{0}}\) and \(\nul{\su{0}}\subseteq\nul{\yuu{1}{n}}\).
 Thus, \rrem{R.AA+B} yields \(\su{0}\su{0}^\mpi\zuu{1}{n}=\zuu{1}{n}\) and \(\yuu{1}{n}\su{0}^\mpi\su{0}=\yuu{1}{n}\).
 Taking additionally into account \rremss{H.R.LL}{H.R.Hblock}, then \eqref{H.L.HsL.A} follows by direct calculation.
\eproof

 \rprop{P0731} shows that the assumption \(\seqs{m}\in\Dtpqu{m}\) in~\zitaa{MR3014197}{\cprop{5.24}{36}} can be dropped.
 Similarly, some of the main results in~\zitaa{MR3014199}{\csec{6}} can be proved without assuming that \(\seqs{m}\) belongs to \(\Dtpqu{m}\).

\section{Some identities for \tbHms{} formed by a sequence and its reciprocal}\label{S1216}

 The main goal of this section is the investigation of the interplay between various \tbHms{}.
 This topic was already studied in~\zitaa{MR3014199}{\csec{6}} for \tnftd{} sequences.
 The application of results from \rsec{S0858} will enable us to verify that central results obtained in~\zitaa{MR3014199}{\csec{6}} are even true in a more general case.
 The first step of our strategy to prove this coincides with that one in~\zita{MR3014199}.
 More precisely, we apply several identities for block Hankel structures built from the \tCP{} of sequences of complex matrices introduced in \rnota{H.D.CP}.
 We work with the following three types of \tbHms{} associated with a sequence of complex \tpqa{matrices}.
 
\bnotal{N.HKG}
 Let \(\seqska \) be a sequence of complex \tpqa{matrices}.
 Then let the \tbHms{} \(\Huo{n}{s}\), \(\Kuo{n}{s}\), and \(\Guo{n}{s}\) be given by \(\Huo{n}{s}\defeq\matauuuo{s_{j+k}}{j}{k}{0}{n}\) for all \(n\in\NO\) with \(2n\leq\kappa\), by \(\Kuo{n}{s}\defeq\matauuuo{s_{j+k+1}}{j}{k}{0}{n}\) for all \(n\in\NO\) with \(2n+1\leq\kappa\), and by \(\Guo{n}{s}\defeq\matauuuo{s_{j+k+2}}{j}{k}{0}{n}\) for all \(n\in\NO\) with \(2n+2\leq\kappa\), \tresp{}
\enota

 If it is clear which sequence \(\seqska\) is meant, then we write \(\Hu{n}\), \(\Ku{n}\) and \(\Gu{n}\) instead of  \(\Huo{n}{s}\), \(\Kuo{n}{s}\) and \(\Guo{n}{s}\), respectively. 
 We write
\(%
 \zdiag{A}{B}
 \defeq\diag\rk{A,B}
\), for arbitrarily given two matrices \(A\in\Cpq\) and \(B\in\Coo{r}{s}\).
 Note that in the following identities, block Toeplitz structures also occur in form of the block triangular matrices given in \rnota{M.N.S}.

\bpropnl{\tcf{}~\zitaa{MR3014199}{\cpropsss{5.2}{138}{5.4}{139}{5.5}{139}}}{P1522}
 Let \(\seqska\) be a sequence of complex \tpqa{matrices} and let \(\seqt{\kappa}\) be a sequence of complex \taaa{q}{r}{matrices}.
 Denote by \(\seq{w_j}{j}{0}{\kappa}\) the \tCPa{\(\seqska\)}{\(\seqt{\kappa}\)}.
 Then
\begin{align*}
 \Huo{n}{w}&=\Huo{n}{s}\SUuo{n}{t}+\rk{\zdiag{\Opq}{\SLuo{n-1}{s}}}\Huo{n}{t}&
&\text{and}&
 \Huo{n}{w}&=\Huo{n}{s}\rk{\zdiag{\Opq}{\SUuo{n-1}{t}}}+\SLuo{n}{s}\Huo{n}{t}
\end{align*}
 for all \(n\in\N\) with \(2n\leq\kappa\),
\[
 \Kuo{n}{w}
 =\Kuo{n}{s}\SUuo{n}{t}+\SLuo{n}{s}\Kuo{n}{t}
\]
 for all \(n\in\NO\) with \(2n+1\leq\kappa\), and
\[
 \Guo{n}{w}
 =\Guo{n}{s}\SUuo{n}{t}+\yuuo{1}{n+1}{s}\zuuo{1}{n+1}{t}+\SLuo{n}{s}\Guo{n}{t}
\]
 for all \(n\in\NO\) with \(2n+2\leq\kappa\).
\eprop

 Now we consider the interplay between \tbHms{} generated by a sequence \(\seqska\) on the one side and by the \trF{} \(\seqa{s^\rez}{\kappa}\) on the other side.
 It will turn out that certain identities, crucial for our subsequent considerations, can be obtained by applying \rprop{P1522} to a sequence \(\seqska\) of complex \tpqa{matrices} together with its \trF{} \(\seq{\su{j}^\rez}{j}{0}{\kappa}\) introduced in \rdefn{D1419}.
 
\breml{H.R.Sblock}
 If \(\seqska\) is a sequence of complex \tpqa{matrices}, then \(\SLu{m}=\tmat{\su{0} & \Ouu{p}{mq} \\\yuu{1}{m} & \SLu{m-1}}\), \(\SLu{m}=\tmat{\SLu{m-1} & \Ouu{mp}{q} \\\zhuu{1}{m} &\su{0} }\), \(\SUu{m}=\tmat{\su{0} & \zuu{1}{m}\\ \Ouu{mp}{q}  & \SUu{m-1}}\) and \(\SUu{m}=\tmat{\SUu{m-1} & \yhuu{1}{m}\\ \Ouu{p}{mq}  & \su{0}}\) for all \(m\in\mn{1}{\kappa}\).
\erem

\bleml{L1506}
 Let \(\seqska\) be a sequence of complex \tpqa{matrices}.
 Denote by \(\seqa{v}{\kappa}\) the \tCPa{\(\seq{s^\rez_j}{j}{0}{\kappa}\)}{\(\seqska \)} and by \(\seqa{w}{\kappa}\) the \tCPa{\(\seqska \)}{\(\seq{s^\rez_j}{j}{0}{\kappa}\)}.
 Then:
\benui
 \il{L1506.a} Let \(j,k\in\NO\) with \(j+k\leq\kappa\).
 Then \(\zhuuo{0}{j}{v}\yuu{k}{k+j}^\rez=\su{k+j}^\rez\) and \(\zhuu{0}{j}^\rez\yuuo{k}{k+j}{w}=\Kronu{0k}\su{j}^\rez\).
 \il{L1506.b} Let \(j,k\in\NO\) and let \(\ell\in\N\) with \(j+k+\ell\leq\kappa\).
 Then \(\mat{\zhuuo{0}{j}{v},\Ouu{q}{\ell q}}\yuu{k}{k+j+\ell}^\rez=\su{k+j}^\rez\) and \(\mat{\zhuu{0}{j}^\rez,\Ouu{q}{\ell p}}\yuuo{k}{k+j+\ell}{w}=\Kronu{0k}\su{j}^\rez\).
 \il{L1506.c} Let \(m,k\in\NO\) with \(m+k\leq\kappa\).
 Then \(\SLuo{m}{v}\yuu{k}{k+m}^\rez=\yuu{k}{k+m}^\rez\) and \(\SLu{m}^\rez\yuuo{k}{k+m}{w}=\Kronu{0k}\yuu{0}{m}^\rez\).
\eenui
\elem
\bproof
 \eqref{L1506.a} Using \rprop{H.L.rsst}, we obtain
\[
 \zhuuo{0}{j}{v}\yuu{k}{k+j}^\rez
 =\sum_{\ell=0}^jv_{j-\ell}\su{k+\ell}^\rez
 =\sum_{\ell=0}^j\Kronu{0,j-\ell}\su{k+\ell}^\rez
 =\su{k+j}^\rez
\]
 and, in view of \(k\geq0\), furthermore
\[
 \zhuu{0}{j}^\rez\yuuo{k}{k+j}{w}
 =\sum_{\ell=0}^j\su{j-\ell}^\rez w_{k+\ell}
 =\sum_{\ell=0}^j\Kronu{0,k+\ell}\su{j-\ell}^\rez
 =
 \begin{cases}
  \su{j}^\rez\tincase{k=0}\\
  \Oqp\tincase{k\geq1}
 \end{cases}.
\]

 \eqref{L1506.b} Using~\eqref{L1506.a}, we get
\[
 \mat{\zhuuo{0}{j}{v},\Ouu{q}{\ell q}}\yuu{k}{k+j+\ell}^\rez
 =\mat{\zhuuo{0}{j}{v},\Ouu{q}{\ell q}}\matp{\yuu{k}{k+j}^\rez}{\yuu{k+j+1}{k+j+\ell}^\rez}
 =\zhuuo{0}{j}{v}\yuu{k}{k+j}^\rez
 =\su{k+j}^\rez
\]
 and
\[
 \mat{\zhuu{0}{j}^\rez,\Ouu{q}{\ell p}}\yuuo{k}{k+j+\ell}{w}
 =\mat{\zhuu{0}{j}^\rez,\Ouu{q}{\ell p}}\matp{\yuuo{k}{k+j}{w}}{\yuuo{k+j+1}{k+j+\ell}{w}}
 =\zhuu{0}{j}^\rez\yuuo{k}{k+j}{w}
 =\Kronu{0k}\su{j}^\rez.
\]
 
 \eqref{L1506.c} Using~\eqref{L1506.b} and~\eqref{L1506.a}, we obtain
\[
 \SLuo{m}{v}\yuu{k}{k+m}^\rez
 =
 \bMat
  \mat{\zhuuo{0}{0}{v},\Ouu{q}{m q}}\\
  \mat{\zhuuo{0}{1}{v},\Ouu{q}{\rk{m-1}q}}\\
  \vdots\\
  \mat{\zhuuo{0}{m-1}{v},\Ouu{q}{q}}\\
  \zhuuo{0}{m}{v}
 \eMat\yuu{k}{k+m}^\rez
 =
 \bMat
  \mat{\zhuuo{0}{0}{v},\Ouu{q}{m q}}\yuu{k}{k+m}^\rez\\
  \mat{\zhuuo{0}{1}{v},\Ouu{q}{\rk{m-1}q}}\yuu{k}{k+m}^\rez\\
  \vdots\\
  \mat{\zhuuo{0}{m-1}{v},\Ouu{q}{q}}\yuu{k}{k+m}^\rez\\
  \zhuuo{0}{m}{v}\yuu{k}{k+m}^\rez
 \eMat
 =
 \bMat
  \su{k+0}^\rez\\
  \su{k+1}^\rez\\
  \vdots\\
  \su{k+m-1}^\rez\\
  \su{k+m}^\rez
 \eMat
 =\yuu{k}{k+m}^\rez 
\]
 and
\[
 \SLu{m}^\rez\yuuo{k}{k+m}{w}
 =
 \bMat
  \mat{\zhuu{0}{0}^\rez,\Ouu{q}{m p}}\\
  \mat{\zhuu{0}{1}^\rez,\Ouu{q}{\rk{m-1}p}}\\
  \vdots\\
  \mat{\zhuu{0}{m-1}^\rez,\Ouu{q}{p}}\\
  \zhuu{0}{m}^\rez
 \eMat\yuuo{k}{k+m}{w}
 =
 \bMat
  \mat{\zhuu{0}{0}^\rez,\Ouu{q}{m p}}\yuuo{k}{k+m}{w}\\
  \mat{\zhuu{0}{1}^\rez,\Ouu{q}{\rk{m-1}p}}\yuuo{k}{k+m}{w}\\
  \vdots\\
  \mat{\zhuu{0}{m-1}^\rez,\Ouu{q}{p}}\yuuo{k}{k+m}{w}\\
  \zhuu{0}{m}^\rez\yuuo{k}{k+m}{w}
 \eMat
 =
 \bMat
  \Kronu{0k}\su{0}^\rez\\
  \Kronu{0k}\su{1}^\rez\\
  \vdots\\
  \Kronu{0k}\su{m-1}^\rez\\
  \Kronu{0k}\su{m}^\rez
 \eMat
 =\Kronu{0k}\yuu{0}{m}^\rez.\qedhere
\]
\eproof

\bnotal{H.N.updo}
 For each \(m\in\N\), let \(\IOquu{0}{m}\defeq\Iu{qm}\) and let \(\OIquu{0}{m}\defeq\Iu{qm}\).
 Furthermore, for all \(\ell,m\in\N\), let \(\IOquu{\ell}{m}\defeq\tmatp{\Iu{mq}}{\Ouu{\ell q}{mq}}\) and let \(\OIquu{\ell}{m}\defeq\tmatp{\Ouu{\ell q}{mq}}{\Iu{mq}}\).
\enota


\bleml{L1419}
 Let \(\seqska\) be a sequence of complex \tpqa{matrices}.
 Denote by \(\seqa{v}{\kappa}\) the \tCPa{\(\seq{s^\rez_j}{j}{0}{\kappa}\)}{\(\seqska \)} and by \(\seqa{w}{\kappa}\) the \tCPa{\(\seqska \)}{\(\seq{s^\rez_j}{j}{0}{\kappa}\)}.
 Then:
\benui
 \il{L1419.a} \(\SLuo{n}{v}\Ku{n}^\rez=\Ku{n}^\rez\) and \(\SLu{n}^\rez\Kuo{n}{w}=\Ouu{(n+1)q}{(n+1)p}\) for all \(n\in\NO\) with \(2n+1\leq\kappa\).
 \il{L1419.b} \(\Hu{n}^\rez-\rk{\zdiag{\Oqq}{\SLuo{n-1}{v}}}\Hu{n}^\rez=\vqu{n}\zuu{0}{n}^\rez\) and \(\SLu{n}^\rez\Huo{n}{w}=\yuu{0}{n}^\rez\vpu{n}^\ad\) for all \(n\in\N\) with \(2n\leq\kappa\).
 \il{L1419.c} \(\SLuo{n}{v}\Gu{n}^\rez=\Gu{n}^\rez\), \(\SLu{n}^\rez\Guo{n}{w}=\Ouu{(n+1)q}{(n+1)p}\), and \(\su{0}^\rez\zuuo{1}{n+1}{w}=\Ouu{q}{(n+1)p}\) for all \(n\in\NO\) with \(2n+2\leq\kappa\).
\eenui
\elem
\bproof
 \eqref{L1419.a} Consider an arbitrary \(n\in\NO\) with \(2n+1\leq\kappa\).
 Using \rlemp{L1506}{L1506.c}, we obtain
\[\begin{split}
 \SLuo{n}{v}\Ku{n}^\rez
 &=\SLuo{n}{v}\mat{\yuu{1}{n+1}^\rez,\yuu{2}{n+2}^\rez,\dotsc,\yuu{n+1}{2n+1}^\rez}\\
 &=\mat{\SLuo{n}{v}\yuu{1}{n+1}^\rez,\SLuo{n}{v}\yuu{2}{n+2}^\rez,\dotsc,\SLuo{n}{v}\yuu{n+1}{2n+1}^\rez}\\
 &=\mat{\yuu{1}{n+1}^\rez,\yuu{2}{n+2}^\rez,\dotsc,\yuu{n+1}{2n+1}^\rez} 
 =\Ku{n}^\rez
\end{split}\]
 and
\[\begin{split}
 \SLu{n}^\rez\Kuo{n}{w}
 &=\SLu{n}^\rez\mat{\yuuo{1}{n+1}{w},\yuuo{2}{n+2}{w},\dotsc,\yuuo{n+1}{2n+1}{w}}\\
 &=\mat{\SLu{n}^\rez\yuuo{1}{n+1}{w},\SLu{n}^\rez\yuuo{2}{n+2}{w},\dotsc,\SLu{n}^\rez\yuuo{n+1}{2n+1}{w}}\\
 &=\mat{\Kronu{01}\yuu{0}{n}^\rez,\Kronu{02}\yuu{0}{n}^\rez,\dotsc,\Kronu{0,n+1}\yuu{0}{n}^\rez} 
 =\Ouu{(n+1)q}{(n+1)p}.
\end{split}\]

 \eqref{L1419.b} Consider an arbitrary \(n\in\N\) with \(2n\leq\kappa\).
 Using~\eqref{L1419.a} and \rlemp{L1506}{L1506.c}, we obtain
\[\begin{split}
 \rk{\zdiag{\Oqq}{\SLuo{n-1}{v}}}\Hu{n}^\rez
 &=
 \bMat
  \Oqq&\Ouu{q}{nq}\\
  \Ouu{nq}{q}&\SLuo{n-1}{v}
 \eMat
 \bMat
  \zuu{0}{n-1}^\rez&\su{n}^\rez\\
  \Ku{n-1}^\rez&\yuu{n+1}{2n}^\rez
 \eMat\\
 &=
 \bMat
  \Ouu{q}{np}&\Oqp\\
  \SLuo{n-1}{v}\Ku{n-1}^\rez&\SLuo{n-1}{v}\yuu{n+1}{2n}^\rez
 \eMat
 =
 \bMat
  \Ouu{q}{np}&\Oqp\\
  \Ku{n-1}^\rez&\yuu{n+1}{2n}^\rez
 \eMat.
\end{split}\]
 Hence, we can conclude
\[\begin{split}
 \Hu{n}^\rez-\rk{\zdiag{\Oqq}{\SLuo{n-1}{v}}}\Hu{n}^\rez
 &=
 \bMat
  \zuu{0}{n-1}^\rez&\su{n}^\rez\\
  \Ku{n-1}^\rez&\yuu{n+1}{2n}^\rez
 \eMat-
 \bMat
  \Ouu{q}{np}&\Oqp\\
  \Ku{n-1}^\rez&\yuu{n+1}{2n}^\rez
 \eMat\\
 &=
 \bMat
  \zuu{0}{n-1}^\rez&\su{n}^\rez\\
  \Ouu{nq}{np}&\Ouu{nq}{p}
 \eMat
 =
 \bMat
  \zuu{0}{n}^\rez\\
  \Ouu{nq}{\rk{n+1}p}
 \eMat
 =\vqu{n}\zuu{0}{n}^\rez.
\end{split}\]
 From \rlemp{L1506}{L1506.c} and \rrem{H.R.Sblock} we get furthermore
\[\begin{split}
 \SLu{n}^\rez\Huo{n}{w}
 &=\SLu{n}^\rez\mat{\yuuo{0}{n}{w},\yuuo{1}{n+1}{w},\dotsc,\yuuo{n}{2n}{w}}
 =\mat{\SLu{n}^\rez\yuuo{0}{n}{w},\SLu{n}^\rez\yuuo{1}{n+1}{w},\dotsc,\SLu{n}^\rez\yuuo{n}{2n}{w}}\\
 &=\mat{\Kronu{00}\yuu{0}{n}^\rez,\Kronu{01}\yuu{0}{n}^\rez,\dotsc,\Kronu{0n}\yuu{0}{n}^\rez}
 =\mat{\yuu{0}{n}^\rez,\Ouu{\rk{n+1}q}{np}}
 =\yuu{0}{n}^\rez\vpu{n}^\ad.
\end{split}\]

 \eqref{L1419.c} Consider an arbitrary \(n\in\NO\) with \(2n+2\leq\kappa\).
 Using \rlemp{L1506}{L1506.c}, we obtain
\[\begin{split}
 \SLuo{n}{v}\Gu{n}^\rez
 &=\SLuo{n}{v}\mat{\yuu{2}{n+2}^\rez,\yuu{3}{n+3}^\rez,\dotsc,\yuu{n+2}{2n+2}^\rez}\\
 &=\mat{\SLuo{n}{v}\yuu{2}{n+2}^\rez,\SLuo{n}{v}\yuu{3}{n+3}^\rez,\dotsc,\SLuo{n}{v}\yuu{n+2}{2n+2}^\rez}\\
 &=\mat{\yuu{2}{n+2}^\rez,\yuu{3}{n+3}^\rez,\dotsc,\yuu{n+2}{2n+2}^\rez} 
 =\Gu{n}^\rez
\end{split}\]
 and
\[\begin{split}
 \SLu{n}^\rez\Guo{n}{w}
 &=\SLu{n}^\rez\mat{\yuuo{2}{n+2}{w},\yuuo{3}{n+3}{w},\dotsc,\yuuo{n+2}{2n+2}{w}}\\
 &=\mat{\SLu{n}^\rez\yuuo{2}{n+2}{w},\SLu{n}^\rez\yuuo{3}{n+3}{w},\dotsc,\SLu{n}^\rez\yuuo{n+2}{2n+2}{w}}\\
 &=\mat{\Kronu{02}\yuu{0}{n}^\rez,\Kronu{03}\yuu{0}{n}^\rez,\dotsc,\Kronu{0,n+2}\yuu{0}{n}^\rez} 
 =\Ouu{(n+1)q}{(n+1)p}.
\end{split}\]
 Furthermore, \rprop{H.L.rsst} yields
\[\begin{split}
 \su{0}^\rez\zuuo{1}{n+1}{w}
 =\su{0}^\rez\mat{w_1,w_2,\dotsc,w_{n+1}}
 &=\mat{\su{0}^\rez w_1,\su{0}^\rez w_2,\dotsc,\su{0}^\rez w_{n+1}}\\
 &=\mat{\Kronu{01}\su{0}^\rez,\Kronu{02}\su{0}^\rez,\dotsc,\Kronu{0,n+1}\su{0}^\rez}
 =\Ouu{q}{(n+1)p}.\qedhere
\end{split}\]
\eproof

\bthml{H.T1431}
 Let \(\seqska\) be a sequence of complex \tpqa{matrices}.
 Then using \rnotass{M.N.S}{H.N.updo} the identity \(\Hu{n}^\rez + \SLu{n}^\rez \Hu{n}\SUu{n}^\rez= \yuu{0}{n}^\rez \vpu{n}^\ad + \vqu{n} \zuu{0}{n}^\rez\) holds for all \(n\in\NO\) with \(2n\leq\kappa\).
\ethm
\bproof
 By virtue of \(\su{0}^\rez=\su{0}^\mpi\), the asserted equation can be obtained for \(n=0\) by direct calculation.
 Now consider an arbitrary \(n\in\N\) with \(2n\leq\kappa\).
 Denote by \(\seqa{v}{\kappa}\) the \tCPa{\(\seq{s^\rez_j}{j}{0}{\kappa}\)}{\(\seqska \)} and by \(\seqa{w}{\kappa}\) the \tCPa{\(\seqska \)}{\(\seq{s^\rez_j}{j}{0}{\kappa}\)}.
 Because of the first equation in \rprop{P1522}, then
\(
 \Hu{n}\SUu{n}^\rez
 =\Huo{n}{w}-\rk{\zdiag{\Opq}{\SLu{n-1}}}\Hu{n}^\rez
\).
 According to \rrem{M.R.s*t}, furthermore \(\SLuo{n-1}{v}=\SLu{n-1}^\rez\SLu{n-1}\).
 Using \rrem{H.R.Sblock}, we obtain
\(
 \SLu{n}^\rez\rk{\zdiag{\Oqq}{\SLu{n-1}}}
 =\zdiag{\Oqq}{\SLuo{n-1}{v}}
\).
 In view of \rlemp{L1419}{L1419.b}, we have
 \(\Hu{n}^\rez-\rk{\zdiag{\Oqq}{\SLuo{n-1}{v}}}\Hu{n}^\rez=\vqu{n}\zuu{0}{n}^\rez\) and \(\SLu{n}^\rez\Huo{n}{w}=\yuu{0}{n}^\rez\vpu{n}^\ad\).
 Hence, we get
\[\begin{split}
 \Hu{n}^\rez+\SLu{n}^\rez\Hu{n}\SUu{n}^\rez
 &=\Hu{n}^\rez+\SLu{n}^\rez\ek*{\Huo{n}{w}-\rk{\zdiag{\Opq}{\SLu{n-1}}}\Hu{n}^\rez}\\
 &=\Hu{n}^\rez+\SLu{n}^\rez\Huo{n}{w}-\SLu{n}^\rez\rk{\zdiag{\Opq}{\SLu{n-1}}}\Hu{n}^\rez\\
 &=\SLu{n}^\rez\Huo{n}{w}+\Hu{n}^\rez-\rk{\zdiag{\Oqq}{\SLuo{n-1}{v}}}\Hu{n}^\rez
 =\yuu{0}{n}^\rez\vpu{n}^\ad+\vqu{n}\zuu{0}{n}^\rez.\qedhere
\end{split}\]
\eproof

\bthml{H.P0836}
 Let \(\seqska\) be a sequence of complex \tpqa{matrices}.
 Then \(\Ku{n}^\rez=-\SLu{n}^\rez\Ku{n}\SUu{n}^\rez\) for all \(n\in\NO\) with \(2n+1\leq\kappa\).
\ethm
\bproof
 Consider an arbitrary \(n\in\NO\) with \(2n+1\leq\kappa\).
 Denote by \(\seqa{v}{\kappa}\) the \tCPa{\(\seq{s^\rez_j}{j}{0}{\kappa}\)}{\(\seqska \)} and by \(\seqa{w}{\kappa}\) the \tCPa{\(\seqska \)}{\(\seq{s^\rez_j}{j}{0}{\kappa}\)}.
 Because of \rprop{P1522}, we have
 \(
  \Ku{n}\SUu{n}^\rez
  =\Kuo{n}{w}-\SLu{n}\Ku{n}^\rez,
 \)
 whereas \rrem{M.R.s*t} yields
 \(
  \SLuo{n}{v}=\SLu{n}^\rez\SLu{n}
 \).
 In view of \rlemp{L1419}{L1419.a}, we have
\(
 \SLuo{n}{v}\Ku{n}^\rez=\Ku{n}^\rez
\) and \(
 \SLu{n}^\rez\Kuo{n}{w}=\Ouu{(n+1)q}{(n+1)p}
\).
 Consequently,
\[
  -\SLu{n}^\rez\Ku{n}\SUu{n}^\rez
  =-\SLu{n}^\rez\rk{\Kuo{n}{w}-\SLu{n}\Ku{n}^\rez}
  =\SLu{n}^\rez\SLu{n}\Ku{n}^\rez-\SLu{n}^\rez\Kuo{n}{w}
  =\SLuo{n}{v}\Ku{n}^\rez-\SLu{n}^\rez\Kuo{n}{w}
  =\Ku{n}^\rez.\qedhere
\]
\eproof

 In our following consideration, Schur complements in \tbHms{} play an essential role.
 For this reason, we recall this construction:
 If \(M=\tmat{A & B\\ C & D}\) is the \tbr{} of a complex \taaa{(p+q)}{(r+s)}{matrix} \(M\) with \taaa{p}{r}{block} \(A\), then the matrix
\beql{E/A}
 M\schca A
 \defeq D-CA^\mpi B
\eeq
 is called the \emph{Schur complement of \(A\) in \(M\)}.
 Concerning a variety of applications of this concept in a lot of areas of mathematics, we refer to~\zitas{MR2160825}.
 The \tbHm{} \(\Hu{n}\) admits the following \tbr{s}:

\breml{H.R.Hblock}
 If \(\seqska\) is a sequence of complex \tpqa{matrices}, then \(\Hu{n}=\tmat{\Hu{n-1} & \yuu{n}{2n-1} \\\zuu{n}{2n-1} & s_{2n}}\) and \(\Hu{n}=\tmat{s_0&\zuu{1}{n}  \\\yuu{1}{n} & \Gu{n-1}}\) for all \(n\in\N\) with \(2n\leq\kappa\).
\erem
  
\bnotal{H.N.LL}
 If \(\seqska\) is a sequence of complex \tpqa{matrices}, then, by virtue of \rrem{H.R.Hblock}, let \(\LLu{0}\defeq\Hu{0}\) and let \(\LLu{n}\defeq\Hu{n}\schca\su{0}\) for all \(n\in\N\) with \(2n\leq\kappa\).
\enota

\breml{H.R.LL}
 If \(\seqska\) is a sequence of complex \tpqa{matrices}, then from \rrem{H.R.Hblock} we infer that \(\LLu{n}=\Gu{n-1}-\yuu{1}{n}\su{0}^\mpi\zuu{1}{n}\) for all \(n\in\N\) with \(2n\leq\kappa\).
\erem

\bthml{H.P0938}
 Let \(\seqska\) be a sequence of complex \tpqa{matrices}.
 Using \rnota{H.N.LL}, then
 \(
  \Gu{n}^\rez
  =-\SLu{n}^\rez\LLu{n+1}\SUu{n}^\rez
 \)
 for all \(n\in\NO\) with \(2n+2\leq\kappa\).
\ethm
\bproof
 Consider an arbitrary \(n\in\NO\) with \(2n+2\leq\kappa\).
 Denote by \(\seqa{v}{\kappa}\) the \tCPa{\(\seq{\su{j}^\rez}{j}{0}{\kappa}\)}{\(\seqska \)} and by \(\seqa{w}{\kappa}\) the \tCPa{\(\seqska\)}{\(\seq{\su{j}^\rez}{j}{0}{\kappa}\)}.
 Because of \rprop{P1522}, then
\beql{H.P0938.1}
 \Gu{n}\SUu{n}^\rez
 =\Guo{n}{w}-\yuu{1}{n+1}\zuu{1}{n+1}^\rez-\SLu{n}\Gu{n}^\rez.
\eeq
 In view of \rprop{P1522}, we have \(\Hu{n+1}\SUu{n+1}^\rez=\Huo{n+1}{w}-\rk{\zdiag{\Opq}{\SLu{n}}}\Hu{n+1}^\rez\).
 Comparing the upper right \taaa{p}{(n+1)p}{block} in the latter equation, in view of \rremss{H.R.Hblock}{H.R.Sblock}, we obtain then \(\su{0}\zuu{1}{n+1}^\rez+\zuu{1}{n+1}\SUu{n}^\rez=\zuuo{1}{n+1}{w}\).
 Since we have by construction \(\su{0}^\rez=\su{0}^\mpi\) and by \rcor{C1511} moreover \(\su{0}^\mpi\su{0}\su{k}^\rez=\su{k}^\rez\) for all \(k\in\mn{0}{\kappa}\), hence,
\beql{H.P0938.2}
 \su{0}^\rez\zuuo{1}{n+1}{w}
 =\su{0}^\mpi\zuuo{1}{n+1}{w}
 =\zuu{1}{n+1}^\rez+\su{0}^\mpi\zuu{1}{n+1}\SUu{n}^\rez
\eeq
 follows.
 According to \rrem{M.R.s*t}, we have \(\SLuo{n}{v}=\SLu{n}^\rez\SLu{n}\).
 In view of \rlemp{L1419}{L1419.c}, we have
\begin{align}\label{H.P0938.4}
 \SLuo{n}{v}\Gu{n}^\rez&=\Gu{n}^\rez,&
 \su{0}^\rez\zuuo{1}{n+1}{w}&=\Ouu{q}{(n+1)p},&
&\text{and}&
 \SLu{n}^\rez\Guo{n}{w}&=\Ouu{(n+1)q}{(n+1)p}.
\end{align}
 By application of \rrem{H.R.LL}, \eqref{H.P0938.1}, \eqref{H.P0938.2}, \(\SLuo{n}{v}=\SLu{n}^\rez\SLu{n}\), and \eqref{H.P0938.4}, we obtain thus
\[\begin{split}
 -\SLu{n}^\rez\LLu{n+1}\SUu{n}^\rez
 &=-\SLu{n}^\rez\rk{\Gu{n}-\yuu{1}{n+1}\su{0}^\mpi\zuu{1}{n+1}}\SUu{n}^\rez
 =-\SLu{n}^\rez\rk{\Gu{n}\SUu{n}^\rez-\yuu{1}{n+1}\su{0}^\mpi\zuu{1}{n+1}\SUu{n}^\rez}\\
 &=-\SLu{n}^\rez\rk{\Guo{n}{w}-\yuu{1}{n+1}\zuu{1}{n+1}^\rez-\SLu{n}\Gu{n}^\rez-\yuu{1}{n+1}\su{0}^\mpi\zuu{1}{n+1}\SUu{n}^\rez}\\
 &=\SLu{n}^\rez\SLu{n}\Gu{n}^\rez+\SLu{n}^\rez\yuu{1}{n+1}\rk{\zuu{1}{n+1}^\rez+\su{0}^\mpi\zuu{1}{n+1}\SUu{n}^\rez}-\SLu{n}^\rez\Guo{n}{w}\\
 &=\SLuo{n}{v}\Gu{n}^\rez+\SLu{n}^\rez\yuu{1}{n+1}\su{0}^\rez\zuuo{1}{n+1}{w}-\SLu{n}^\rez\Guo{n}{w}
 =\Gu{n}^\rez.\qedhere
\end{split}\]
\eproof

\section{Matricial Hamburger moment sequences and \hHT{ation}}\label{H.S}
 We recall classes of sequences of complex \tqqa{matrices} corresponding to solvability criteria for matricial moment problems on \(\Omega=\R\): 
 For each \(n\in\NO\), denote by \(\Hggqu{2n}\) the set of all sequences \(\seqs{2n}\) of complex \tqqa{matrices} for which the corresponding \tbHm{} \(\Hu{n}=\matauuo{\su{j+k}}{j,k}{0}{n}\)
 is \tnnH{}.
 Furthermore, denote by \(\Hggqinf\) the set of all sequences \(\seqsinf\) of complex \tqqa{matrices} satisfying \(\seqs{2n}\in\Hggqu{2n}\) for all \(n\in\NO\).
 The sequences belonging to \(\Hggqu{2n}\) or \(\Hggqinf\) are said to be \emph{\tHnnd}.
 (Note that in~\zita{MR2570113} and subsequent papers of the corresponding authors, the sequences belonging to \(\Hggqu{2\kappa}\) were called \emph{Hankel non-negative definite}.
 Our terminology here differs for the sake of consistency.)
 Finite \tHnnd{} sequences constitute exactly the class of prescribed data \(\seqs{2n}\) for which a slightly modified version of the Problem~\mprob{\R}{2n}{=} is solvable (see, \teg{}~\zitas{MR1624548,MR2570113,MR1395706}).

 We continue with some later used observations on the arithmetic of \tHnnd{} sequences.
 Against to the background of \rrem{H.R.Hblock}, we use in the sequel the following notation:

\bnotal{H.N.L}
 If \(\seqska \) is a sequence of complex \tpqa{matrices}, then let \(\Lu{0}\defeq\Hu{0}\) and let \(\Lu{n}\defeq\Hu{n}\schca \Hu{n-1}\) for all \(n\in\N\) with \(2n\leq\kappa\).
\enota

 Obviously, \(\Hggqu{0}\) is exactly the set of sequences \((s_j)_{j=0}^0\) with \(\su{0}\in\Cggq\).
 In view of \rrem{H.R.Hblock}, we obtain from \rrem{L.AEP} furthermore:
 
\breml{H.R.Hggtr}
 If \(\seqs{2\kappa}\in\Hggqu{2\kappa}\), then \(\seqs{2n}\in\Hggqu{2n}\) and \(\Lu{n}\in\Cggq\) for all \(n\in\mn{0}{\kappa}\).
\erem

 In the following two results, we are using the equivalence relation ``\(\ldu\)'' introduced in \rnota{A.N.ldusim} (see also \rrem{A.R.ldueq}):
 
\bleml{L2108}
 Let \(n\in\NO\) and let \(\seqs{2n}\in\Hggqu{2n}\).
 Then \(\Hu{n}\ldu\diag\rk{\Lu{0},\Lu{1},\dotsc,\Lu{n}}\).
\elem
\bproof
 The case \(n=0\) is trivial.
 If \(n\geq1\), then the assertion is an easy consequence of~\zitaa{MR2805417}{\cprop{4.17}{465}}.
\eproof

\bleml{L3219}
 Let \(n\in\NO\), let \(\seqs{2n}\in\Hggqu{2n}\), and let \(A_0,A_1,\dotsc,A_n\) be complex \tqqa{matrices}.
 If \(\Hu{n}\ldu\diag\rk{A_{0},A_{1},\dotsc,A_{n}}\), then \(A_j=\Lu{j}\) for all \(j\in\mn{0}{n}\).
\elem
\bproof
 In view of \rrem{A.R.ldueq} and \rlem{L2108}, this follows from \rrem{091.YM}.
\eproof

 We write \(A\kp B\defeq\mat{a_{jk}B}_{\substack{j=1,\dotsc,m\\k=1,\dotsc,n}}\) for the Kronecker product of two matrices \(A=\mat{a_{jk}}_{\substack{j=1,\dotsc,m\\k=1,\dotsc,n}}\in\Coo{m}{n}\) and \(B\in\Cpq\).
 
\breml{H.R.HxB}
 Let \(B\in\Cpq\) and let \(\seqska\) be a sequence of complex numbers.
 Let the sequence \(\seq{x_j}{j}{0}{\kappa}\) be given by \(x_j\defeq\su{j}B\).
 Then \(\Huo{n}{x}=\Huo{n}{s}\kp B\) for all \(n\in\NO\) with \(2n\leq\kappa\) and \(\Kuo{n}{x}=\Kuo{n}{s}\kp B\) for all \(n\in\NO\) with \(2n+1\leq\kappa\).
\erem

\breml{H.R.sxBHgg}
 If \(B\in\Cggq\) and \(\seqs{2\kappa}\in\Hgguu{1}{2\kappa}\), then \(\seq{\su{j}B}{j}{0}{2\kappa}\in\Hggqu{2\kappa}\) by \rremss{A.R.AxB>=0}{H.R.HxB}.
\erem

 Now we introduce a class of finite sequences from \(\Cqq\) which turns out to be closely related to the class of \tHnnd{} sequences: 
 Let \(n\in\NO\).
 Denote by \(\Hggequ{2n}\) the set of all sequences \(\seqs{2n}\) of complex \tqqa{matrices} for which there exists a pair \((\su{2n+1},\su{2n+2})\) of complex \tqqa{matrices} such that the sequence \(\seqs{2n+2}\) belongs to \(\Hggqu{2n+2}\).
 Denote by \(\Hggequ{2n+1}\) the set of all sequences \(\seqs{2n+1}\) of complex \tqqa{matrices} for which there exists a complex \tqqa{matrix} \(\su{2n+2}\) such that the sequence \(\seqs{2n+2}\) belongs to \(\Hggqu{2n+2}\).
 Furthermore, let \(\Hggeqinf\defeq\Hggqinf\).
 The sequences belonging to \(\Hggequ{2n}\), \(\Hggequ{2n+1}\), or \(\Hggeqinf\) are said to be \emph{\tHnnde}.
 (Note that in~\zita{MR2570113} and subsequent papers of the corresponding authors, the sequences belonging to \(\Hggeqka\) were called \emph{Hankel non-negative definite extendable}.)

\bthmnl{\zitaa{MR2805417}{\cthm{6.6}{486}} (\tcf{}~\zitaa{MR1624548}{\cthm{3.1}{211}})}{H.P.MPsolv}
 Let \(\seqska \) be a sequence of complex \tqqa{matrices}.
 Then \(\MggqkappaRsg\neq\emptyset\) if and only if \(\seqska\in\Hggeqka\).
\ethm

 Now we consider the class of \tHnnde{} sequences and \tHnnd{} sequences against the background of \tftd{} sequences and \tnftd{} sequences, respectively.
 
\bpropnl{\zitaa{MR3014199}{\cprop{4.24}{136}}}{H.R.He<D}
 For each \(\kappa\in\NOinf\), we have \(\Hggequ{\kappa}\subseteq\Dqqu{\kappa}\).
\eprop
  
\bpropnl{\tcf{}~\zitaa{MR3014199}{\cprop{4.25}{137}}}{P4-25}
 For each \(n\in\NO\), the inclusion \(\Hggqu{2n}\subseteq\Dtqqu{2n}\) holds.
\eprop

 In order to describe the relations between the elements of \tHnnd{} sequences (even in the general case of not necessarily invertible \tnnH{} \tbHms{} \(\Hu{n}\)), in~\zita{MR2570113} the so-called canonical Hankel parametrization was introduced and further discussed in~\zitas{MR2805417,MR3014199}.
 We slightly reformulate this notion in a more convenient form.
 Therefore, we need several special matrices.
 
\bnotal{N.Lambda}
 Let \(\seqska \) be a sequence of complex \tpqa{matrices}.
 Then let \(\Tripu{0}\defeq\Opq\), \(\Sigmau{0}\defeq\Opq\), \(\Tripu{n}\defeq\zuu{n}{2n-1}\Hu{n-1}^\mpi\yuu{n}{2n-1}\), and \(\Sigmau{n}\defeq\zuu{n}{2n-1}\Hu{n-1}^\mpi\Ku{n-1}\Hu{n-1}^\mpi\yuu{n}{2n-1}\) for all \(n\in\N\) with \(2n-1\leq\kappa\).
 Furthermore, let \(\Mu{0}\defeq\Opq\), \(\Nu{0}\defeq\Opq\), \(\Mu{n}\defeq\zuu{n}{2n-1}\Hu{n-1}^\mpi\yuu{n+1}{2n}\), and \(\Nu{n}\defeq\zuu{n+1}{2n}\Hu{n-1}^\mpi\yuu{n}{2n-1}\) for all \(n\in\N\) with \(2n\leq\kappa\).
 Moreover, let \(\Lambdau{n}\defeq\Mu{n}+\Nu{n}-\Sigmau{n}\) for all \(n\in\NO\) with \(2n\leq\kappa\).
\enota
 
 In addition to \rrem{H.R.HxB}, we have:

\breml{H.R.yzxB}
 Let \(B\in\Cpq\) and let \(\seqska\) be a sequence of complex numbers.
 Let \(\seq{x_j}{j}{0}{\kappa}\) be given by \(x_j\defeq\su{j}B\).
 Then \(\yuuo{m}{n}{x}=\yuuo{m}{n}{s}\kp B\) and \(\zuuo{m}{n}{x}=\zuuo{m}{n}{s}\kp B\) for all \(m,n\in\mn{0}{\kappa}\).
\erem

 Using \rremsss{A.R.AxB*CxD}{A.R.AxB+}{A.R.wxB}, we obtain from \rremss{H.R.HxB}{H.R.yzxB}:

\breml{H.R.MNxB}
 Let \(B\in\Cpq\) and let \(\seqska\) be a sequence of complex numbers.
 Let the sequence \(\seq{x_j}{j}{0}{\kappa}\) be given by \(x_j\defeq\su{j}B\).
 Then \(\Tripuo{n}{x}=\Tripuo{n}{s}B\) and \(\Sigmauo{n}{x}=\Sigmauo{n}{s}B\) for all \(n\in\NO\) with \(2n-1\leq\kappa\).
 Furthermore, \(\Muo{n}{x}=\Muo{n}{s}B\) and \(\Nuo{n}{x}=\Nuo{n}{s}B\) and, consequently, \(\Lambdauo{n}{x}=\Lambdauo{n}{s}B\) for all \(n\in\NO\) with \(2n\leq\kappa\).
\erem

\bdefnl{102.HPN}
 Let \(\seqska \) be a sequence of complex \tpqa{matrices}.
 Using \rnota{N.Lambda} let the sequence \(\seq{\hpu{j}}{j}{0}{\kappa}\) be given by \(\hpu{2k}\defeq\su{2k}-\Tripu{k}\) for all \(k\in\NO\) with \(2k\leq\kappa\) and by \(\hpu{2k+1}\defeq\su{2k+1}-\Lambdau{k}\) for all \(k\in\NO\) with \(2k+1\leq\kappa\).
 Then we call \(\seq{\hpu{j}}{j}{0}{\kappa}\) the \emph{\thpfa{\(\seqska \)}}.
\edefn

 The \thp{s} were introduced in a slightly different form as a pair of two sequences in~\zita{MR2570113} under the notion \emph{canonical Hankel parametrization}.
 Several properties of a sequence connected to \tHnnd{ness} can be characterized in terms of its \thp{s}.
 They also occur in connection with three term recurrence relations for systems of orthogonal matrix polynomials with respect to \tnnH{} measures on \(\rk{\R,\BsAR}\) (see~\zita{MR2805417},~\zitaa{MR2570113}{\csec{3}}).
 
\breml{H.R.h2L}
 Let \(\seqska\) be a sequence of complex \tpqa{matrices}.
 For all \(n\in\NO\) with \(2n\leq\kappa\), then \(\hpu{2n}=\Lu{n}\).
 In particular, if \(n\geq1\), then \(\hpu{2n}\) is the Schur complement \(\Hu{n}\schca \Hu{n-1}\) of \(\Hu{n-1}\) in \(\Hu{n}\).
\erem

\breml{H.R.hxB}
 If \(B\in\Cpq\) and \(\seqska\) is a sequence of complex numbers, then \rrem{H.R.MNxB} shows that \(\seq{\hpu{j}B}{j}{0}{\kappa}\) coincides with the \thpfa{\(\seq{\su{j}B}{j}{0}{\kappa}\)}.
\erem

 Using the \trF{} given in \rdefn{D1419}, we introduce now a transformation of sequences of complex matrices, which for matricial Hamburger moment sequences corresponds to the elementary step of a Schur--Nevanlinna type algorithm:

\bdefnl{103.S31}
 Let \(\seqska\) be a sequence of complex \tpqa{matrices} and assume \(\kappa\geq2\).
 Then we call the sequence \(\seqt{\kappa-2}\) given by \(t_j\defeq-s_0 s_{j+2}^\rez s_0\) the \emph{\tHTv{\(\seqska\)}}.
\edefn

 In~\zitaa{MR3014199}{\cdefn{8.1}{159}}, the \tHT{} was introduced under the name \emph{first Schur transform}.
 For \tHnnde{} sequences \(S_0,S_1,S_2,\dotsc,S_{2n}\) it coincides with the sequence \(S_0^{(1)},\dotsc,S_{2n-2}^{(1)}\) implicitly given via~\zitaa{MR1624548}{\ceqpp{2.9}{206}{207}} by use of Drazin inverses (\tcf{}~\zitaa{MR3014199}{\cpagesbs{166}{167}}).
 For sequences \(\su{0},\su{1},\su{2},\dotsc,\su{2n-2}\) of \tH{} matrices with invertible \tnnH{} \tbHm{} \(\Hu{n-1}\), one can find a similar construction already in~\zitaa{zbMATH03875641}{\cpage{445}}, equivalent up to normalization.
 The matrices \(\su{0}^{(2)},\su{1}^{(2)},\dotsc,\su{2n-4}^{(2)}\) considered there are related to the \tHT{} \(\seqt{2n-4}\) of \(\seqs{2n-2}\) via \(\su{j}^{(2)}=\su{0}^\inv t_j\su{0}^\inv=-\su{j+2}^\rez\).

 We are now going to obtain representations of \tbHms{} built form the \tHT{} of a sequence.
 
\bpropl{H.R.HTH}
 Assume \(\kappa\geq2\) and let \(\seqska\) be a sequence of complex \tpqa{matrices} with \tHT{} \(\seqt{\kappa-2}\).
 For all \(n\in\NO\) with \(2n+2\leq\kappa\), then \(\Huo{n}{t}=-\sdiag{\su{0}}{n}\Gu{n}^\rez\sdiag{\su{0}}{n}\) and \(\Huo{n}{t}=\sdiag{\su{0}}{n}\SLu{n}^\rez\LLu{n+1}\SUu{n}^\rez\sdiag{\su{0}}{n}\).
\eprop
\bproof
 Taking into account \rdefn{103.S31} the assertions follow immediately from \rthm{H.P0938}.
\eproof

 In the case of \tnftd{} sequences we want to improve the last mentioned factorization of \(\Huo{n}{t}\) in \rprop{103.TS1A} below.
 For this reason, we introduce corresponding regularized triangular factors belonging to the classes \(\nudpu{m}\) and \(\nodqu{m}\) given in \rnota{ab.N1308}:

\bnotal{H.N.D}
 If \(\seqska\) is a sequence of complex \tpqa{matrices}, then, for all \(m\in\mn{0}{\kappa}\), let
\begin{align*}
 \Dlu{m}&\defeq\sdiag{\su{0}}{m}\SLu{m}^\rez+\sdiag{\Ip-\su{0}\su{0}^\mpi}{m}&
&\text{and}&
 \Dru{m}&\defeq\SUu{m}^\rez\sdiag{\su{0}}{m}+\sdiag{\Iq-\su{0}^\mpi\su{0}}{m}.
\end{align*}
\enota

\breml{H.R.DinLU}
 Let \(\seqska\) be a sequence of complex \tpqa{matrices}.
 Because of \(\su{0}^\rez=\su{0}^\mpi\), for all \(m\in\mn{0}{\kappa}\), then \(\Dlu{m}\) is a block Toeplitz matrix belonging to \(\nudpu{m}\) and \(\Dru{m}\) is a block Toeplitz matrix belonging to \(\nodqu{m}\).
 Furthermore, \(\det\Dlu{m}=1\) and \(\det\Dru{m}=1\) for all \(m\in\mn{0}{\kappa}\).
\erem

\bnotal{H.N.epXi}
 Let \(\seqska\) be a sequence of complex \tpqa{matrices}.
 Then let \(\defu{m}\defeq\su{m}-\su{0}\su{0}^\mpi\su{m}\su{0}^\mpi\su{0}\) for all \(m\in\mn{0}{\kappa}\) and, using \rnota{H.N.updo}, let \(\Defuu{n}{m}\defeq\OIpuu{n}{1}\defu{m}\OIquu{n}{1}^\ad\)
 for all \(n\in\NO\) and all \(m\in\mn{0}{\kappa}\).
\enota

 Obviously, we have \(\defu{0}=\Opq\) and hence \(\Defuu{n}{0}=\Ouu{(n+1)p}{(n+1)q}\) for all \(n\in\NO\).
 The following result is a slight extension of~\zitaa{MR3014199}{\cprop{8.19}{163}}.
 
\bpropl{103.TS1A}
 Let \(n\in \NO\) and let \(\seqs{2n+2}\in\Dtpqu{2n+2}\) with \tHT{} \(\seqt{2n}\).
 In view of \rnotass{H.N.D}{H.N.epXi}, then
\begin{align*}
 \Huo{n}{t}&=\sdiag{\su{0}}{n}\SLu{n}^\mpi\LLu{n+1}\SUu{n}^\mpi\sdiag{\su{0}}{n},&
 \Huo{n}{t}&=\Dlu{n}\rk{\LLu{n+1}-\Defuu{n}{2n+2}}\Dru{n},
\end{align*}
 and \(\rank\rk{\Huo{n}{t}}=\rank(\LLu{n+1}-\Defuu{n}{2n+2})\).
 If \(p=q\), then furthermore
\begin{align*}
 \det(\Huo{n}{t})&=(\det\su{0})(\det\su{0})^\mpi\det\LLu{n+1}&
&\text{and}&
 \det(\Huo{n}{t})&=\det\rk{\LLu{n+1}-\Defuu{n}{2n+2}}.
\end{align*}
\eprop
\bproof
 According to \rprop{101.S216}, we have \(\SLu{n}^\rez=\SLu{n}^\mpi\) and \(\SUu{n}^\rez=\SUu{n}^\mpi\).
 Consequently, the first asserted equation follows from \rprop{H.R.HTH} and the matrices \(\Dlu{n}\) and \(\Dru{n}\) admit the representations
\begin{align*}
 \Dlu{n}&=\sdiag{\su{0}}{n}\SLu{n}^\mpi+\sdiag{\Ip-\su{0}\su{0}^\mpi}{n}&
&\text{and}&
 \Dru{n}&=\SUu{n}^\mpi\sdiag{\su{0}}{n}+\sdiag{\Iq-\su{0}^\mpi\su{0}}{n}.
\end{align*}
 Hence, the second asserted equation is valid, by virtue of~\zitaa{MR3014199}{\cprop{8.19}{163}, \cfo{8.5}{}}.
 According to \rrem{H.R.DinLU}, the latter one implies the third one and, if \(p=q\), furthermore the fifth asserted equation.
 If \(p=q\), the fourth asserted equation follows from \(\su{0}^\rez=\su{0}^\mpi\) and \rprop{H.R.HTH}.
\eproof

 We are now going to iterate the \tHTion{} introduced in \rdefn{103.S31}:

\bdefnl{103.S321}
 Let \(\seqska\) be a sequence of complex \tpqa{matrices}.
 Let the sequence \(\seq{\su{j}^\HTa{0}}{j}{0}{\kappa}\) be given by \(\su{j}^\HTa{0}\defeq\su{j}\).
 If \(\kappa\geq2\), then for all \(k\in\N\) with \(2k\leq\kappa\), let the sequence \(\seq{\su{j}^\HTa{k}}{j}{0}{\kappa-2k}\) be recursively defined to be the \tHT{} of the sequence \(\seq{\su{j}^\HTa{k-1}}{j}{0}{\kappa-2(k-1)}\).
 For all \(k\in\NO\) with \(2k\leq\kappa\), we call the sequence \(\seq{\su{j}^\HTa{k}}{j}{0}{\kappa-2k}\) the \emph{\tnHTv{k}{\(\seq{\su{j}}{j}{0}{\kappa}\)}}.
\edefn

 From the next observations it becomes clear why the \tHTion{} is important.
 The \tHTion{} preserves \tHnnd{} extendability:

\bpropnl{\zitaa{MR3014199}{\cprop{9.4}{168}}}{H.C.HTinHgge}
 Let \(m\in\NO\) and let \(\seqs{m}\in\Hggequ{m}\).
 For all \(k\in\NO\) with \(2k\leq m\), then \(\seq{\su{j}^\HTa{k}}{j}{0}{m-2k}\in\Hggequ{m-2k}\).
\eprop

 One of the main results concerning the \tHT{} of \tHnnd{} sequences can be reformulated in our terminology as follows:
 
\bthmnl{\tcf{}~\zitaa{MR3014199}{\cthm{9.15}{178}}}{H.P.T915}
 Let \(\seqs{2\kappa}\in\Hggequ{2\kappa}\) with \thpf{} \(\hpseqo{2\kappa}\).
 Then \(\hpu{2k}=\su{0}^\HTa{k}\) for all \(k\in\mn{0}{\kappa}\) and \(\hpu{2k+1}=\su{1}^\HTa{k}\) for all \(k\in\mn{0}{\kappa-1}\).
\ethm

\section{Matricial Stieltjes moment sequences}\label{K.S}
 Let \(\ug,\obg\in\R\).
 In this section, we summarize results on \(\rhl\)\nobreakdash-Stieltjes moment sequences and \(\lhl\)\nobreakdash-Stieltjes moment sequences for our subsequent considerations.
 We start with some classes of sequences of complex matrices corresponding to solvability criteria for matricial Stieltjes moment problems on the half-line \(\rhl\):
 Let \(\Kggqu{0}\defeq\Hggqu{0}\).
 For each \(n\in\N\), denote by \(\Kggqu{2n}\) the set of all sequences \(\seqs{2n}\) of complex \tqqa{matrices} for which the \tbHms{} \(\Hu{n}\) and \(-\ug\Hu{n-1}+\Ku{n-1}\) are both \tnnH{}.
 For each \(n\in\NO\), denote by \(\Kggqu{2n+1}\) the set of all sequences \(\seqs{2n+1}\) of complex \tqqa{matrices} for which the \tbHms{} \(\Hu{n}\) and \(-\ug\Hu{n}+\Ku{n}\) are both \tnnH{}.
 Furthermore, denote by \(\Kggqinf\) the set of all sequences \(\seqsinf\) of complex \tqqa{matrices} satisfying \(\seqs{m}\in\Kggqu{m}\) for all \(m\in\NO\).
 The sequences belonging to \(\Kggqu{0}\), \(\Kggqu{2n}\), \(\Kggqu{2n+1}\), or \(\Kggqinf\) are said to be \emph{\tKnnd}.
 (Note that in~\zita{MR3014201} and subsequent papers of the corresponding authors, the sequences belonging to \(\Kggqu{\kappa}\) were called \emph{\(\ug\)\nobreakdash-Stieltjes (right-sided) non-negative definite}.
 Our terminology here differs in order to find a common notation for the three matricial power moment problems under study.)
 Finite \tKnnd{} sequences constitute exactly the class of prescribed data \(\seqs{m}\) for which a slightly modified version of the Problem~\mprob{\rhl}{m}{=} is solvable (see~\zitas{MR2735313}).

\bnotal{K.N.sa}
 Suppose \(\kappa\geq1\).
 Let \(\seqska\) be a sequence of complex \tpqa{matrices}.
 Then let the sequence \(\seqsa{\kappa-1}\) be given by \(\sau{j}\defeq-\ug\su{j}+\su{j+1}\).
 For each matrix \(X_k=X_k^{\ok{s}}\) built from the sequence \(\seqska\), we denote (if possible) by \(X_\aur{k}\defeq X_k^{\ok{\saus}}\) the corresponding matrix built from the sequence \(\seqsa{\kappa-1}\) instead of \(\seqska\).
\enota

 In view of \rnota{N.HKG}, we get in particular \(\Hau{n}=-\ug\Hu{n}+\Ku{n}\) for all \(n\in\NO\) with \(2n+1\leq\kappa\).
 In the classical case \(\ug=0\), we have furthermore \(\sau{j}=\su{j+1}\) for all \(j\in\mn{0}{\kappa-1}\).
 From \rrem{H.R.Hggtr} we obtain:

\breml{K.R.Kggtr}
 Let \(\seqska\in\Kggqu{\kappa}\).
 For all \(m\in \mn{0}{\kappa}\), then \(\seqs{m}\in\Kggqu{m}\).
 Furthermore, \(\Lu{n}\in\Cggq\) for all \(n\in\NO\) with \(2n\leq\kappa\) and \(\Lau{n}\in\Cggq\) for all \(n\in\NO\) with \(2n+1\leq\kappa\).
\erem

 For each \(m\in\NO\), denote by \(\Kggequ{m}\) the set of all sequences \(\seqs{m}\) of complex \tqqa{matrices} for which there exists a complex \tqqa{matrix} \(\su{m+1}\) such that the sequence \(\seqs{m+1}\) belongs to \(\Kggqu{m+1}\).
 Furthermore, let \(\Kggeqinf\defeq\Kggqinf\).
 The sequences belonging to \(\Kggequ{m}\) or \(\Kggeqinf\) are said to be \emph{\tKnnde}.
 It is readily checked that \(\Kggequ{0}=\Kggqu{0}\) and, for all \(n\in\N\), furthermore \(\Kggequ{2n}=\setaca{\seqs{2n}\in\Hggqu{2n}}{\seqsa{2n-1}\in\Hggequ{2n-1}}\) and
\beql{Kgg2n-1e}
 \Kggequ{2n-1}
 =\setaca*{\seqs{2n-1}\in\Hggequ{2n-1}}{\seqsa{2(n-1)}\in\Hggqu{2(n-1)}}.
\eeq

\breml{K.R.Ketr}
 If \(\seqska\in\Kggeqka\), then \rrem{K.R.Kggtr} shows that \(\seqs{m}\in\Kggequ{m}\) for all \(m\in\mn{0}{\kappa}\).
\erem

\bthmnl{\zitaa{MR3014201}{\cthm{1.6}{214}}}{dm2.T10}
 Let \(\seqska \) be a sequence of complex \tqqa{matrices}.
 Then \(\MggqKsg{\kappa}\neq\emptyset\) if and only if \(\seqska \in\Kggequu{\kappa}{\ug}\).
\ethm

 Now we consider, for an arbitrarily fixed number \(\obg\in\R\), a class of sequences of complex matrices, which is dual to the class \(\Kggqka\) and which takes the corresponding role for the moment problem on the left half-line \(\lhl\).
 (This class was already studied in~\zita{MR3014201}.):
 Let \(\Lggqu{0}\defeq\Hggqu{0}\).
 For each \(n\in\N\), denote by \(\Lggqu{2n}\) the set of all sequences \(\seqs{2n}\) of complex \tqqa{matrices} for which the \tbHms{} \(\Hu{n}\) and \(\obg\Hu{n-1}-\Ku{n-1}\) are both \tnnH{}.
 For each \(n\in\NO\), denote by \(\Lggqu{2n+1}\) the set of all sequences \(\seqs{2n+1}\) of complex \tqqa{matrices} for which the \tbHms{} \(\Hu{n}\) and \(\obg\Hu{n}-\Ku{n}\) are both \tnnH{}.
 Furthermore, denote by \(\Lggqinf\) the set of all sequences \(\seqsinf\) of complex \tqqa{matrices} satisfying \(\seqs{m}\in\Lggqu{m}\) for all \(m\in\NO\).
 The sequences belonging to \(\Lggqu{0}\), \(\Lggqu{2n}\), \(\Lggqu{2n+1}\), or \(\Lggqinf\) are said to be \emph{\tLnnd}.
 (Note that in~\zita{MR3014201} the sequences belonging to \(\Lggqu{\kappa}\) were called \emph{\(\obg\)\nobreakdash-Stieltjes left-sided non-negative definite}.)
 
\bnotal{K.N.sb}
 Suppose \(\kappa\geq1\).
 Let \(\seqska\) be a sequence of complex \tpqa{matrices}.
 Then let the sequence \(\seqsb{\kappa-1}\) be given by \(\sub{j}\defeq\obg\su{j}-\su{j+1}\).
 For each matrix \(X_k=X_k^{\ok{s}}\) built from the sequence \(\seqska\), we denote, where applicable, by \(X_\lub{k}\defeq X_k^{\ok{\subs}}\) the corresponding matrix built from the sequence \(\seqsb{\kappa-1}\) instead of \(\seqska\).
\enota

 In view of \rnota{N.HKG}, we get in particular \(\Hub{n}=\obg\Hu{n}-\Ku{n}\) for all \(n\in\NO\) with \(2n+1\leq\kappa\).
 Similar to \rrem{K.R.Kggtr}, we obtain from \rrem{H.R.Hggtr}:
 
\breml{K.R.Lggtr}
 If \(\seqska\in \Lggqu{\kappa}\), then \(\seqs{m}\in\Lggqu{m}\) for all \(m\in \mn{0}{\kappa}\) and, furthermore, \(\Lu{n}\in\Cggq\) for all \(n\in\NO\) with \(2n\leq\kappa\) and \(\Lub{n}\in\Cggq\) for all \(n\in\NO\) with \(2n+1\leq\kappa\).
\erem

 For each \(m\in\NO\), denote by \(\Lggequ{m}\) the set of all sequences \(\seqs{m}\) of complex \tqqa{matrices} for which there exists a complex \tqqa{matrix} \(\su{m+1}\) such that the sequence \(\seqs{m+1}\) belongs to \(\Lggqu{m+1}\).
 Furthermore, let \(\Lggeqinf\defeq\Lggqinf\).
 The sequences belonging to \(\Lggequ{m}\) or \(\Lggeqinf\) are said to be \emph{\tLnnde}.
 It is readily checked that \(\Lggequ{0}=\Lggqu{0}\) and, for all \(n\in\N\), furthermore \(\Lggequ{2n}=\setaca{\seqs{2n}\in\Hggqu{2n}}{\seq{\sub{j}}{j}{0}{2n-1}\in\Hggequ{2n-1}}\) and
\beql{Lgg2n-1e}
 \Lggequ{2n-1}
 =\setaca*{\seqs{2n-1}\in\Hggequ{2n-1}}{\seq{\sub{j}}{j}{0}{2(n-1)}\in\Hggqu{2(n-1)}}.
\eeq

\breml{K.R.Letr}
 If \(\seqska\in\Lggeqka\), then \rrem{K.R.Lggtr} shows that \(\seqs{m}\in\Lggequ{m}\) for all \(m\in\mn{0}{\kappa}\).
\erem

\bthmnl{\zitaa{MR3014201}{\cthm{1.8}{215}}}{121.T1609}
 Let \(\obg\in\R\) and let \(\seqska \) be a sequence of complex \tqqa{matrices}.
 Then \(\MggqLsg{\kappa}\neq\emptyset\) if and only if \(\seqska\in\Lggequu{\kappa}{\obg}\).
\ethm

\section{Matricial \(\ab\)-Hausdorff moment sequences}\label{F.S.abms}
 In the remaining part of this paper, let \(\ug\) and \(\obg\) be two arbitrarily given real numbers satisfying \(\ug<\obg\).
 In this section, we recall a collection of results on the matricial Hausdorff moment problem and corresponding moment sequences of \tnnH{} \tqqa{measures} on the interval \(\ab\), which are mostly taken from~\zitas{MR3775449,MR3979701}.
 To state a solvability criterion for the matricial \(\ab\)\nobreakdash-Hausdorff moment problem, we first extend \rnotass{K.N.sa}{K.N.sb}:
 
\bnotal{F.N.sa}
 Suppose \(\kappa\geq1\).
 Let \(\seqska\) be a sequence of complex \tpqa{matrices}.
 Then let the sequences \(\seqsa{\kappa-1}\) and \(\seqsb{\kappa-1}\) be given by
\begin{align*}
 \sau{j}& \defeq-\ug\su{j}+\su{j+1}&
&\text{and}&
 \sub{j}& \defeq\obg\su{j}-\su{j+1},
\end{align*}
 \tresp{}
 Furthermore, if \(\kappa\geq2\), let the sequence \(\seqsab{\kappa-2}\) be given by
 \[
  \sab{j}
  \defeq-\ug\obg\su{j}+(\ug+\obg)\su{j+1}-\su{j+2}.
 \]
 For each matrix \(X_k=X_k^{\ok{s}}\) built from the sequence \(\seqska\), we denote (if possible) by \(X_\aur{k}\defeq X_k^{\ok{\saus}}\), by \(X_\lub{k}\defeq X_k^{\ok{\subs}}\), and by \(X_\aub{k}\defeq X_k^{\ok{\sabs}}\) the corresponding matrix built from the sequences \(\seqsa{\kappa-1}\), \(\seqsb{\kappa-1}\), and \(\seqsab{\kappa-2}\) instead of \(\seqska\), \tresp{}
\enota

 In view of \rnota{N.HKG}, we get in particular
\begin{align}\label{F.G.HaHbHK}
 \Hau{n}&=-\ug\Hu{n}+\Ku{n}&
 &\text{and}&
 \Hub{n}&=\obg\Hu{n}-\Ku{n}
\end{align}
 for all \(n\in\NO\) with \(2n+1\leq\kappa\) and
\beql{F.G.HabHKG}
 \Hab{n}
 =-\ug\obg\Hu{n}+(\ug+\obg)\Ku{n}-\Gu{n}
\eeq
 for all \(n\in\NO\) with \(2n+2\leq\kappa\).
 Note that in the classical case where \(\ug=0\) and \(\obg=1\) hold true, we have furthermore \(\sau{j}=\su{j+1}\) and \(\sub{j}=\su{j}-\su{j+1}\) for all \(j\in\mn{0}{\kappa-1}\) and \(\sab{j}=\su{j+1}-\su{j+2}\) for all \(j\in\mn{0}{\kappa-2}\).
 The following observations are immediate from the definitions:

\breml{F.R.c=ab}
 If \(\seqska\) is a sequence of complex \tpqa{matrices}, then \(\sab{j}=-\ug\sub{j}+\sub{j+1}\) and \(\sab{j}=\obg\sau{j}-\sau{j+1}\) for all \(j\in\mn{0}{\kappa-2}\).
\erem

 In following computations, often the length
\[%
 \ba
 \defeq\obg-\ug
\]
 of the interval \(\ab\) occurs.
 Observe that \(\ba>0\), by virtue of our assumption \(\ug<\obg\).
 
\breml{F.R.012}
 Let \(\seqska\) be a sequence of complex \tpqa{matrices}.
 Then \(\ba\su{j}=\sau{j}+\sub{j}\) and \(\ba\su{j+1}=\obg\sau{j}+\ug\sub{j}\) for all \(j\in\mn{0}{\kappa-1}\).
 Furthermore, \(\ba\su{j+2}=\obg^2\sau{j}+\ug^2\sub{j}-\ba\sab{j}\) for all \(j\in\mn{0}{\kappa-2}\).
\erem

 We denote by \(\Cgq\defeq\setaca{M\in\Cqq}{v^\ad Mv\in(0,\infp)\text{ for all }v\in\Cq\setminus\set{\Ouu{q}{1}}}\)
 the set of \tpH{} matrices from \(\Cqq\).

\bdefnl{D0745}
 Let \(\Fggqu{0}\) (\tresp{}, \(\Fgqu{0}\)) be the set of all sequences \((s_j)_{j=0}^0\) with \(\su{0}\in\Cggq\) (\tresp{}, \(\su{0}\in\Cgq\)).
 For each \(n\in\N\), denote by \(\Fggqu{2n}\) (\tresp{}, \(\Fgqu{2n}\)) the set of all sequences \(\seqs{2n}\) of complex \tqqa{matrices}, for which the \tbHms{} \(\Hu{n}\) and \(\Hab{n-1}\) are both non-negative (\tresp{}, positive) \tH{}.
 For each \(n\in\NO\), denote by \(\Fggqu{2n+1}\) (\tresp{}, \(\Fgqu{2n+1}\)) the set of all sequences \(\seqs{2n+1}\) of complex \tqqa{matrices} for which the \tbHms{} \(\Hau{n}\) and \(\Hub{n}\) are both non-negative (\tresp{}, positive) \tH{}.
 Furthermore, denote by \(\Fggqinf\) (\tresp{}, \(\Fgqinf\)) the set of all sequences \(\seqsinf \) of complex \tqqa{matrices} such that \(\seqs{m}\) belongs to \(\Fggqu{m}\) (\tresp{}, \(\Fgqu{m}\)) for all \(m\in\NO\).
 The sequences belonging to \(\Fggqu{0}\), \(\Fggqu{2n}\), \(\Fggqu{2n+1}\), or \(\Fggqinf\) (\tresp{}, \(\Fgqu{0}\), \(\Fgqu{2n}\), \(\Fgqu{2n+1}\), or \(\Fgqinf\)) are said to be \emph{\(\ab\)\nobreakdash-non-negative} (\tresp{}, \emph{positive}) \emph{definite}.
\edefn

 Note that in~\zita{MR3775449} the sequences belonging to \(\Fggqu{\kappa}\) (\tresp{}, \(\Fgqu{\kappa}\)) were called \emph{\(\ab\)\nobreakdash-Hausdorff non-negative} (\tresp{}, \emph{positive}) \emph{definite}.
 The following result essentially influences our subsequent considerations.
 
\bthmnl{\tcf{}~\zitaa{MR2222521}{\cthm{1.3}{127}} and~\zitaa{MR2342899}{\cthm{1.3}{106}}}{I.P.ab}
 Let \(\seqska\) be a sequence of complex \tqqa{matrices}.
 Then \(\MggqFksg\neq\emptyset\) if and only if \(\seqska\in\Fggqka\).
\ethm

 Since \(\Omega=\ab\) is bounded, one can easily see that \(\MggqF=\MggquF{\infi}\), \tie{}, the power moment \(\mpm{\sigma}{j}\) defined by \eqref{I.G.mom} exists for all \(j\in\NO\).
 If \(\sigma\in\MggqF\), then we call \(\seqmpm{\sigma}\) given by \eqref{I.G.mom} the \emph{\tfpmfa{\(\sigma\)}}.
 Given the complete sequence of prescribed power moments \(\seqsinf\), the moment problem on the compact interval \(\Omega=\ab\) differs from the moment problems on the unbounded sets \(\Omega=\rhl\) and \(\Omega=\R\) in having necessarily a unique solution, assumed that a solution exists.
 The following result is well known and can be proved, in view of \rthm{I.P.ab}, using the corresponding result in the scalar case \(q=1\) (see, \zita{MR1544592} or~\zitaa{MR0184042}{\cthm{2.6.4}{74}}).

\bpropl{I.P.ab8}
 If \(\seqsinf\in\Fggqinf\), then the set \(\MggqFsg{\infi}\) consists of exactly one element.
\eprop

 Assuming \(\kappa=\infi\), we can summarize \rprop{I.P.ab8} and \rthm{I.P.ab}:

\bpropl{I.P.ab8Fgg}
 The mapping \(\Xi_\ab\colon\MggqF\to\Fggqinf\) given by \(\sigma\mapsto\seqmpm{\sigma}\) is well defined and bijective.
\eprop

 Observe that, for each \(\ell\in\NO\), the \emph{\(\ell\)\nobreakdash-th moment space}
\beql{msp}
 \MomqF{\ell}
 \defeq\setaca*{\col\rk{\mpm{\sigma}{0},\mpm{\sigma}{1},\dotsc,\mpm{\sigma}{\ell}}}{\sigma\in\MggqF}
\eeq
 \emph{on \(\ab\)} corresponds to the class \(\Fggqu{\ell}\) via
\[
 \MomqF{\ell}
 =\setaca*{\yuu{0}{\ell}}{\seqs{\ell}\in\Fggqu{\ell}}
\]
 (\tcf{}~\zitaa{MR1883272}{\cthm{2.6}{177}}).
 For all \(\iota\in\NOinf\) and each non-empty set \(\mathcal{X}\), let
\(
 \seqset{\iota}{\mathcal{X}}
 \defeq\setaca{\seq{X_j}{j}{0}{\iota}}{X_j\in\mathcal{X}\text{ for all }j\in\mn{0}{\iota}}
\).
 Obviously, \(\Fggqu{0}\) coincides with the set of all sequences \((s_j)_{j=0}^0\) with \(\su{0}\in\Cggq\).
 Furthermore, using \rnota{F.N.sa}, we have
\beql{Fgg2n}
 \Fggqu{2n}
 =\setaca*{\seqs{2n}\in\Hggqu{2n}}{\seq{\sab{j}}{j}{0}{2(n-1)}\in\Hggqu{2(n-1)}}
\eeq
 for all \(n\in\N\) and
\beql{Fgg2n+1}
 \Fggqu{2n+1}
 =\setaca*{\seqs{2n+1}\in\seqset{2n+1}{\Cqq}}{\set{\seqsa{2n},\seqsb{2n}}\subseteq\Hggqu{2n}}
\eeq
 for all \(n\in\NO\).
 The following \rpropsss{ab.R0933}{ab.R0951}{F.L.sabF}, which are proved in a purely algebraic way in~\zita{MR3775449}, can also be obtained immediately from \rthmss{I.P.ab}{H.P.MPsolv}.
 
\bpropnl{\tcf{}~\zitaa{MR3775449}{\cprop{7.7(a)}{18}}}{ab.R0933}
 If \(\seqska  \in \Fggqka \), then \(\seqs{m} \in \Fggqu{m}\) for all \(m\in \mn{0}{\kappa}\).
\eprop

 In view of \rprop{ab.R0933}, the definition of the class \(\Fggqinf\) seems to be natural.
 Taking additionally into account \eqref{Fgg2n} and \eqref{Fgg2n+1}, we have furthermore
\begin{equation}\label{Fgg8}
 \begin{split}
  \Fggqinf
  &=\setaca*{\seqsinf\in\Hggqinf}{\seqsabinf\in\Hggqinf}\\
  &=\setaca*{\seqsinf\in\seqset{\infi}{\Cqq}}{\set{\seqsainf,\seqsbinf}\subseteq\Hggqinf}.
 \end{split}
\end{equation}
 
\bpropnl{\tcf{}~\zitaa{MR3775449}{\cprop{7.10}{19}}}{ab.R0951}
 If \(m\in\NO\), then \(\Fggqu{m}\subseteq\Hggequ{m}\).
\eprop

\bpropnl{\zitaa{MR3775449}{\cprop{9.1}{27}}}{F.L.sabF}
 Let \(\seqska\in\Fggqka\).
 If \(\kappa\geq1\), then \(\set{\seqsa{\kappa-1},\seqsb{\kappa-1}}\subseteq\Fggqu{\kappa-1}\).
 Furthermore, if \(\kappa\geq2\), then \(\seqsab{\kappa-2}\in\Fggqu{\kappa-2}\).
\eprop

\bpropl{F.R.Fgg<D}
\benui
 \il{F.R.Fgg<D.a} \(\Fggqka\subseteq\Dqqu{\kappa}\).
 \il{F.R.Fgg<D.b} If \(\kappa\geq1\) and if \(\seqska\in\Fggqka\), then \(\set{\seqsa{\kappa-1},\seqsb{\kappa-1}}\subseteq\Dqqu{\kappa-1}\).
 \il{F.R.Fgg<D.c} If \(\kappa\geq2\) and if \(\seqska\in\Fggqka\), then \(\seqsab{\kappa-2}\in\Dqqu{\kappa-2}\).
\eenui
\eprop
\bproof
 \eqref{F.R.Fgg<D.a} Because of \rprop{H.R.He<D}, this follows in the case \(\kappa<\infi\) from \rprop{ab.R0951} and in the case \(\kappa=\infp\) from \eqref{Fgg8}, since \(\Hggqinf=\Hggeqinf\) by definition.
 
 \eqref{F.R.Fgg<D.b},~\eqref{F.R.Fgg<D.c} Combine~\eqref{F.R.Fgg<D.a} and \rprop{F.L.sabF}. 
\eproof
 
 The set \(\CHq\defeq\setaca{M\in\Cqq}{M^\ad=M}\) of \tH{} matrices from \(\Cqq\) is a partially ordered vector space over the field \(\R\) with positive cone \(\Cggq\).
 For two complex \tqqa{matrices} \(A\) and \(B\), we write \(A\lleq B\) or \(B\lgeq A\) if \(A,B\in\CHq\) and \(B-A\in\Cggq\) are fulfilled.
 The above mentioned partial order \(\lleq\) 
 is sometimes called \emph{L\"owner semi-ordering} in \(\CHq\).
 
\blemnl{\zitaa{MR3979701}{\clem{5.7}{15}}}{F.R.Fgg-s}
 Let \(\seqska  \in \Fggqka\).
 Then \(\su{j}\in\CHq\) for all \(j\in\mn{0}{\kappa}\) and \(\su{2k}\in\Cggq\) for all \(k\in\NO\) with \(2k\leq\kappa\).
 Furthermore, \(\ug\su{2k}\lleq\su{2k+1}\lleq\obg\su{2k}\) for all \(k\in\NO\) with \(2k+1\leq\kappa\).
\elem

\breml{F.R.Fgg-r}
 If \(\seqska\in\Fggqka\), then \rrem{A.R.rA<rB} and \rlem{F.R.Fgg-s} yield \(\ran{\sau{2k}}\cup\ran{\sub{2k}}\subseteq\ran{\su{2k}}\) and \(\nul{\su{2k}}\subseteq\nul{\sau{2k}}\cap\nul{\sub{2k}}\) for all \(k\in\NO\) with \(2k\leq\kappa-1\).
\erem

 Finite sequences from \(\Fggqu{m}\) can always be extended to sequences from \(\Fggqu{\ell}\) for all \(\ell\in\minf{m+1}\), which becomes clear from \rthm{I.P.ab} and the fact that a \tnnH{} measure on the bounded set \(\ab\) possesses power moments of all non-negative orders.
 One of the main results in~\zita{MR3775449} states that the possible one-step extensions \(\su{m+1}\in\Cqq\) of a sequence \(\seqs{m}\) to an \tFnnd{} sequence \(\seqs{m+1}\) fill out a matricial interval.
 In order to give an exact description of this interval, we are now going to introduce several matrices and recall their role in the corresponding extension problem for \tFnnd{} sequences, studied in~\zita{MR3775449}.
 
\bdefnl{D0750}
 If \(\seqska\) is a sequence of complex \tpqa{matrices}, then (using the \rnotasss{N.Lambda}{K.N.sa}{K.N.sb}) the sequences \(\seq{\umg{j}}{j}{0}{\kappa}\) and \(\seq{\omg{j}}{j}{0}{\kappa}\) given by \(  \umg{2k}\defeq\ug\su{2k}+\Tripa{k}\) and \(\omg{2k}\defeq\obg\su{2k}-\Tripb{k}\) for all \(k\in\NO\) with \(2k\leq\kappa\) and by \(  \umg{2k+1}\defeq\Tripu{k+1}\) and \(\omg{2k+1}\defeq-\ug\obg\su{2k}+(\ug+\obg)\su{2k+1}-\Tripab{k}\) for all \(k\in\NO\) with \(2k+1\leq\kappa\) are called the \emph{\tflep{\(\seqska\)}} and the \emph{\tfrep{\(\seqska\)}}, \tresp{}
\edefn

 By virtue of \rnota{N.Lambda}, we have in particular
\begin{align}\label{F.G.uo01}
 \umg{0}&=\ug\su{0},&
 \omg{0}&=\obg\su{0},&
 \umg{1}&=\su{1}\su{0}^\mpi\su{1},&
&\text{and}&
 \omg{1}&=-\ug\obg\su{0}+(\ug+\obg)\su{1}.
\end{align}

\breml{ab.L0911}
 If \(\seqska\in\Fggqka\), then \rlem{F.R.Fgg-s} and \rrem{A.R.A++*} yield \(\set{\umg{j},\omg{j}}\subseteq\CHq\) for all \(j\in\mn{0}{\kappa}\).
\erem
 
 Observe that for arbitrarily given \tH{} \tqqa{matrices} \(A\) and \(B\), the (closed) matricial interval
\[
 \matint{A}{B}
 \defeq\setaca{X\in\CHq}{A\lleq X\lleq B}
\]
 is non-empty if and only if \(A\lleq B\).
 
\bthmnl{\zitaa{MR3775449}{\cthm{11.2(a)}{44}}}{165.T112}
 If \(m\in\NO\) and \(\seqs{m}\in\Fggqu{m}\), then the matricial interval \(\matint{\umg{m}}{\omg{m}}\) is non-empty and coincides with the set of all complex \tqqa{matrices} \(\su{m+1}\) for which \(\seqs{m+1}\) belongs to \(\Fggqu{m+1}\).
\ethm

\bdefnl{D1861}
 If \(\seqska \) is a sequence of complex \tpqa{matrices}, then we call \(\seq{\dia{j}}{j}{0}{\kappa}\) be given by
\(
  \dia{j}
  \defeq\omg{j}-\umg{j}
\)
 the \emph{\tfdfa{\(\seqska \)}}.
\edefn

 By virtue of \eqref{F.G.uo01}, we have in particular
\begin{align}\label{F.G.d01}
 \dia{0}&=\ba\su{0}&
 &\text{and}&
 \dia{1}&=-\ug\obg\su{0}+(\ug+\obg)\su{1}-\su{1}\su{0}^\mpi\su{1}.
\end{align}

\breml{F.R.diatr}
 Let \(\seqska\) be a sequence of complex \tpqa{matrices} with \tfdf{} \(\seqdiaka \).
 For each \(k\in\mn{0}{\kappa}\), the matrix \(\dia{k}\) is built  from the matrices \(\su{0},\su{1},\dotsc,\su{k}\).
 In particular, for each \(m\in\mn{0}{\kappa}\), the \tfdfa{\(\seqs{m}\)} coincides with \(\seqdia{m}\).
\erem

 Using \rrem{H.R.MNxB}, we can conclude:

\breml{F.R.dxB}
 If \(B\in\Cpq\) and if \(\seqska\) is a sequence of complex numbers, then \(\seq{\dia{j}B}{j}{0}{\kappa}\) coincides with the \tfdfa{\(\seq{\su{j}B}{j}{0}{\kappa}\)}.
\erem

\bdefnl{D2972}
 Let \(\seqska\) be a sequence of complex \tpqa{matrices}.
 Then the sequence \(\seq{\usc{j}}{j}{0}{\kappa}\) given by \(\usc{0}\defeq\su{0}\) and by \(\usc{j}\defeq\su{j}-\umg{j-1}\) is called the \emph{\tflsc{\(\seqska\)}}.
 Furthermore, if \(\kappa\geq1\), then the sequence \(\seq{\osc{j}}{j}{1}{\kappa}\) given by \(\osc{j}\defeq\omg{j-1}-\su{j}\) is called the \emph{\tfusc{\(\seqska\)}}.
\edefn

 Because of \eqref{F.G.uo01}, we have in particular
\begin{align}\label{F.G.AB1B2}
 \usc{1}&=\sau{0},&
 \osc{1}&=\sub{0},&
 &\text{and}&
 \osc{2}&=\sab{0}.
\end{align}

\breml{F.R.ABL}
 Let \(\seqska \) be a sequence of complex \tpqa{matrices}.
 Then \(\usc{2n}=\Lu{n}\) (\tresp{}, \(\usc{2n+1}=\Lau{n}\)) for all \(n\in\NO\) with \(2n\leq\kappa\) (\tresp{}, \(2n+1\leq\kappa\)).
 In particular, if \(n\geq1\), then \(\usc{2n}\) is the Schur complement of \(\Hu{n-1}\) in \(\Hu{n}\) and \(\usc{2n+1}\) is the Schur complement of \(\Hau{n-1}\) in \(\Hau{n}\). 
 Furthermore, \(\osc{2n+1}=\Lub{n}\) (\tresp{}, \(\osc{2n+2}=\Lab{n}\)) for all \(n\in\NO\) with \(2n+1\leq\kappa\) (\tresp{}, \(2n+2\leq\kappa\)).
 In particular, if \(n\geq1\), then \(\osc{2n+1}\) is the Schur complement of \(\Hub{n-1}\) in \(\Hub{n}\) and \(\osc{2n+2}\) is the Schur complement of \(\Hab{n-1}\) in \(\Hab{n}\).
\erem

\breml{ab.R1420}
 If \(\seqska\) is a sequence of complex \tpqa{matrices}, then \(\dia{j}=\usc{j+1}+\osc{j+1}\) for all \(j\in\mn{0}{\kappa-1}\).
\erem

\bpropnl{\zitaa{MR3775449}{\cprop{10.15}{200}}}{ab.C0929}
 If \(\seqska\) belongs to \(\Fggqka\) (\tresp{} \(\Fgqu{\kappa}\)), then \(\dia{j}\) is non-negative (\tresp{}, positive) \tH{} for all \(j\in\mn{0}{\kappa}\).
\eprop

\bpropnl{\zitaa{MR3775449}{\cprop{10.18}{38}}}{ab.R13371422}
 Let \(\seqska\in\Fggqka\).
 Then \(\ran{\dia{0}}=\ran{\usc{0}}\) and \(\nul{\dia{0}}=\nul{\usc{0}}\).
 Furthermore, \(\ran{\dia{j}}=\ran{\usc{j}}\cap\ran{\osc{j}}\) and \(\nul{\dia{j}}=\nul{\usc{j}}+\nul{\osc{j}}\) for all \(j\in\mn{1}{\kappa}\), and \(\ran{\dia{j}}=\ran{\usc{j+1}}+\ran{\osc{j+1}}\) and \(\nul{\dia{j}}=\nul{\usc{j+1}}\cap\nul{\osc{j+1}}\) for all \(j\in\mn{0}{\kappa-1}\).
\eprop

\bpropnl{\zitaa{MR3775449}{\ccor{10.20}{39}}}{ab.C1101}
 Let \(\seqska\in\Fggqka\).
 For all \(j\in\mn{1}{\kappa}\), then
 \(\ran{\dia{j}}\subseteq\ran{\dia{j-1}}\), \(\nul{\dia{j-1}}\subseteq\nul{\dia{j}}\) and \(\ran{\usc{j}}\subseteq\ran{\usc{j-1}}\), \(\nul{\usc{j-1}}\subseteq\nul{\usc{j}}\).
 For all \(j\in\mn{2}{\kappa}\), furthermore, \(\ran{\osc{j}}\subseteq\ran{\osc{j-1}}\) and \(\nul{\osc{j-1}}\subseteq\nul{\osc{j}}\).
\eprop

 For sequences belonging to the class \(\Dpqka\), the matrix \(\dia{1}\) admits a symmetric multiplicative representation:
 
\breml{ab.L0907}
 Suppose \(\kappa\geq1\).
 If \(\seqska \in\Dpqka\), then one can easily see from \rrem{R.AA+B} and \eqref{F.G.d01} that \(\dia{1}=\sau{0}\su{0}^\mpi\sub{0}\) and \(\dia{1}=\sub{0}\su{0}^\mpi\sau{0}\).
\erem

 In~\zitaa{MR3979701}{\cdefn{6.1}{19}}, we subsumed the Schur complements mentioned in \rrem{F.R.ABL} to a new parameter sequence:

\bdefnl{D0752}
 Let \(\seqska \) be a sequence of complex \tpqa{matrices}.
 Let the sequence \(\fpseqka\) be given by \(\fpu{0}\defeq\usc{0}\), by \(\fpu{4k+1}\defeq\usc{2k+1}\) and \(\fpu{4k+2}\defeq\osc{2k+1}\) for all \(k\in\NO\) with \(2k+1\leq\kappa\), and by \(\fpu{4k+3}\defeq\osc{2k+2}\) and \(\fpu{4k+4}\defeq\usc{2k+2}\) for all \(k\in\NO\) with \(2k+2\leq\kappa\). 
 Then we call \(\fpseqka\) the \emph{\tfpfa{\(\seqska \)}}.
\edefn

 In view of \eqref{F.G.AB1B2} and \eqref{F.G.uo01}, we have in particular
\begin{align}
 \fpu{0}&=\su{0},&
 \fpu{1}&=\sau{0}=\su{1}-\ug\su{0},&
&\text{and}&
 \fpu{2}&=\sub{0}=\obg\su{0}-\su{1}.\label{F.G.f012}
\end{align}
 Note that \tFnnd{ness} can be characterized in terms of \tfp{s}, as well as rank constellations among the \tnnH{} \tbHms{} \(\Hu{n}\), \(\Hau{n}\), \(\Hub{n}\), and \(\Hab{n}\) (\tcf~\zitaa{MR3979701}{\cpropss{6.13}{22}{6.14}{22}}).

\breml{F.R.fp12}
 If \(\seqska\) is a sequence of complex \tpqa{matrices}, then \(\set{\fpu{2m-1},\fpu{2m}}=\set{\usc{m},\osc{m}}\) for all \(m\in\mn{1}{\kappa}\).
\erem

\breml{F.R.fxB}
 If \(B\in\Cpq\) and \(\seqska\) is a sequence of complex numbers, then from \rrem{H.R.MNxB} one can see that \(\seq{\fpu{j}B}{j}{0}{2\kappa}\) coincides with the \tfpfa{\(\seq{\su{j}B}{j}{0}{\kappa}\)}.
\erem

\breml{F.R.f2n-1}
 If \(\seqska\) is a sequence of complex \tpqa{matrices}, then \rremss{ab.R1420}{F.R.fp12} show that \(\fpu{2k-1}=\dia{k-1}-\fpu{2k}\) for all \(k\in\mn{1}{\kappa}\).
\erem

 Let \(\seqska\in\Fggqka\) with associated \tfpf{} \(\fpseqka\) (see \rdefn{D0752}) and associated \tfdf{} \(\seqdiaka\).
 Then we are going to describe an important connection between these objects, which is described with the aid of parallel sum of matrices.
 For this reason, we recall first the construction of parallel sum:
 If \(A\) and \(B\) are two complex \tpqa{matrices}, then the matrix
\beql{ps}
 A\ps B
 \defeq A(A+B)^\mpi B
\eeq
 is called the \emph{parallel sum of \(A\) and \(B\)}.
 We have the following connection between \tfp{s} and \tfd{s} associated to a sequence belonging to \(\Fggqka\):

\bpropl{F.R.d+=f}
 If \(\seqska\in\Fggqka\), then \(\dia{0}=\ba\fpu{0}\) and, for all \(k\in\mn{1}{\kappa}\), furthermore \(\dia{k}=\ba\rk{\fpu{2k-1}\ps\fpu{2k}}\) and \(\dia{k}=\ba\rk{\fpu{2k}\ps\fpu{2k-1}}\).
\eprop
\bproof
 In view of \rrem{F.R.fp12}, this is an immediate consequence of~\zitaa{MR3775449}{\cthm{10.14}{32}}.
\eproof

 For each matrix \(A\in\Cggq\), we use \(Q=A^\varsqrt\) to denote the  \tnnH{} square root of \(A\).
 To uncover relations between the \tfp{s} \(\fpseqka\) and obtain a parametrization of the set \(\Fggqka\), we introduced in~\zitaa{MR3979701}{\cdefn{6.21}{24} and \cnota{6.28}{26}} another parameter sequence \(\seqciaka\) and a corresponding class \(\esqkad\) of sequences of complex matrices.
 (Observe that these constructions are well defined due to \rprop{ab.C0929} and \rrem{A.R.XAX}.)

\bdefnl{D0754}
 Let \(\seqska\in\Fggqka\) with \tfpf{} \(\fpseqka\) and \tfdf{} \(\seqdiaka\).
 Then we call \(\seqciaka\) given by
\begin{align*}
 \cia{0}&\defeq\fpu{0}&
&\text{and by}&
 \cia{j}&\defeq\rk{\dia{j-1}^\varsqrt}^\mpi\fpu{2j}\rk{\dia{j-1}^\varsqrt}^\mpi
\end{align*}
 for each \(j\in\mn{1}{\kappa}\) the \emph{\tfcfa{\(\seqska\)}}.
\edefn

 With the Euclidean scalar product \(\ipE{\cdot}{\cdot}\colon\x{\Cq}\to\C\) given by \(\ipE{x}{y}\defeq y^\ad x\), which is \(\C\)\nobreakdash-linear in its first argument, the vector space \(\Cq\) over the field \(\C\) becomes a unitary space.
 Let \(\mathcal{U}\) be an arbitrary non-empty subset of \(\Cq\).
 The orthogonal complement \(\mathcal{U}^\orth\defeq\setaca{v\in\Cq}{\ipE{v}{u}=0\text{ for all }u\in\mathcal{U}}\) of \(\mathcal{U}\) is a linear subspace of the unitary space \(\Cq\).
 If \(\mathcal{U}\) is a linear subspace itself, the unitary space \(\Cq\) is the orthogonal sum of \(\mathcal{U}\) and \(\mathcal{U}^\orth\).
 In this case, we write \(\OPu{\mathcal{U}}\) for the transformation matrix corresponding to the orthogonal projection onto \(\mathcal{U}\) with respect to the standard basis of \(\Cq\), \tie{}, \(\OPu{\mathcal{U}}\) is the uniquely determined matrix \(P\in\Cqq\) satisfying \(P^2=P=P^\ad\) and \(\ran{P}=\mathcal{U}\).

\bnotal{F.N.01seq}
 For each \(\eta\in[0,\infp)\), let \(\es{q}{\kappa}{\eta}\) be the set of all sequences \(\seq{e_k}{k}{0}{\kappa}\) from \(\Cggq\) which fulfill the following condition.
 If \(\kappa\geq1\), then \(e_k\lleq\OPu{\ran{d_{k-1}}}\) for all \(k\in\mn{1}{\kappa}\), where the sequence \(\seq{d_k}{k}{0}{\kappa}\) is recursively given by \(d_0\defeq\eta e_0\) and
\[
 d_k
 \defeq\eta d_{k-1}^\varsqrt  e_k^\varsqrt (\OPu{\ran{d_{k-1}}}-e_k) e_k^\varsqrt  d_{k-1}^\varsqrt.
\]
\enota

 Regarding \rthm{165.T112}, in the case \(q=1\) (\tcf{}~\zitaa{MR1468473}{\csec{1.3}}), the (classical) canonical moments \(p_1,p_2,p_3,\dotsc\) of a point in the moment space corresponding to a probability measure \(\mu\) on \(\ab=[0,1]\) are given in our notation by %
\begin{align}\label{F.G.p=A/d}
 p_k
 &=\frac{\su{k}-\umg{k-1}}{\omg{k-1}-\umg{k-1}}
 =\frac{\usc{k}}{\dia{k-1}},&k&\in\N,
\end{align}
 where the sequence \(\seqsinf\) of power moments \(\su{j}\defeq\int_{[0,1]}x^j\mu\rk{\dif x}\) associated with \(\mu\) is \txyFnnd{0}{1} with \(\su{0}=1\).
 Observe that the \tfup{0}{1}s \(\seqciainf\) of \(\seqsinf\) are connected to the canonical moments via
\begin{align}\label{F.G.epq}
 p_1&=1-\cia{1},&
 p_2&=\cia{2},&
 p_3&=1-\cia{3},&
 p_4&=\cia{4},&
 p_5&=1-\cia{5},&
 &\dotsc
\end{align}
 By virtue of \rrem{ab.R1420}, we have, for all \(\ell\in\N\), indeed
\[\begin{split}
 \fpu{2\rk{2\ell-1}}
 =\fpu{4\ell-2}
 =\fpu{4\rk{\ell-1}+2}
 &=\osc{2\rk{\ell-1}+1}\\
 &=\dia{2\rk{\ell-1}}-\usc{2\rk{\ell-1}+1}
 =\dia{\rk{2\ell-1}-1}-\usc{2\ell-1}
\end{split}\]
 and
\[
 \fpu{2\rk{2\ell}}
 =\fpu{4\ell}
 =\fpu{4\rk{\ell-1}+4}
 =\usc{2\rk{\ell-1}+2}
 =\usc{2\ell},
\]
 implying together with \eqref{F.G.p=A/d} the relations \(\cia{2\ell-1}=1-p_{2\ell-1}\) and \(\cia{2\ell}=p_{2\ell}\).
 The quantities \(q_k=1-p_k\) occur in the classical framework as well (see, \teg{}~\zitaa{MR1468473}{\csec{1.3}}).

\bthmnl{\zitaa{MR3979701}{\cthm{6.30}{29}}}{F.T.Fggcia}
 The mapping \(\Sigma_{\ug,\obg}\colon\Fggqka\to\esqkad\) given by \(\seqska\mapsto\seqciaka\) is well defined and bijective.
\ethm

\bexaml{F.E.e=1/2}
 Let \(\lambda\in(0,1)\), let \(B\in\Cggq\), and let \(P\defeq\OPu{\ran{B}}\).
 We show that the infinite sequence \(B,\lambda P,\lambda P,\lambda P,\dotsc\) belongs to \(\es{q}{\infi}{\ba}\) and, by virtue of \rthm{F.T.Fggcia}, hence it is the \tfcfa{some} sequence \(\seqsinf\in\Fggqinf\):
  
 Let \(\eta\defeq\ba\), let the sequence \(\seq{e_k}{k}{0}{\infi}\) be given by \(e_0\defeq B\) and by \(e_k\defeq\lambda P\) for all \(k\in\N\), and let the sequence \(\seq{d_k}{k}{0}{\infi}\) be given by \(d_k\defeq\eta\ek{\eta\lambda\rk{1-\lambda}}^kB\).
 Then \(\seq{e_k}{k}{0}{\infi}\) and \(\seq{d_k}{k}{0}{\infi}\) are sequences of \tnnH{} \tqqa{matrices} fulfilling \(e_k\lleq P=\OPu{\ran{d_{k-1}}}\) for all \(k\in\N\).
 Observe that \(PB^\varsqrt=B^\varsqrt\).
 For all \(k\in\N\), we thus have
\[\begin{split}
 &\eta d_{k-1}^\varsqrt e_{k}^\varsqrt\rk{\OPu{\ran{d_{k-1}}}-e_{k}}e_{k}^\varsqrt d_{k-1}^\varsqrt\\
 &=\eta\rk*{\eta\ek*{\eta\lambda\rk{1-\lambda}}^{k-1}B}^\varsqrt\rk{\lambda P}^\varsqrt\rk{P-\lambda P}\rk{\lambda P}^\varsqrt\rk*{\eta\ek*{\eta\lambda\rk{1-\lambda}}^{k-1}B}^\varsqrt\\
 &=\eta\lambda\rk{1-\lambda}\rk*{\eta\ek*{\eta\lambda\rk{1-\lambda}}^{k-1}B}
 =\eta\ek{\eta\lambda\rk{1-\lambda}}^{k}B
 =d_k.
\end{split}\]
 Taking additionally into account \(d_0=\eta B=\eta e_0\), we can conclude then that \(\seq{d_k}{k}{0}{\infi}\) is exactly the sequence associated to \(\seq{e_k}{k}{0}{\infi}\) in \rnota{F.N.01seq}.
 Hence, the sequence \(\seq{e_k}{k}{0}{\infi}\) belongs to \(\es{q}{\infi}{\eta}\).
 By virtue of \rthm{F.T.Fggcia}, there exists a uniquely determined sequence \(\seqsinf\in\Fggqinf\) with \tfcf{} \(\seq{e_j}{j}{0}{\infi}\).
 From~\zitaa{MR3979701}{\cprop{6.32}{30}} we can furthermore infer that \(\seq{d_j}{j}{0}{\infi}\) is the \tfdfa{\(\seqsinf\)}.
 Denote by \(\fpseq{\infi}\) the \tfpfa{\(\seqsinf\)}.
 Then \(\fpu{0}=e_0=B\).
 Using~\zitaa{MR3979701}{\cprop{6.27}{26}}, for all \(j\in\N\) we get moreover
\[\begin{split}
 \fpu{2j-1}
 &= d_{j-1}^\varsqrt \rk{\OPu{\ran{d_{j-1}}}-e_j} d_{j-1}^\varsqrt\\
 &=\rk*{\eta\ek*{\eta\lambda\rk{1-\lambda}}^{j-1}B}^\varsqrt\rk{P-\lambda P}\rk*{\eta\ek*{\eta\lambda\rk{1-\lambda}}^{j-1}B}^\varsqrt\\
 &=\rk{1-\lambda}\rk*{\eta\ek*{\eta\lambda\rk{1-\lambda}}^{j-1}B}
 =\lambda^{j-1}\ek*{\eta\rk{1-\lambda}}^jB
\end{split}\]
 and, analogously,
\[\begin{split}
 \fpu{2j}
 &= d_{j-1}^\varsqrt e_j d_{j-1}^\varsqrt\\
 &=\rk*{\eta\ek*{\eta\lambda\rk{1-\lambda}}^{j-1}B}^\varsqrt\rk{\lambda P}\rk*{\eta\ek*{\eta\lambda\rk{1-\lambda}}^{j-1}B}^\varsqrt
 =\rk{\lambda\eta}^j\rk{1-\lambda}^{j-1}B.
\end{split}\]
\eexam

 The following result supplements~\zitaa{MR3979701}{\crem{6.24}{27}}:
 
\bleml{F.L.exB}
 Let \(B\in\Cggq\) and let \(\seqska\in\Fgguuuu{1}{\kappa}{\ug}{\obg}\).
 Let the sequence \(\seq{x_j}{j}{0}{\kappa}\) be given by \(x_j\defeq\su{j}B\) and denote by \(\seq{\mathfrak{p}_j}{j}{0}{\kappa}\) the \tfcfa{\(\seq{x_j}{j}{0}{\kappa}\)}.
 Then \(\seq{x_j}{j}{0}{\kappa}\in\Fggqka\) and, furthermore, \(\mathfrak{p}_{0}=\cia{0}B\) and \(\mathfrak{p}_{j}=\cia{j}\OPu{\ran{B}}\) for all \(j\in\mn{1}{\kappa}\).
\elem
\bproof
 Using \rrem{H.R.sxBHgg} and taking into account \eqref{Fgg2n}, \eqref{Fgg2n+1}, and \eqref{Fgg8}, it is readily checked that \(\seq{x_j}{j}{0}{\kappa}\) belongs to \(\Fggqka\).
 Denote by \(\seq{\mathfrak{l}_j}{j}{0}{\kappa}\) the \tfdfa{\(\seq{x_j}{j}{0}{\kappa}\)} and by \(\seq{\mathfrak{g}_j}{j}{0}{\kappa}\) the \tfpfa{\(\seq{x_j}{j}{0}{\kappa}\)}.
 According to \rrem{F.R.fxB}, then \(\mathfrak{g}_j=\fpu{j}B\) for all \(j\in\mn{0}{2\kappa}\).
 In particular, \(\mathfrak{p}_0=\mathfrak{g}_0=\fpu{0}B=\cia{0}B\).
 Now assume \(\kappa\geq1\).
 Observe that \(\seqdiaka\) and \(\fpseqka\) are sequences of complex numbers.
 Using \rrem{F.R.dxB}, we obtain \(\mathfrak{l}_j=\dia{j}B\) for all \(j\in\mn{0}{\kappa}\).
 Consequently, we have \(\rk{\mathfrak{l}_j^\varsqrt}^\mpi=\rk{\dia{j}^\varsqrt}^\mpi\rk{B^\varsqrt}^\mpi\) for all \(j\in\mn{0}{\kappa}\), by virtue of \rrem{A.R.l*A}.
 From \rrem{ab.R1052} we infer furthermore \(B^\varsqrt\rk{B^\varsqrt}^\mpi=\OPu{\ran{B^\varsqrt}} =\OPu{\ran{B}}\).
 Taking additionally into account \rrem{A.R.A++*}, we can conclude \(\OPu{\ran{B}}^\ad=\ek{B^\varsqrt\rk{B^\varsqrt}^\mpi}^\ad=\rk{B^\varsqrt}^\mpi B^\varsqrt\).
 For all \(j\in\mn{1}{\kappa}\), hence,
\[\begin{split}
 \mathfrak{p}_j
 &=\rk{\mathfrak{l}_{j-1}^\varsqrt}^\mpi\mathfrak{g}_{2j}\rk{\mathfrak{l}_{j-1}^\varsqrt}^\mpi
 =\ek*{\rk{\dia{j-1}^\varsqrt}^\mpi\rk{B^\varsqrt}^\mpi}\rk{\fpu{2j}B}\ek*{\rk{\dia{j-1}^\varsqrt}^\mpi\rk{B^\varsqrt}^\mpi}\\
 &=\ek*{\rk{\dia{j-1}^\varsqrt}^\mpi\fpu{2j}\rk{\dia{j-1}^\varsqrt}^\mpi}\ek*{\rk{B^\varsqrt}^\mpi B\rk{B^\varsqrt}^\mpi}
 =\cia{j}\OPu{\ran{B}}^\ad\OPu{\ran{B}}
 =\cia{j}\OPu{\ran{B}}.\qedhere
\end{split}\]
\eproof

 In the remaining part of this section, we recall two special properties for sequences \(\seqska\in\Fggqka\), already considered in~\zita{MR3775449} and characterized in~\zita{MR3979701} in terms of their \tfc{s} \(\seqciaka\).
 In particular, the geometry of a matricial interval suggests that in addition to the end points \(\umg{m}\) and \(\omg{m}\) of the extension interval in \rthm{165.T112}, the center of the matricial interval \(\matint{\umg{m}}{\omg{m}}\) is of particular interest:
 
\bdefnl{D0756}
 If \(\seqska \) is a sequence of complex \tpqa{matrices}, then we call \(\seq{\mi{j}}{j}{0}{\kappa}\) given by \(\mi{j}\defeq\frac{1}{2}\rk{\umg{j}+\omg{j}}\) the \emph{\tfmfa{\(\seqska \)}}.
\edefn

 It should be mentioned that the choice \(\su{j+1}=\mi{j}\) corresponds to the maximization of the width \(\dia{j+1}\) of the corresponding section of the moment space \eqref{msp}, as was shown in~\zitaa{MR3775449}{\cprop{10.23}{39}}.

\bdefnnl{\tcf{}~\zitaa{MR3775449}{\cdefn{10.33}{42}}}{F.D.abZ}
 Let \(\seqska\) be a sequence of complex \tpqa{matrices}, assume \(\kappa\geq1\), and let \(k\in\mn{1}{\kappa}\).
 Then \(\seqska\) is said to be \emph{\tabZo{k}} if \(\su{j}=\mi{j-1}\) for all \(j\in\mn{k}{\kappa}\).
\edefn

\bnotal{F.N.Pd}
 If \(\seqska\) is a sequence of complex \tpqa{matrices}, then, for all \(j\in\mn{0}{\kappa}\), let \(\Pd{j}\defeq\OPu{\ran{\dia{j}}}\) be the matrix corresponding to the orthogonal projection onto \(\ran{\dia{j}}\).
\enota

 The following result is a slight modification of~\zitaa{MR3979701}{\cprop{6.40}{32}}:
 
\bpropl{F.P.abZe}
 Let \(\seqska\in\Fggqka\), assume \(\kappa\geq1\), and let \(k\in\mn{1}{\kappa}\).
 Then the following statements are equivalent:
\baeqi{0}
 \il{F.P.abZe.i} \(\seqska\) is \tabZo{k}.
 \il{F.P.abZe.ii} \(\cia{j}=\frac{1}{2}\Pd{j-1}\) for all \(j\in\mn{k}{\kappa}\).
 \il{F.P.abZe.iii} \(\cia{j}=\frac{1}{2}\Pd{k-1}\) for all \(j\in\mn{k}{\kappa}\).
\eaeqi
 If~\rstat{F.P.abZe.i} is fulfilled, then \(\dia{j}=\rk{\ba/4}^{j-k+1}\dia{k-1}\) and \(\Pd{j}=\Pd{k-1}\) for all \(j\in\mn{k}{\kappa}\), where \(\ba\defeq\obg-\ug\).
\eprop
\bproof
 Suppose that~\rstat{F.P.abZe.iii} is fulfilled.
 In particular, then \(\cia{k}=\frac{1}{2}\Pd{k-1}\).
 Consequently, \(\dia{k}=\rk{\ba/4}\dia{k-1}\), by virtue of~\zitaa{MR3979701}{\cthm{6.34(c)}{31}}.
 In view of \rnota{F.N.Pd} and\(\ba>0\), hence, \(\Pd{k}=\Pd{k-1}\).
 Thus, there exists an integer \(m\in\mn{k}{\kappa}\) such that
 \(\Pd{\ell}=\Pd{k-1}\) for all \(\ell\in\mn{k}{m}\).
 Assume \(m<\kappa\).
 Because of~\rstat{F.P.abZe.iii}, then \(\cia{m+1}=\frac{1}{2}\Pd{k-1}=\frac{1}{2}\Pd{m}\).
 Consequently, \(\dia{m +1}=\rk{\ba/4}\dia{m }\), by virtue of~\zitaa{MR3979701}{\cthm{6.34(c)}{31}}.
 In view of \(\ba>0\), hence, \(\Pd{m +1}=\Pd{m }=\Pd{k-1}\).
 Thus, by mathematical induction, we have proved \(\Pd{j}=\Pd{k-1}\) for all \(j\in\mn{k}{\kappa}\).
 In particular,~\rstat{F.P.abZe.ii} holds true.
 The application of~\zitaa{MR3979701}{\cprop{6.40}{32}} completes the proof.
\eproof

\bexaml{F.E.Z1}
 Let \(\seqsinf\in\Fggqinf\) be \tabZo{1} and let \(\ba\defeq\obg-\ug\).
 According to \rprop{F.P.abZe}, then \(\cia{j}=\frac{1}{2}\Pd{0}\) and \(\dia{j}=\rk{\ba/4}^j\dia{0}\) for all \(j\in\N\).
 From \eqref{F.G.d01} and \eqref{F.G.f012}, we obtain \(\dia{j}=\ba\rk{\ba/4}^j\su{0}\) for all \(j\in\NO\) and, furthermore, \(\cia{0}=\su{0}\) and \(\cia{j}=\frac{1}{2}\OPu{\ran{\su{0}}}\) for all \(j\in\N\).
\eexam

 On the contrary, the interval \(\matint{\umg{m}}{\omg{m}}\) in \rthm{165.T112} shrinks to a single point if \(\dia{m}=\Oqq\), which is equivalent to have trivial intersection \(\ran{\usc{m}}\cap\ran{\osc{m}}=\set{\Ouu{q}{1}}\) for the column spaces of the distances \(\usc{m}=\su{m}-\umg{m-1}\) and \(\osc{m}=\omg{m-1}-\su{m}\) of \(\su{m}\in\matint{\umg{m-1}}{\omg{m-1}}\) to the preceding left and right interval endpoint (\tcf{}~\zitaa{MR3775449}{\clem{10.26}{40}}):
 
\bdefnnl{\tcf{}~\zitaa{MR3775449}{\cdefnssp{10.24}{10.25}{40}}}{F.D.abCD}
 Let \(\seqska\) be a sequence of complex \tpqa{matrices} and let \(k\in\mn{0}{\kappa}\).
 Then \(\seqska\) is said to be \emph{\tabCDo{k}} if \(\dia{k}=\Opq\).
\edefn

\section{The \hFTion{} for sequences of complex matrices}\label{F.S.FT}
 In this section, we consider a particular transformation of sequences of matrices, which will be the elementary step of a particular algorithm.
 First we introduce a particular matrix polynomial and list some simple observations on it.
 
\bnotal{H.N.T}
 For each \(n\in\NO\), let \(\Tq{n}\defeq\matauuo{\Kronu{j,k+1}\Iq}{j,k}{0}{n}\).
\enota
 
 In particular, we have \(\Tq{0}=\Oqq\) and \(\Tq{n}=\smat{\Ouu{q}{nq}&\Oqq\\ \Iu{nq}&\Ouu{nq}{q}}\) for all \(n\in\N\).
 For each \(n\in\NO\), the block matrix \(\Tq{n}\) has the shape of the first of the two matrices in \rnota{M.N.S} with the matrix \(\Oqq\) in the \tqqa{block} main diagonal.
 Thus, \(\det\rk{\Iu{(n+1)q}-z\Tq{n}}=1\) holds true for all \(z\in\C\).
 Hence, the following matrix-valued function is well defined:
 
\bnotal{H.N.Res}
 For all \(n\in\NO\), let \(\Rq{n}\colon\C\to\Coo{(n+1)q}{(n+1)q}\) be given by \[\Rqa{n}{z}\defeq\rk{\Iu{(n+1)q}-z\Tq{n}}^\inv.\]
\enota

 By virtue of \(\Tq{n}^{n+1}=\Ouu{\rk{n+1}q}{\rk{n+1}q}\), it is readily checked that \(\Rqa{n}{z}=\sum_{\ell=0}^nz^\ell\Tq{n}^\ell\) holds true for all \(z\in\C\).
 In particular, \(\Rq{n}\) is a matrix polynomial of degree \(n\) with invertible values, admitting the following \tbr{}: 
 
\bexaml{H.E.ResS}
 Let \(n\in\NO\), let \(z\in\C\), and let the sequence \(\seqs{n}\) be given by \(\su{j}\defeq z^j\Iq\).
 Using \rnota{M.N.S}, then \(\SLu{n}=\Rqa{n}{z}\) and \(\SUu{n}=\ek{\Rqa{n}{\ko z}}^\ad\).
\eexam
 
 Referring to the classes of particular block triangular matrices introduced in \rnota{ab.N1308}, we can conclude from \rexam{H.E.ResS}:
 
\breml{H.R.RinLU}
 Let \(n\in\NO\) and let \(z\in\C\).
 In view of \rnota{ab.N1308}, then \(\Rqa{n}{z}\) is a block Toeplitz matrix belonging to \(\nudpu{n}\) and \(\ek{\Rqa{n}{\ko z}}^\ad\) is a block Toeplitz matrix belonging to \(\nodqu{n}\).
 In particular, \(\det\Rqa{n}{z}=1\), by virtue of \rrem{R1553}.
\erem
  
 Given a sequence \(\seqska\) of complex \tpqa{matrices} and a bounded interval \(\ab\) of \(\R\) we have introduced in \rnota{F.N.sa} three particular sequences associated with \(\seqska\).
 These three sequences correspond to the intervals \(\rhl\), \(\lhl\), and \(\ab\), respectively.
 Now we slightly modify each of these sequences.
 We start with the case of the interval \(\rhl\).
 
\bdefnl{114.D1455}
 Let \(\seqska\) be a sequence of complex \tpqa{matrices}.
 Further let \(\sau{-1}\defeq\su{0}\) and, in the case \(\kappa\geq1\), let \(\seqsa{\kappa-1}\) be given by \rnota{F.N.sa}.
 Then we call the sequence \(\seqspaka\) given by \(\spa{j}\defeq\sau{j-1}\) the \emph{\tamodv{\seqska }}.
 For each matrix \(X_k=X_k^{\ok{s}}\) built from the sequence \(\seqska\), we denote by \(X_\aur{k}'\defeq X_k^{\ok{\spas}}\) the corresponding matrix built from the sequence \(\seqspaka\) instead of \(\seqska\).
\edefn

 In the classical case \(\ug=0\), the sequences \(\seqska \) and \(\seqspaka\) coincide.
 For an arbitrary \(\ug\in\R\), the sequence \(\seqska\) is reconstructible from \(\seqspaka\) as well.
 The corresponding relations can be written in a convenient form, using the particular block Toeplitz matrices introduced in \rnotass{M.N.S}{H.N.Res}:

\breml{ab.R1115a}
 If \(\seqska\) is a sequence of complex \tpqa{matrices}, then \(\SLpa{m}=\ek{\Rpa{m}{\ug}}^\inv\SLu{m}\) and \(\SUpa{m}=\SUu{m}\ek{\Rqa{m}{\ug}}^\invad\) for all \(m\in\mn{0}{\kappa}\).
\erem
 
 Using~\zitaa{MR3611479}{\clem{4.6}{157}}, we obtain the analogous relation between the corresponding \trF{s} associated via \rdefn{D1419}:

\blemnl{\tcf{}~\zitaa{MR3611479}{\crem{4.7}{158}}}{ab.R1835a}
 Let \(\seqska \) be a sequence of complex \tpqa{matrices}.
 Denote by \(\seqr{\kappa}\) the \trFa{\(\seqspaka\)}.
 For all \(m\in\mn{0}{\kappa}\), then \(\SLuo{m}{r}=\Rqa{m}{\ug}\SLu{m}^\rez\) and \(\SUuo{m}{r}=\SUu{m}^\rez\ek{\Rpa{m}{\ug}}^\ad\).
\elem

 The analogue of \rdefn{114.D1455} for the interval \(\lhl\) looks as follows:

\bdefnl{ab.N1137b}
 Let \(\seqska\) be a sequence of complex \tpqa{matrices}.
 Further let \(\sub{-1}\defeq-\su{0}\) and, in the case \(\kappa\geq1\), let \(\seqsb{\kappa-1}\) be given by \rnota{F.N.sa}.
 Then we call the sequence \(\seqspbka\) given by \(\spb{j}\defeq\sub{j-1}\) the \emph{\tbmodv{\seqska}}.
 For each matrix \(X_k=X_k^{\ok{s}}\) built from the sequence \(\seqska\), we denote by \(X_\lub{k}'\defeq X_k^{\ok{\spbs}}\) the corresponding matrix built from the sequence \(\seqspbka\) instead of \(\seqska\).
\edefn

 In particular, if \(\obg=0\), then \(\seqspbka\) coincides with the sequence \(\seq{-\su{j}}{j}{0}{\kappa}\).
 For an arbitrary \(\obg\in\R\), the sequence \(\seqska\) is reconstructible from \(\seqspbka\) as well.
 The corresponding relations can be written in a convenient form using the particular block Toeplitz matrices introduced in \rnotass{M.N.S}{H.N.Res}:
 
\breml{ab.R1115b}
 If \(\seqska \) is a sequence of complex \tpqa{matrices}, then \(\SLpb{m}=-\ek{\Rpa{m}{\obg}}^\inv\SLu{m}\) and \(\SUpb{m}=-\SUu{m}\ek{\Rqa{m}{\obg}}^\invad\) for all \(m\in\mn{0}{\kappa}\).
\erem

 By virtue of \rrem{H.R.rezhom}, we can conclude from \rlem{ab.R1835a} the analogous relation between the \trF{s} associated via \rdefn{D1419} to \(\seqska\) and to \(\seqspbka\):
 
\bleml{ab.R1835b}
 Let \(\seqska \) be a sequence of complex \tpqa{matrices}.
 Denote by \(\seqr{\kappa}\) the \trFa{\(\seqspbka\)}.
 Then \(\SLuo{m}{r}=-\Rqa{m}{\obg}\SLu{m}^\rez\) and \(\SUuo{m}{r}=-\SUu{m}^\rez\ek{\Rpa{m}{\obg}}^\ad\) for all \(m\in\mn{0}{\kappa}\).
\elem

 In addition to the \tamod{} \(\seqspaka\) and the \tbmod{} \(\seqspbka\), we introduce, in view of \rnota{F.N.sa}, the corresponding construction associated to the sequence \(\seqsab{\kappa-2}\):
 
\bdefnl{ab.N1137c}
 Let \(\seqska\) be a sequence of complex \tpqa{matrices} and let \(\sab{-2}\defeq-\su{0}\).
 In the case \(\kappa\geq1\), let \(\sab{-1}\defeq(\ug+\obg)\su{0}-\su{1}\).
 Then we call the sequence \(\seqspcka \) given by \(\spc{j}\defeq\sab{j-2}\) the \emph{\tabmodv{\seqska}}.
 For each matrix \(X_k=X_k^{\ok{s}}\) built from the sequence \(\seqska\), we denote by \(X_\aub{k}'\defeq X_k^{\ok{\spcs}}\) the corresponding matrix built from the sequence \(\seqspcka \) instead of \(\seqska\).
\edefn

\breml{R0801}
 It is readily checked that the sequence \(\seqspcka \) coincides with the \tamodv{\seqspbka} as well as with the \tbmodv{\seqspaka}.
\erem

 Using \rremss{ab.R1115b}{ab.R1115a} and \rlemss{ab.R1835b}{ab.R1835a}, \tresp{}, we obtain in particular:
 
\breml{ab.R1115c}
 If \(\seqska \) is a sequence of complex \tpqa{matrices}, then \(\SLpc{m}=-\ek{\Rpa{m}{\ug}}^\inv\ek{\Rpa{m}{\obg}}^\inv\SLu{m}\) and \(\SUpc{m}=-\SUu{m}\ek{\Rqa{m}{\obg}}^\invad\ek{\Rqa{m}{\ug}}^\invad\) for all \(m\in\mn{0}{\kappa}\).
\erem

 Using the \tCP{} given in \rnota{H.D.CP} and the \trFa{the \tbmodv{\seqsa{\kappa-1}}}, we introduce now a transformation of sequences of complex matrices, which in the context of the moment problem on the interval \(\ab\) turns out to play the same role as the \tHT{} for the moment problem on \(\R\):
 
\bdefnl{ab.N0940}
 Suppose \(\kappa\geq1\).
 Let \(\seqska\) be a sequence of complex \tpqa{matrices}.
 Denote by \(\seqapb{\kappa-1}\) the \tbmodv{\seqsa{\kappa-1}} and by \(\seq{x_j}{j}{0}{\kappa-1}\) the \tCPa{\(\seqsb{\kappa-1}\)}{\(\seq{\apb{j}^\rez}{j}{0}{\kappa-1}\)}.
 Then we call the sequence \(\seqt{\kappa-1}\) given by \(\tu{j}\defeq-\sau{0}\su{0}^\mpi x_j\sau{0}\) the \emph{\tFTv{\(\seqska\)}}.
\edefn

 Since, in the classical case \(\ug=0\) and \(\obg=1\), the sequence \(\seqsa{\kappa-1}\) coincides with the shifted sequence \(\seq{\su{j+1}}{j}{0}{\kappa-1}\), the \tFuuT{0}{1} is given by \(\tu{j}=-\su{1}\su{0}^\mpi x_j\su{1}\) with the \tCP{} \(\seq{x_j}{j}{0}{\kappa-1}\) of \(\seqsb{\kappa-1}\) and \(\seq{\apb{j}^\rez}{j}{0}{\kappa-1}\), where the sequence \(\seqsb{\kappa-1}\) is given by \(\sub{j}=\su{j}-\su{j+1}\) and the sequence \(\seqapb{\kappa-1}\) is given by \(\apb{0}=-\su{1}\) and by \(\apb{j}=\su{j}-\su{j+1}\) for \(j\in\mn{1}{\kappa-1}\).
 
\breml{F.R.FTtr}
 Assume \(\kappa\geq1\) and let \(\seqska\) be a sequence of complex \tpqa{matrices} with \tFT{} \(\seqt{\kappa-1}\).
 Then one can see from \rrem{M.R.reztr} that for each \(k\in\mn{0}{\kappa-1}\), the matrix \(t_k\) is built from the matrices \(\su{0},\su{1},\dotsc,\su{k+1}\).
 In particular, for all \(m\in\mn{1}{\kappa}\), the \tFTv{\(\seqs{m}\)} coincides with \(\seqt{m-1}\).
\erem

 The main goal of this section can be described as follows:
 Let \(\kappa\geq1\) and let \(\seqska\in\Fggqka\).
 Denote by \(\seqa{t}{\kappa-1}\) the \tFTv{\(\seqska\)}.
 Then we are going to prove that \(\seqa{t}{\kappa-1}\in\Fggqu{\kappa-1}\).
 Our main strategy to realize this aim is based on using appropriate identities for \tbHms{} associated with the sequences \(\seqska\), \(\seqa{t}{\kappa-1}\) and also the three sequences \(\seqsa{\kappa-1}\), \(\seqsb{\kappa-1}\) and \(\seqsab{\kappa-2}\) built from \(\seqska\) via \rnota{F.N.sa}.
 Clearly, we will concentrate on the finite sections of the sequences \(\seqska\).
 There will be two different cases, namely the sections \(\seqs{2n}\) on the one hand and the sections \(\seqs{2n+1}\) on the other hand.
 We start with a collection of identities which will be later used in many proofs.

\bleml{R1433}
 Let \(n\in\NO\) and let \(z,w\in\C\).
 Using \rnota{H.N.Res}, then
\begin{align}
 \ek*{\Rpa{n}{z}}^\inv\ek*{\Rpa{n}{w}}^\inv&=\ek*{\Rpa{n}{w}}^\inv\ek*{\Rpa{n}{z}}^\inv,\label{SR4}\\
 \ek*{\Rpa{n}{z}}^\inv\Rpa{n}{w}&=\Rpa{n}{w}\ek*{\Rpa{n}{z}}^\inv,\label{SR3}\\
 \Rpa{n}{z}\Rpa{n}{w}&=\Rpa{n}{w}\Rpa{n}{z}.\label{SR5}
\end{align}
\elem
\bproof
 The identity \eqref{SR4} is an immediate consequence of \rnota{H.N.Res}.
 Formula \eqref{SR3} follows from \eqref{SR4} and implies \eqref{SR5}. 
\eproof

\breml{H.R.updoT}
 Let \(\ell,m\in\N\).
 Using \rnotass{H.N.updo}{H.N.T}, direct computations give us:
\benui
 \il{H.R.updoT.a} \(\IOquu{\ell}{m}\IOquu{\ell}{m}^\ad=\zdiag{\Iu{mq}}{\Ouu{\ell}{\ell}}\) and \(\OIquu{\ell}{m}\OIquu{\ell}{m}^\ad=\zdiag{\Ouu{\ell}{\ell}}{\Iu{mq}}\).
 \il{H.R.updoT.b} \(\IOquu{\ell}{m}^\ad\IOquu{\ell}{m}=\Iu{mq}\) and \(\OIquu{\ell}{m}^\ad\OIquu{\ell}{m}=\Iu{mq}\).
 \il{H.R.updoT.c} If \(n\in\N\), then \(\IOquu{1}{n}^\ad\OIquu{1}{n}=\Tq{n-1}\) and \(\OIquu{1}{n}\IOquu{1}{n}^\ad=\Tq{n}\).
\eenui
\erem
 
\bleml{L1301}
 Let \(\seqska\) be a sequence of complex \tpqa{matrices}.
\benui
 \il{L1301.a} Let \(m\in\mn{1}{\kappa}\).
 In view of \rnotass{M.N.S}{H.N.updo}, then
\begin{align}
 \vpu{m}^\ad\SLu{m}&=\su{0}\vqu{m}^\ad,&&&
 \SUu{m}\vqu{m}&=\vpu{m}\su{0},\label{SD1}\\
 \IOpuu{1}{m}^\ad\SLu{m}&=\SLu{m-1}\IOquu{1}{m}^\ad,&&&
 \SUu{m}\IOquu{1}{m}&=\IOpuu{1}{m}\SUu{m-1},\label{SD2}\\
 \SLu{m}\OIquu{1}{m}&=\OIpuu{1}{m}\SLu{m-1},&&\text{and}&
 \OIpuu{1}{m}^\ad\SUu{m}&=\SUu{m-1}\OIquu{1}{m}^\ad.\label{SD3}
\end{align}
 \il{L1301.b} Let \(n\in\mn{0}{\kappa}\).
 In view of \rnota{H.N.T}, then \(\Tp{n}\SLu{n}=\SLu{n}\Tq{n}\) and \(\SUu{n}\Tq{n}^\ad=\Tp{n}^\ad\SUu{n}\).
 \il{L1301.c} Let \(n\in\mn{0}{\kappa}\) and let \(z\in\C\).
 In view of \rnota{H.N.Res}, then
\begin{align}
 \ek*{\Rpa{n}{z}}^\inv\SLu{n}&=\SLu{n}\ek*{\Rqa{n}{z}}^\inv,&&&
 \SUu{n}\ek*{\Rqa{n}{z}}^\invad&=\ek*{\Rpa{n}{z}}^\invad\SUu{n},\label{SR1}\\
 \Rpa{n}{z}\SLu{n}&=\SLu{n}\Rqa{n}{z},&&\text{and}&\SUu{n}\ek*{\Rqa{n}{z}}^\ad&=\ek*{\Rpa{n}{z}}^\ad\SUu{n}.\label{SR2}
\end{align}
\eenui
\elem
\bproof
 \eqref{L1301.a} The identities \eqref{SD1}--\eqref{SD3} follow from \rnota{H.N.updo} and \rrem{H.R.Sblock}.
 
 \eqref{L1301.b} In view of \(\Tp{0}=\Opp\) and  \(\Tq{0}=\Oqq\), the case \(n=0\) is trivial.
 Suppose \(n\geq1\).
 Using \rremp{H.R.updoT}{H.R.updoT.c}, \eqref{SD2}, and \eqref{SD3}, we obtain
\[
 \Tp{n}\SLu{n}
 =\OIpuu{1}{n}\IOpuu{1}{n}^\ad\SLu{n}
 =\OIpuu{1}{n}\SLu{n-1}\IOquu{1}{n}^\ad
 =\SLu{n}\OIquu{1}{n}\IOquu{1}{n}^\ad
 =\SLu{n}\Tq{n}
\]
 and
\[
 \SUu{n}\Tq{n}^\ad
 =\SUu{n}\IOquu{1}{n}\OIquu{1}{n}^\ad
 =\IOpuu{1}{n}\SUu{n-1}\OIquu{1}{n}^\ad
 =\IOpuu{1}{n}\OIpuu{1}{n}^\ad\SUu{n}
 =\Tp{n}^\ad\SUu{n}.
\]

 \eqref{L1301.c} The identities \eqref{SR1} follow from \rnota{H.N.Res} and \eqref{L1301.b}.
 The identities \eqref{SR2} follow from \eqref{SR1}.
\eproof

\breml{R1317}
 Let \(A\) be a complex \tpqa{matrix}.
 Then \rlem{L1301} yields:
\benui
 \il{R1317.a} Let \(m\in\N\).
 Denote by \(\sdiag{A}{m}\) the block diagonal matrix built via \eqref{sdiag} from \(A\).
 In view of \rnotass{H.N.updo}{H.N.T}, then
\begin{align}
 \vpu{m}^\ad\sdiag{A}{m}&=A\vqu{m}^\ad,&&&
 \sdiag{A}{m}\vqu{m}&=\vpu{m}A,\label{DA1}\\
 \IOpuu{1}{m}^\ad\sdiag{A}{m}&=\sdiag{A}{m-1}\IOquu{1}{m}^\ad,&&&
 \sdiag{A}{m}\IOquu{1}{m}&=\IOpuu{1}{m}\sdiag{A}{m-1},\notag\\
 \sdiag{A}{m}\OIquu{1}{m}&=\OIpuu{1}{m}\sdiag{A}{m-1},&&\text{and}&
 \OIpuu{1}{m}^\ad\sdiag{A}{m}&=\sdiag{A}{m-1}\OIquu{1}{m}^\ad.\label{DA3}
\end{align}
 \il{R1317.b} Let \(n\in\NO\).
 Denote by \(\sdiag{A}{n}\) the block diagonal matrix built via \eqref{sdiag} from \(A\).
 In view of \rnota{H.N.T}, then \(\Tp{n}\sdiag{A}{n}=\sdiag{A}{n}\Tq{n}\) and \(\sdiag{A}{n}\Tq{n}^\ad=\Tp{n}^\ad\sdiag{A}{n}\).
 
 \il{R1317.c} Let \(n\in\NO\) and let \(z\in\C\).
 Denote by \(\sdiag{A}{n}\) the block diagonal matrix built via \eqref{sdiag} from \(A\).
 In view of \rnota{H.N.Res}, then
\begin{align}
 \ek*{\Rpa{n}{z}}^\inv\sdiag{A}{n}&=\sdiag{A}{n}\ek*{\Rqa{n}{z}}^\inv,&&&\sdiag{A}{n}\ek*{\Rqa{n}{z}}^\invad&=\ek*{\Rpa{n}{z}}^\invad\sdiag{A}{n},\label{RA1}\\
 \Rpa{n}{z}\sdiag{A}{n}&=\sdiag{A}{n}\Rqa{n}{z},&&\text{and}&
 \sdiag{A}{n}\ek*{\Rqa{n}{z}}^\ad&=\ek*{\Rpa{n}{z}}^\ad\sdiag{A}{n}.\label{RA2}
\end{align}
\eenui
\erem

\breml{R1326}
 Let \(A\) be a complex \tpqa{matrix} and let \(n\in\NO\).
 Denote by \(\sdiag{A}{n}\) the block diagonal matrix built via \eqref{sdiag} from \(A\).
\benui
 \il{R1326.a} It holds \(\rk{\sdiag{A}{n}}^\ad=\sdiag{A^\ad}{n}\) and, by virtue of \rrem{A.R.diag+}, moreover \(\rk{\sdiag{A}{n}}^\mpi=\sdiag{A^\mpi}{n}\).
 \il{R1326.b} Let \(r\in\N\), let \(B\in\Coo{q}{r}\), and let \(\lambda\in\C\).
 Then \(\sdiag{\lambda A}{n}=\lambda\sdiag{A}{n}\) and \(\sdiag{AB}{n}=\sdiag{A}{n}\sdiag{B}{n}\).
\eenui
\erem
 
 Now we turn our attention to the \([\ug,\infp)\)- and \tbmod{} and to the \tabmod{} of the \tFT{} and the interplay of these sequences and their \trF{s} with forming \tCP{s}:

\bleml{ab.R1110}
 Let \(\seqska\) be a sequence of complex \tpqa{matrices} and assume \(\kappa\geq1\).
 Let the sequences \(\seqsa{\kappa-1}\) and \(\seqsb{\kappa-1}\) be given by \rnota{F.N.sa}.
 Denote by \(\seqt{\kappa-1}\) the \tFTv{\(\seqska\)} and by \(\seqbpa{\kappa-1}\) the \tamodv{\seqsb{\kappa-1}}.
 \benui
  \il{ab.R1110.a} Denote by \(\seqapb{\kappa-1}\) the \tbmodv{\seqsa{\kappa-1}}, by \(\seq{\mathbf{x}_j}{j}{0}{\kappa-1}\) the \tCPa{\(\seqbpa{\kappa-1}\)}{\(\seq{\apb{j}^\rez}{j}{0}{\kappa-1}\)}, and by \(\seqtpa{\kappa-1}\) the \tamodv{\seqt{\kappa-1}}.
  For all \(j\in\mn{0}{\kappa-1}\), then \(\tpa{j}=-\sau{0}\su{0}^\mpi \mathbf{x}_j\sau{0}\).
  \il{ab.R1110.b} Denote by \(\seq{\mathbf{y}_j}{j}{0}{\kappa-1}\) the \tCPa{\(\seqsb{\kappa-1}\)}{\(\seq{\sau{j}^\rez}{j}{0}{\kappa-1}\)} and by \(\seqtpb{\kappa-1}\) the \tbmodv{\seqt{\kappa-1}}.
  For all \(j\in\mn{0}{\kappa-1}\), then \(\tpb{j}=-\sau{0}\su{0}^\mpi \mathbf{y}_j\sau{0}\).
  \il{ab.R1110.c} Denote by \(\seq{\mathbf{z}_j}{j}{0}{\kappa-1}\) the \tCPa{\(\seqbpa{\kappa-1}\)}{\(\seq{\sau{j}^\rez}{j}{0}{\kappa-1}\)} and by \(\seqtpc{\kappa-1}\) the \tabmodv{\seqt{\kappa-1}}.
  For all \(j\in\mn{0}{\kappa-1}\), then \(\tpc{j}=-\sau{0}\su{0}^\mpi \mathbf{z}_j\sau{0}\).
 \eenui
\elem
\bproof
 We consider an arbitrary \(k\in\mn{0}{\kappa-1}\).
 Denote by \(\seq{\mathbf{h}_j}{j}{0}{\kappa-1}\) the \trFa{\(\seqapb{\kappa-1}\)}.
 In view of \rdefn{ab.N0940} and \rrem{M.R.s*t}, then
 \beql{ab.R1110.1}
  \SLuo{k}{t}
  =\sdiag{-\sau{0}\su{0}^\mpi}{k}\SLuo{k}{b}\SLuo{k}{\mathbf{h}}\sdiag{\sau{0}}{k}.
 \eeq
 According to \rrem{ab.R1115a}, we have
 \beql{ab.R1110.3}
  \ek*{\Rpa{k}{\ug}}^\inv\SLuo{k}{b}
  =\SLuo{k}{\bpas}.
 \eeq
 Denote by \(\seq{r_j}{j}{0}{\kappa-1}\) the \trFa{\(\seqsa{\kappa-1}\)}.
 Then \rlem{ab.R1835b} yields
 \beql{ab.R1110.2}
  -\ek*{\Rqa{k}{\obg}}^\inv\SLuo{k}{\mathbf{h}}
  =\SLuo{k}{r}.
 \eeq
 
 \eqref{ab.R1110.a} According to \rrem{ab.R1115a}, we have \(\SLuo{k}{\tpas}=\ek{\Rpa{k}{\ug}}^\inv\SLuo{k}{t}\).
 By virtue of \eqref{ab.R1110.1}, \eqref{RA1}, \eqref{ab.R1110.3}, and \rrem{M.R.s*t}, then
\[
  \SLuo{k}{\tpas}
  =\ek*{\Rpa{k}{\ug}}^\inv\sdiag{-\sau{0}\su{0}^\mpi}{k}\SLuo{k}{b}\SLuo{k}{\mathbf{h}}\sdiag{\sau{0}}{k}
  =\sdiag{-\sau{0}\su{0}^\mpi}{k}\SLuo{k}{\bpas}\SLuo{k}{\mathbf{h}}\sdiag{\sau{0}}{k}
  =\sdiag{-\sau{0}\su{0}^\mpi}{k}\SLuo{k}{\mathbf{x}}\sdiag{\sau{0}}{k}.
\]
 follows.
 Consequently, \(\tpa{k}=-\sau{0}\su{0}^\mpi \mathbf{x}_k\sau{0}\) holds true.
 
 \eqref{ab.R1110.b} According to \rrem{ab.R1115b}, we have \(\SLuo{k}{\tpbs}=-\ek{\Rpa{k}{\obg}}^\inv\SLuo{k}{t}\).
 By virtue of \eqref{ab.R1110.1}, \eqref{RA1}, \eqref{SR1}, \eqref{ab.R1110.2}, and \rrem{M.R.s*t}, we get then
\[
  \SLuo{k}{\tpbs}
  =-\ek{\Rpa{k}{\obg}}^\inv\sdiag{-\sau{0}\su{0}^\mpi}{k}\SLuo{k}{b}\SLuo{k}{\mathbf{h}}\sdiag{\sau{0}}{k}
  =\sdiag{-\sau{0}\su{0}^\mpi}{k}\SLuo{k}{b}\SLuo{k}{r}\sdiag{\sau{0}}{k}
  =\sdiag{-\sau{0}\su{0}^\mpi}{k}\SLuo{k}{\mathbf{y}}\sdiag{\sau{0}}{k}.
\]
 Consequently, \(\tpb{k}=-\sau{0}\su{0}^\mpi \mathbf{y}_k\sau{0}\) holds true.

 \eqref{ab.R1110.c} According to \rrem{ab.R1115c}, we have \(\SLuo{k}{\tpcs}=-\ek{\Rpa{k}{\ug}}^\inv\ek{\Rpa{k}{\obg}}^\inv\SLuo{k}{t}\).
 By virtue of \eqref{ab.R1110.1}, \eqref{RA1}, \eqref{SR1}, \eqref{ab.R1110.3}, \eqref{ab.R1110.2}, and \rrem{M.R.s*t}, then
\begin{multline*}
  \SLuo{k}{\tpcs}
  =-\ek{\Rpa{k}{\ug}}^\inv\ek{\Rpa{k}{\obg}}^\inv\sdiag{-\sau{0}\su{0}^\mpi}{k}\SLuo{k}{b}\SLuo{k}{\mathbf{h}}\sdiag{\sau{0}}{k}\\
  =\sdiag{-\sau{0}\su{0}^\mpi}{k}\SLuo{k}{\bpas}\SLuo{k}{r}\sdiag{\sau{0}}{k}
  =\sdiag{-\sau{0}\su{0}^\mpi}{k}\SLuo{k}{\mathbf{z}}\sdiag{\sau{0}}{k}
\end{multline*}
 follows.
 Consequently, \(\tpc{k}=-\sau{0}\su{0}^\mpi \mathbf{z}_k\sau{0}\) holds true.
\eproof

 In the remaining part of this section, we will derive representations for the \tbHms{} built from the \tFT{} of a sequence via \rnota{N.HKG}, \eqref{F.G.HaHbHK}, and  \eqref{F.G.HabHKG}.
 Therefore, we first supplement \rnota{H.N.D}:

\bnotal{F.N.Dd}
 Let \(\kappa,\tau\in\NOinf\) and let \(\seqska\) and \(\seqt{\tau}\) be two sequences of complex \tpqa{matrices}.
 Denote by \(\seqr{\kappa}\) the \trFa{\(\seqska\)}.
 For each \(m\in\NO\) with \(m\leq\min\set{\kappa,\tau}\), then let
\begin{align*}
 \Dloou{s}{t}{m}
 &\defeq\sdiag{\su{0}}{m}\SLuo{m}{r}\SLuo{m}{t}\sdiag{t_0^\mpi}{m}+\sdiag{\Ip-\su{0}\su{0}^\mpi t_0t_0^\mpi}{m}
\shortintertext{and}
 \Droou{s}{t}{m}
 &\defeq\sdiag{t_0^\mpi}{m}\SUuo{m}{t}\SUuo{m}{r}\sdiag{\su{0}}{m}+\sdiag{\Iq-t_0^\mpi t_0\su{0}^\mpi\su{0}}{m}.
\end{align*}
\enota

 In particular, if \(p=q\) and if the sequence \(\seqt{\tau}\) is given by \(t_j\defeq\Kronu{j0}\Iq\), then \(\Dloou{s}{t}{m}=\Dluo{m}{s}\) and \(\Droou{s}{t}{m}=\Druo{m}{s}\).
 Referring to the classes of particular block triangular matrices introduced in \rnota{ab.N1308}, by virtue of \(r_0=\su{0}^\rez=\su{0}^\mpi\) and \rrem{M.R.s*t}, we have:
 
\breml{F.R.DdinLU}
 Let \(\kappa,\tau\in\NOinf\) and let \(\seqska\) and \(\seqt{\tau}\) be two sequences from \(\Cpq\).
 For all \(m\in\mn{0}{\kappa}\cap\mn{0}{\tau}\), then \(\Dloou{s}{t}{m}\) is a block Toeplitz matrix belonging to \(\nudpu{m}\) and \(\Droou{s}{t}{m}\) is a block Toeplitz matrix belonging to \(\nodqu{m}\).
\erem

\breml{F.R.detDd}
 Let \(\kappa,\tau\in\NOinf\) and let \(\seqska\) and \(\seqt{\tau}\) be two sequences from \(\Cpq\).
 In view of \rremss{F.R.DdinLU}{R1553}, for all \(m\in\mn{0}{\kappa}\cap\mn{0}{\tau}\), then  \(\det\Dloou{s}{t}{m}=1\) and \(\det\Droou{s}{t}{m}=1\).
\erem

 By virtue of \rprop{H.R.rez*}, \rremp{R1326}{R1326.a}, and \rrem{A.R.A++*}, the following two remarks can be verified:

\breml{H.R.D*}
 Let \(\seqska\) be a sequence of complex \tpqa{matrices} and let the sequence \(\seqt{\kappa}\) be given by \(t_j\defeq\su{j}^\ad\).
 For all \(m\in\mn{0}{\kappa}\), then \(\Dlu{m}^\ad=\Druo{m}{t}\).
\erem

\breml{F.R.Dd*}
 Let \(\kappa,\tau\in\NOinf\) and let \(\seqska\) and \(\seqt{\tau}\) be two sequences from \(\Cpq\).
 Let the sequences \(\seq{u_j}{j}{0}{\kappa}\) and \(\seq{v_j}{j}{0}{\tau}\) be given by \(u_j\defeq\su{j}^\ad\) and \(v_j\defeq\tu{j}^\ad\), \tresp{}
 For all \(m\in\mn{0}{\kappa}\cap\mn{0}{\tau}\), then \(\rk{\Dloou{s}{t}{m}}^\ad=\Droou{u}{v}{m}\).
\erem

 Using \rremp{R1326}{R1326.b} and \rrem{R.AA+B}, we obtain the following two remarks:

\breml{ab.R1602}
 Let \(\seqska\) be a sequence of complex \tpqa{matrices} and let \(m\in\mn{0}{\kappa}\).
 For all \(B\in\Coo{p}{u}\) with \(\ran{B}\subseteq\ran{\su{0}}\), then \(\Dlu{m}\sdiag{B}{m}=\sdiag{\su{0}}{m}\SLu{m}^\rez\sdiag{B}{m}\).
 For all \(C\in\Coo{v}{q}\) with \(\nul{\su{0}}\subseteq\nul{C}\), furthermore \(\sdiag{C}{m}\Dru{m}=\sdiag{C}{m}\SUu{m}^\rez\sdiag{\su{0}}{m}\).
\erem

\breml{ab.R1520}
 Let \(\kappa,\tau\in\NOinf\), let \(\seqska\) and \(\seqt{\tau}\) be two sequences of complex \tpqa{matrices}, and let \(m\in\NO\) with \(m\leq\min\set{\kappa,\tau}\).
 Denote by \(\seqr{\kappa}\) the \trFa{\(\seqska\)}.
 For all \(B\in\Coo{p}{u}\) with \(\ran{B}\subseteq\ran{t_0}\cap\ran{\su{0}}\), then \(\Dloou{s}{t}{m}\sdiag{B}{m}=\sdiag{\su{0}}{m}\SLuo{m}{r}\SLuo{m}{t}\sdiag{t_0^\mpi B}{m}\).
 For all \(C\in\Coo{v}{q}\) with \(\nul{t_0}\cup\nul{\su{0}}\subseteq\nul{C}\), furthermore \(\sdiag{C}{m}\Droou{s}{t}{m}=\sdiag{Ct_0^\mpi}{m}\SUuo{m}{t}\SUuo{m}{r}\sdiag{\su{0}}{m}\).
\erem

\bleml{ab.P1728a}
 Suppose \(\kappa\geq1\).
 Let \(\seqska\in\Dpqka\).
 Denote by \(\seq{t_j}{j}{0}{\kappa-1}\) the \tFTv{\(\seqska\)}.
 Then \(t_0=\dia{1}\), where \(\dia{1}\) is given by \rdefn{D1861}.
\elem
\bproof
 Denote by \(\seqapb{\kappa-1}\) the \tbmodv{\seqsa{\kappa-1}}, by \(\seq{\mathbf{h}_j}{j}{0}{\kappa-1}\) the  \trFa{\(\seq{\apb{j}}{j}{0}{\kappa-1}\)}, and by \(\seq{x_j}{j}{0}{\kappa-1}\) the \tCPa{\(\seqsb{\kappa-1}\)}{\(\seq{\mathbf{h}_j}{j}{0}{\kappa-1}\)}.
 Then \(\mathbf{h}_0=\apb{0}^\mpi=\rk{-\sau{0}}^\mpi\).
 Consequently, \(x_0=\sub{0}\mathbf{h}_0=-\sub{0}\sau{0}^\mpi\).
 Taking into account \rrem{ab.L0907}, we get then
 \(
  \tu{0}
  =-\sau{0}\su{0}^\mpi x_0\sau{0}
  =\sau{0}\su{0}^\mpi\sub{0}\sau{0}^\mpi\sau{0}
  =\sub{0}\su{0}^\mpi\sau{0}\sau{0}^\mpi\sau{0}
  =\dia{1}
 \).
\eproof

 From \rlem{ab.P1728a} we obtain for the \tFT{} \(\seqt{\kappa-1}\) of a sequence \(\seqska \) belonging to \(\Dpqka\) in particular \(\Huo{0}{t}=\dia{1}\).
 To obtain in \rprop{F.P.FTH} below, a convenient formula for the \tbHms{} \(\Huo{n}{t}\) for an arbitrary \(n\in\N\) with \(2n+1\leq\kappa\), we need a series of auxiliary results:
 
\bleml{ab.L1746}
 Suppose \(\kappa\geq1\).
 Let \(\seqska \in\Dpqka\) and let the sequence \(\seq{h_j}{j}{0}{\kappa-1}\) be given by \(h_j\defeq\su{j+1}\).
 Then \(\SLu{m}\sdiag{\su{0}^\mpi}{m}\SLuo{m}{h}=\SLuo{m}{h}\sdiag{\su{0}^\mpi}{m}\SLu{m}\) and \(\SUu{m}\sdiag{\su{0}^\mpi}{m}\SUuo{m}{h}=\SUuo{m}{h}\sdiag{\su{0}^\mpi}{m}\SUu{m}\) for all \(m\in\mn{0}{\kappa-1}\).
\elem
\bproof
 Consider an arbitrary \(m\in\mn{0}{\kappa-1}\).
 According to \rnota{M.N.S}, we have \(\SLuo{m}{h}\Tq{m}=\SLu{m}-\sdiag{\su{0}}{m}=\Tp{m}\SLuo{m}{h}\).
 Because of \(\seqska \in\Dpqka\), furthermore \(\ran{h_j}\subseteq\ran{\su{0}}\) and \(\nul{\su{0}}\subseteq\nul{h_j}\) hold true for all \(j\in\mn{0}{\kappa-1}\).
 By virtue of \rrem{R.AA+B}, we get then \(\sdiag{\su{0}\su{0}^\mpi}{m}\SLuo{m}{h}=\SLuo{m}{h}=\SLuo{m}{h}\sdiag{\su{0}^\mpi\su{0}}{m}\).
 Consequently, using additionally \rremp{R1326}{R1326.b} and \rremp{R1317}{R1317.b}, we obtain thus
\[\begin{split}
 \SLu{m}\sdiag{\su{0}^\mpi}{m}\SLuo{m}{h}
 &=\rk*{\SLuo{m}{h}\Tq{m}+\sdiag{\su{0}}{m}}\sdiag{\su{0}^\mpi}{m}\SLuo{m}{h}\\
 &=\SLuo{m}{h}\Tq{m}\sdiag{\su{0}^\mpi}{m}\SLuo{m}{h}+\sdiag{\su{0}\su{0}^\mpi}{m}\SLuo{m}{h}
 =\SLuo{m}{h}\sdiag{\su{0}^\mpi}{m}\Tp{m}\SLuo{m}{h}+\SLuo{m}{h}\sdiag{\su{0}^\mpi\su{0}}{m}\\
 &=\SLuo{m}{h}\sdiag{\su{0}^\mpi}{m}\rk*{\Tp{m}\SLuo{m}{h}+\sdiag{\su{0}}{m}}
 =\SLuo{m}{h}\sdiag{\su{0}^\mpi}{m}\SLu{m}.
\end{split}\]
 The second identity can be checked analogously.
\eproof

 Observe that, with the sequence \(\seq{h_j}{j}{0}{\kappa-1}\) given in \rlem{ab.L1746}, we have \(\SLau{m}=-\ug\SLu{m}+\SLuo{m}{h}\), \(\SUau{m}=-\ug\SUu{m}+\SUuo{m}{h}\), \(\SLub{m}=\obg\SLu{m}-\SLuo{m}{h}\), and \(\SUub{m}=\obg\SUu{m}-\SUuo{m}{h}\)
 for all  \(m\in\mn{0}{\kappa-1}\).
 In view of \rrem{F.R.c=ab}, we obtain:

\breml{F.R.abp=c-1}
 Suppose \(\kappa\geq1\).
 Let \(\seqska\) be a sequence of complex \tpqa{matrices}.
 Denote by \(\seqbpa{\kappa-1}\) the \tamodv{\seqsb{\kappa-1}} and by \(\seqapb{\kappa-1}\) the \tbmodv{\seqsa{\kappa-1}}.
 Then \(\bpa{0}=\sub{0}\) and \(\apb{0}=-\sau{0}\).
 For all \(j\in\mn{1}{\kappa-1}\), furthermore \(\bpa{j}=\sab{j-1}=\apb{j}\).
\erem

\breml{ab.L1118}
 Suppose \(\kappa\geq1\).
 Let \(\seqska\) be a sequence of complex \tpqa{matrices}.
 Denote by \(\seqbpa{\kappa-1}\) the \tamodv{\seqsb{\kappa-1}} and by \(\seqapb{\kappa-1}\) the \tbmodv{\seqsa{\kappa-1}}.
 From \rremss{F.R.abp=c-1}{F.R.012} we can conclude that \(\bpa{0}-\apb{0}=\ba\su{0}\) and \(\bpa{j}-\apb{j}=\Opq\) for all \(j\in\mn{1}{\kappa-1}\).
\erem

\bleml{ab.L1846b}
 Suppose \(\kappa\geq2\).
 Let \(\seqska \) be a sequence of complex \tpqa{matrices} and let \(n\in\mn{0}{\kappa-2}\).
 Denote by \(\seqbpa{\kappa-1}\) the \tamodv{\seqsb{\kappa-1}} and by \(\seqapb{\kappa-1}\) the \tbmodv{\seqsa{\kappa-1}}.
 For all \(k\in\mn{1}{\kappa-n-1}\), then
\beql{ab.L1846b.B}
 \SLau{n}\sdiag{\su{0}^\mpi}{n}\yuuo{k}{k+n}{\bpas}-\SLuo{n}{\bpas}\sdiag{\su{0}^\mpi}{n}\yauu{k}{k+n}
 =\obg\yauu{0}{n}\su{0}^\mpi\sau{k-1}-\ba\sdiag{\su{0}\su{0}^\mpi}{n}\yauu{k}{k+n}
\eeq
 and
\beql{ab.L1846b.A}
 \SLub{n}\sdiag{\su{0}^\mpi}{n}\yuuo{k}{k+n}{\apbs}-\SLuo{n}{\apbs}\sdiag{\su{0}^\mpi}{n}\yuub{k}{k+n}
 =\ba\sdiag{\su{0}\su{0}^\mpi}{n}\yuub{k}{k+n}-\ug\yuub{0}{n}\su{0}^\mpi\sub{k-1}.
\eeq
\elem
\bproof
 By way of example, we only show \eqref{ab.L1846b.A}.
 First consider an arbitrary \(k\in\mn{1}{\kappa-1}\).
 In view of \rremss{F.R.abp=c-1}{F.R.c=ab}, then \(\bpa{0}=\sub{0}\) and \(\apb{k}=\sab{k-1}=-\ug\sub{k-1}+\sub{k}\).
 In particular,
\(
 \sub{0}\su{0}^\mpi \apb{k}
 =-\ug\sub{0}\su{0}^\mpi\sub{k-1}+\bpa{0}\su{0}^\mpi\sub{k}
\).
 Using \rrem{ab.L1118}, consequently,
\(
 \sub{0}\su{0}^\mpi \apb{k}-\apb{0}\su{0}^\mpi\sub{k}
 =-\ug\sub{0}\su{0}^\mpi\sub{k-1}+\ba\su{0}\su{0}^\mpi\sub{k},
\)
 implying \eqref{ab.L1846b.A} in the case \(n=0\).

 Now assume \(\kappa\geq3\) and \(n\geq1\).
 Consider an arbitrary \(k\in\mn{1}{\kappa-n-1}\).
 From \rrem{ab.L1118} \rremp{R1326}{R1326.b} we infer \(\SLuo{n}{\bpas}-\SLuo{n}{\apbs}=\ba\sdiag{\su{0}}{n}\).
 Hence, taking into account \rremp{R1326}{R1326.b}, we get
\[
 \SLuo{n}{\apbs}\sdiag{\su{0}^\mpi}{n}\yuub{k}{k+n}
 =\SLuo{n}{\bpas}\sdiag{\su{0}^\mpi}{n}\yuub{k}{k+n}-\ba\sdiag{\su{0}\su{0}^\mpi}{n}\yuub{k}{k+n}.
\]
 The application of \rrem{ab.R1115a} to the sequence \(\seqsb{\kappa-1}\) yields furthermore \(\SLuo{n}{\bpas}=\ek{\Rpa{n}{\ug}}^\inv\SLub{n}\), which by virtue of \eqref{SR1} and \eqref{RA1} implies
\[
 \SLuo{n}{\bpas}\sdiag{\su{0}^\mpi}{n}\yuub{k}{k+n}
 =\SLub{n}\sdiag{\su{0}^\mpi}{n}\ek*{\Rpa{n}{\ug}}^\inv\yuub{k}{k+n}.
\]
 Thus, we can conclude
\begin{multline}\label{ab.L1846b.1}
 \SLub{n}\sdiag{\su{0}^\mpi}{n}\yuuo{k}{k+n}{\apbs}-\SLuo{n}{\apbs}\sdiag{\su{0}^\mpi}{n}\yuub{k}{k+n}\\
 =\SLub{n}\sdiag{\su{0}^\mpi}{n}\rk*{\yuuo{k}{k+n}{\apbs}-\ek*{\Rpa{n}{\ug}}^\inv\yuub{k}{k+n}}+\ba\sdiag{\su{0}\su{0}^\mpi}{n}\yuub{k}{k+n}.
\end{multline}
 Using \rremss{F.R.abp=c-1}{F.R.c=ab}, we obtain
\(
 \yuuo{k}{k+n}{\apbs}
 =\yuuo{k-1}{k+n-1}{c}
 =\smat{
  \sab{k-1}\\
  \yuuo{k}{k+n-1}{c}
 }
\)
 and
\[
 \ek*{\Rpa{n}{\ug}}^\inv\yuub{k}{k+n}
 =
 \bMat
  \sub{k}\\
  \yuub{k+1}{k+n}-\ug\yuub{k}{k+n-1}
 \eMat
 =
 \bMat
  \sab{k-1}+\ug\sub{k-1}\\
  \yuuo{k}{k+n-1}{c}
 \eMat,
\]
 implying
\beql{ab.L1846b.2}
 \yuuo{k}{k+n}{\apbs}-\ek*{\Rpa{n}{\ug}}^\inv\yuub{k}{k+n}
 =-\ug\vpu{n}\sub{k-1}.
\eeq
 By virtue of \rrem{H.R.Sblock}, we infer \(\SLub{n}\vqu{n}=\yuub{0}{n}\).
 Hence, in view of \eqref{DA1}, we get
\[
 \SLub{n}\sdiag{\su{0}^\mpi}{n}\rk{-\ug\vpu{n}\sub{k-1}}
 =-\ug\SLub{n}\vqu{n}\su{0}^\mpi\sub{k-1}
 =-\ug\yuub{0}{n}\su{0}^\mpi\sub{k-1}.
\]
 Taking into account \eqref{ab.L1846b.2} and \eqref{ab.L1846b.1}, then \eqref{ab.L1846b.A} follows.
\eproof

\bleml{ab.L1413b}
 Suppose \(\kappa\geq2\).
 Let \(\seqska\) be a sequence of complex \tpqa{matrices} and let \(n\in\mn{1}{\kappa-1}\).
 Denote by \(\seqapb{\kappa-1}\) the \tbmodv{\seqsa{\kappa-1}}.
 For all \(k\in\mn{0}{\kappa-n-1}\), then
\beql{ab.L1413b.A}
 \SLub{n}\sdiag{\su{0}^\mpi}{n}\OIpuu{1}{n}\yuuo{k+1}{k+n}{\apbs}-\SLuo{n}{\apbs}\sdiag{\su{0}^\mpi}{n}\yuub{k}{k+n}
 =
 \bMat
  \sau{0}\su{0}^\mpi\sub{k}\\
  \ba\sdiag{\su{0}\su{0}^\mpi}{n}\yuub{k+1}{k+n}-\yuub{1}{n}\su{0}^\mpi\sub{k}
 \eMat.
\eeq
\elem
\bproof
 Let \(k\in\mn{0}{\kappa-n-1}\).
 Denote by \(X\) the upper \tpqa{block} and by \(Y\) the lower \taaa{np}{q}{block} of the matrix on
 the left-hand side of \eqref{ab.L1413b.A}.
 In view of \rrem{H.R.Sblock}, we get
\[
 \SLub{n}\sdiag{\su{0}^\mpi}{n}\OIpuu{1}{n}\yuuo{k+1}{k+n}{\apbs}
 =\SLub{n}\bMat\Ouu{q}{np}\\ \sdiag{\su{0}^\mpi}{n-1}\eMat\yuuo{k+1}{k+n}{\apbs}
 =\bMat\Ouu{p}{nq}\\ \SLub{n-1}\sdiag{\su{0}^\mpi}{n-1}\yuuo{k+1}{k+n}{\apbs}\eMat
\]
 and
\[\begin{split}
 \SLuo{n}{\apbs}\sdiag{\su{0}^\mpi}{n}\yuub{k}{k+n}
 &=
 \bMat
  \apb{0} & \Ouu{p}{nq} \\
  \yuuo{1}{n}{\apbs} & \SLuo{n-1}{\apbs}
 \eMat\sdiag{\su{0}^\mpi}{n}
 \bMat
  \sub{k}\\
  \yuub{k+1}{k+n}
 \eMat\\
 &=
 \bMat
  \apb{0}\su{0}^\mpi\sub{k}\\
  \yuuo{1}{n}{\apbs}\su{0}^\mpi\sub{k}+\SLuo{n-1}{\apbs}\sdiag{\su{0}^\mpi}{n-1}\yuub{k+1}{k+n}
 \eMat.
\end{split}\]
 Consequently,
\begin{align*}
 X&=-\apb{0}\su{0}^\mpi\sub{k}&
&\text{and}&
 Y&=\SLub{n-1}\sdiag{\su{0}^\mpi}{n-1}\yuuo{k+1}{k+n}{\apbs}-\yuuo{1}{n}{\apbs}\su{0}^\mpi\sub{k}-\SLuo{n-1}{\apbs}\sdiag{\su{0}^\mpi}{n-1}\yuub{k+1}{k+n}.
\end{align*}
 According to \rrem{F.R.abp=c-1}, then \(X=\sau{0}\su{0}^\mpi\sub{k}\).
 The application of \rlem{ab.L1846b} with \(n-1\) instead of \(n\) and \(k+1\) instead of \(k\) implies
\begin{multline*}
 \SLub{n-1}\sdiag{\su{0}^\mpi}{n-1}\yuuo{k+1}{k+n}{\apbs}-\SLuo{n-1}{\apbs}\sdiag{\su{0}^\mpi}{n-1}\yuub{k+1}{k+n}\\
 =\ba\sdiag{\su{0}\su{0}^\mpi}{n-1}\yuub{k+1}{k+n}-\ug\yuub{0}{n-1}\su{0}^\mpi\sub{k}.
\end{multline*}
 Hence, \(
 Y
 =\ba\sdiag{\su{0}\su{0}^\mpi}{n-1}\yuub{k+1}{k+n}-\rk{\yuuo{1}{n}{\apbs}+\ug\yuub{0}{n-1}}\su{0}^\mpi\sub{k}
 \).
 Thus, \eqref{ab.L1413b.A} holds true, since \rremss{F.R.abp=c-1}{F.R.c=ab} yield
\[
 \yuuo{1}{n}{\apbs}+\ug\yuub{0}{n-1}
 =\yuuo{0}{n-1}{c}+\ug\yuub{0}{n-1}
 =\yuub{1}{n}.\qedhere
\]
\eproof

\bleml{ab.L1555b1}
 Suppose \(\kappa\geq2\).
 Let \(\seqska \in\Dpqka\).
 Denote by \(\seqapb{\kappa-1}\) the \tbmodv{\seqsa{\kappa-1}}.
 For all \(n\in\mn{1}{\kappa-1}\), then
\beql{ab.L1555b1.A}
 \SLub{n}\sdiag{\su{0}^\mpi}{n}\OIpuu{1}{n}\yuuo{1}{n}{\apbs}-\SLuo{n}{\apbs}\sdiag{\su{0}^\mpi}{n}\yuub{0}{n}
 =\yuub{0}{n}\su{0}^\mpi\sau{0}.
\eeq
\elem
\bproof
 Let \(n\in\mn{1}{\kappa-1}\).
 Obviously, for \(k=0\) the matrix on the left-hand side of \eqref{ab.L1413b.A} coincides with the matrix on the left-hand side of \eqref{ab.L1555b1.A}.
 Thus, taking into account \rlem{ab.L1413b} and the \tbr{} \(\yuub{0}{n}=\tmat{\sub{0}\\\yuub{1}{n}}\), it is sufficient to prove that \(\sau{0}\su{0}^\mpi\sub{0}=\sub{0}\su{0}^\mpi\sau{0}\) and \(\ba\sdiag{\su{0}\su{0}^\mpi}{n}\yuub{1}{n}-\yuub{1}{n}\su{0}^\mpi\sub{0}=\yuub{1}{n}\su{0}^\mpi\sau{0}\) are valid.
 Obviously, \(\sau{0}\su{0}^\mpi\sub{0}=\sub{0}\su{0}^\mpi\sau{0}\) is fulfilled, by virtue of \rlem{ab.L0907}.
 Because of \(\seqska \in\Dpqka\) and \rrem{A.R.rs+}, we have \(\ran{\sub{j}}\subseteq\ran{\su{0}}\) and \(\nul{\su{0}}\subseteq\nul{\sub{j}}\) for all \(j\in\mn{0}{\kappa-1}\).
 Using \rrem{R.AA+B}, we infer then
\(
 \ba\sdiag{\su{0}\su{0}^\mpi}{n}\yuub{1}{n}
 =\ba\yuub{1}{n}
 =\yuub{1}{n}\su{0}^\mpi\rk{\ba\su{0}}
\).
 By virtue of \rrem{F.R.012}, hence
\(
 \ba\sdiag{\su{0}\su{0}^\mpi}{n}\yuub{1}{n}-\yuub{1}{n}\su{0}^\mpi\sub{0}
 =\yuub{1}{n}\su{0}^\mpi\rk{\ba\su{0}-\sub{0}}
 =\yuub{1}{n}\su{0}^\mpi\sau{0}
\).
\eproof

\bleml{ab.L1030b1}
 Suppose \(\kappa\geq3\).
 Let \(\seqska \in\Lggeqka\).
 Denote by \(\seqapb{\kappa-1}\) the \tbmodv{\seqsa{\kappa-1}}.
 For all \(n\in\N\) with \(2n+1\leq\kappa\), then
\beql{ab.L1030b1.A}
 \SLub{n}\sdiag{\su{0}^\mpi}{n}\OIquu{1}{n}\Guo{n-1}{\apbs}-\SLuo{n}{\apbs}\sdiag{\su{0}^\mpi}{n}\Hub{n}\OIquu{1}{n}
 =
 \bMat
  \sau{0}\su{0}^\mpi\zuub{1}{n}\\
  \ba\LLub{n}+\yuub{1}{n}\sub{0}^\mpi\dia{1}\sub{0}^\mpi\zuub{1}{n}
 \eMat.
\eeq
\elem
\bproof
 Consider an arbitrary \(n\in\N\) with \(2n+1\leq\kappa\).
 According to \rlem{ab.L1413b}, we have \eqref{ab.L1413b.A} for all \(k\in\mn{1}{n}\).
 Taking into account \(\Guo{n-1}{\apbs}=\mat{\yuuo{2}{n+1}{\apbs},\yuuo{3}{n+2}{\apbs},\dotsc,\yuuo{n+1}{2n}{\apbs}}\) as well as \(\Hub{n}\OIquu{1}{n}=\mat{\yuub{1}{n+1},\yuub{2}{n+2},\dotsc,\yuub{n}{2n}}\) and \(\Gub{n-1}=\mat{\yuub{2}{n+1},\yuub{3}{n+2},\dotsc,\yuub{n+1}{2n}}\), these \(n\) equations can be subsumed  as columns in the single equation
\[
 \SLub{n}\sdiag{\su{0}^\mpi}{n}\OIquu{1}{n}\Guo{n-1}{\apbs}-\SLuo{n}{\apbs}\sdiag{\su{0}^\mpi}{n}\Hub{n}\OIquu{1}{n}
 =
 \bMat
  \sau{0}\su{0}^\mpi\zuub{1}{n}\\
  \ba\sdiag{\su{0}\su{0}^\mpi}{n}\Gub{n-1}-\yuub{1}{n}\su{0}^\mpi\zuub{1}{n}
 \eMat.
\]
 Using \rrem{K.R.Letr}, we conclude that the sequence \(\seqs{2n+1}\) belongs to \(\Lggequ{2n+1}\) and, because of \eqref{Lgg2n-1e}, hence belongs to \(\Hggequ{2n+1}\).
 According to \rprop{H.R.He<D}, thus \(\seqs{2n+1}\in\Dqqu{2n+1}\).
 By virtue of \rrem{A.R.rs+}, then \(\ran{\sub{j}}\subseteq\ran{\su{0}}\) follows for all \(j\in\mn{0}{2n}\).
 Consequently, \rremp{R.AA+B}{R.AA+B.a} yields \(\su{0}\su{0}^\mpi\sub{0}=\sub{0}\) and \(\sdiag{\su{0}\su{0}^\mpi}{n}\Gub{n-1}=\Gub{n-1}\).
 In view of \rrem{H.R.LL}, we get therefore
\(
 \ba\sdiag{\su{0}\su{0}^\mpi}{n}\Gub{n-1}
 =\ba\rk{\LLub{n}+\yuub{1}{n}\sub{0}^\mpi\zuub{1}{n}}
\)
 and, hence,
\[
 \ba\sdiag{\su{0}\su{0}^\mpi}{n}\Gub{n-1}-\yuub{1}{n}\su{0}^\mpi\zuub{1}{n}
 =\ba\LLub{n}+\yuub{1}{n}\rk{\ba\sub{0}^\mpi-\su{0}^\mpi}\zuub{1}{n}.
\]
 In view of \eqref{Lgg2n-1e}, we have \(\seqsb{2n}\in\Hggqu{2n}\).
 Consequently, \rprop{P4-25} yields \(\seqsb{2n}\in\Dtqqu{2n}\).
 Using \rrem{R.AA+B}, then we infer \(\sub{0}\sub{0}^\mpi\zuub{1}{n}=\zuub{1}{n}\) and \(\yuub{1}{n}\sub{0}^\mpi\sub{0}=\yuub{1}{n}\).
 Taking into account \(\su{0}\su{0}^\mpi\sub{0}=\sub{0}\), we have, by virtue of \rremss{F.R.012}{ab.L0907}, furthermore
\[
 \sub{0}\rk{\ba\sub{0}^\mpi-\su{0}^\mpi}\sub{0}
 =\ba\sub{0}-\sub{0}\su{0}^\mpi\sub{0}
 =\ba\su{0}\su{0}^\mpi\sub{0}-\sub{0}\su{0}^\mpi\sub{0}
 =\rk{\ba\su{0}-\sub{0}}\su{0}^\mpi\sub{0}
 =\sau{0}\su{0}^\mpi\sub{0}
 =\dia{1}.
\]
 Consequently, we obtain \eqref{ab.L1030b1.A} from
\[
 \yuub{1}{n}\rk{\ba\sub{0}^\mpi-\su{0}^\mpi}\zuub{1}{n}
 =\yuub{1}{n}\sub{0}^\mpi\sub{0}\rk{\ba\sub{0}^\mpi-\su{0}^\mpi}\sub{0}\sub{0}^\mpi\zuub{1}{n}
 =\yuub{1}{n}\sub{0}^\mpi\dia{1}\sub{0}^\mpi\zuub{1}{n}.\qedhere
\]
\eproof

\bleml{ab.L0959b1}
 Let \(n\in\N\) and let \(\seqs{2n+1}\in\Fggqu{2n+1}\).
 Denote by \(\seqapb{2n}\) the \tbmodv{\seqsa{2n}}.
 Then
\begin{multline}\label{ab.L0959b1.A}
  \SLub{n}\OIquu{1}{n}\sdiag{\su{0}^\mpi}{n-1}\OIquu{1}{n}^\ad\Huo{n}{\apbs}-\SLuo{n}{\apbs}\sdiag{\su{0}^\mpi}{n}\Hub{n}\\
  =\yuub{0}{n}\sub{0}^\mpi\dia{1}\sub{0}^\mpi\zuub{0}{n}+\ba\OIquu{1}{n}\LLub{n}\OIquu{1}{n}^\ad.
\end{multline}
\elem
\bproof
 Let
 \[
  \SLub{n}\OIquu{1}{n}\sdiag{\su{0}^\mpi}{n-1}\OIquu{1}{n}^\ad\Huo{n}{\apbs}-\SLuo{n}{\apbs}\sdiag{\su{0}^\mpi}{n}\Hub{n}
  =\mat{A,B}
 \]
 be the \tbr{} of the matrix on the left-hand side of \eqref{ab.L0959b1.A} with \taaa{(n+1)q}{q}{block} \(A\).
 In view of \rrem{H.R.Hblock}, we have
\begin{align*}
 \OIquu{1}{n}^\ad\Huo{n}{\apbs}&=\mat{\yuuo{1}{n}{\apbs},\Guo{n-1}{\apbs}}&
&\text{and}&
 \Hub{n}&=\mat{\yuub{0}{n},\Hub{n}\OIquu{1}{n}}.
\end{align*}
 Consequently, we obtain
\begin{align*}
 A&=\SLub{n}\OIquu{1}{n}\sdiag{\su{0}^\mpi}{n-1}\yuuo{1}{n}{\apbs}-\SLuo{n}{\apbs}\sdiag{\su{0}^\mpi}{n}\yuub{0}{n}
\shortintertext{and}
 B&=\SLub{n}\OIquu{1}{n}\sdiag{\su{0}^\mpi}{n-1}\Guo{n-1}{\apbs}-\SLuo{n}{\apbs}\sdiag{\su{0}^\mpi}{n}\Hub{n}\OIquu{1}{n}.
\end{align*}
 Since \rprop{F.R.Fgg<D} shows that \(\seqs{2n+1}\in\Dqqu{2n+1}\), \rlem{ab.L1555b1} in combination with \eqref{DA3} yields
\[
 A
 =\SLub{n}\sdiag{\su{0}^\mpi}{n}\OIquu{1}{n}\yuuo{1}{n}{\apbs}-\SLuo{n}{\apbs}\sdiag{\su{0}^\mpi}{n}\yuub{0}{n}
 =\yuub{0}{n}\su{0}^\mpi\sau{0}.
\]
 Because of~\zitaa{MR3775449}{\cprop{9.2}{28}}, the sequence \(\seqs{2n+1}\) belongs to \(\Lggequ{2n+1}\).
 Thus, \rlem{ab.L1030b1} in combination with \eqref{DA3} shows that \(B\) coincides with the block matrix on the right-hand side of \eqref{ab.L1030b1.A}.
 Hence, we have
\[
 \mat{A,B}
 =
 \begin{pmat}[{|}]
  \sub{0}\su{0}^\mpi\sau{0}&\sau{0}\su{0}^\mpi\zuub{1}{n}\cr\-
  \yuub{1}{n}\su{0}^\mpi\sau{0} &\ba\LLub{n}+\yuub{1}{n}\sub{0}^\mpi\dia{1}\sub{0}^\mpi\zuub{1}{n}\cr
 \end{pmat}
 =M+\ba\OIquu{1}{n}\LLub{n}\OIquu{1}{n}^\ad
\]
 where
\[
 M\defeq
 \bMat
  \sub{0}\su{0}^\mpi\sau{0}&\sau{0}\su{0}^\mpi\zuub{1}{n}\\
  \yuub{1}{n}\su{0}^\mpi\sau{0} &\yuub{1}{n}\sub{0}^\mpi\dia{1}\sub{0}^\mpi\zuub{1}{n}
 \eMat.
\]
 From \rrem{ab.L0907} we infer \(\sub{0}\sub{0}^\mpi\dia{1}\sub{0}^\mpi\sub{0}=\sub{0}\su{0}^\mpi\sau{0}\).
 As in the proof of \rlem{ab.L1030b1}, we can conclude \(\sub{0}\sub{0}^\mpi\zuub{1}{n}=\zuub{1}{n}\) and \(\yuub{1}{n}\sub{0}^\mpi\sub{0}=\yuub{1}{n}\).
 Taking into account \rrem{ab.L0907}, we obtain hence \(\sub{0}\sub{0}^\mpi\dia{1}\sub{0}^\mpi\zuub{1}{n}=\sau{0}\su{0}^\mpi\zuub{1}{n}\) and \(\yuub{1}{n}\sub{0}^\mpi\dia{1}\sub{0}^\mpi\sub{0}=\yuub{1}{n}\su{0}^\mpi\sau{0}\).
 Consequently,
\[
 M
 =
 \bMat
  \sub{0}\sub{0}^\mpi\dia{1}\sub{0}^\mpi\sub{0}&\sub{0}\sub{0}^\mpi\dia{1}\sub{0}^\mpi\zuub{1}{n}\\
  \yuub{1}{n}\sub{0}^\mpi\dia{1}\sub{0}^\mpi\sub{0}&\yuub{1}{n}\sub{0}^\mpi\dia{1}\sub{0}^\mpi\zuub{1}{n}
 \eMat
 =\yuub{0}{n}\sub{0}^\mpi\dia{1}\sub{0}^\mpi\zuub{0}{n}
\]
 follows, implying \eqref{ab.L0959b1.A}.
\eproof

 We are now able to basically reduce the \tbHm{} \(\Huo{n}{t}\) built from the \tFT{} \(\seqt{2n}\) of a sequence \(\seqs{2n+1}\in\Fggqu{2n+1}\) to a block diagonal matrix consisting of the matrices \(\dia{1}\) and \(\LLub{n}\) given via \eqref{F.G.d01} and \rnota{H.N.LL}, \tresp{}
 In the following proof, we make use of the \tHTion{} introduced in \rdefn{103.S31}:
 
\bpropl{F.P.FTH}
 Let \(n\in\N\) and let \(\seqs{2n+1}\in\Fggqu{2n+1}\) with \tFT{} \(\seqt{2n}\).
 Then the \tbHm{} \(\Huo{n}{t}\) admits the representations
\beql{F.P.FTH.A}
 \Huo{n}{t}
 =\Rqa{n}{\obg}\sdiag{\sau{0}}{n}\SLau{n}^\mpi\rk{\yuub{0}{n}\sub{0}^\mpi\dia{1}\sub{0}^\mpi\zuub{0}{n}+\ba\OIquu{1}{n}\LLub{n}\OIquu{1}{n}^\ad}\SUau{n}^\mpi\sdiag{\sau{0}}{n}\ek*{\Rqa{n}{\obg}}^\ad
\eeq
and
\beql{F.P.FTH.B}
 \Huo{n}{t}
 =\Rqa{n}{\obg}\Dloou{a}{b}{n}\cdot\diag\rk{\dia{1},\ba\Dlub{n-1}\LLub{n}\Drub{n-1}}\cdot\Droou{a}{b}{n}\ek*{\Rqa{n}{\obg}}^\ad.
\eeq
 In particular, \(\rank\Huo{n}{t}=\rank\dia{1}+\rank\LLub{n}\) and \(\det\Huo{n}{t}=\ba^{nq}\det(\dia{1})\det(\LLub{n})\).
\eprop
\bproof
 Denote by \(\seqapb{2n}\) the \tbmodv{\seqsa{2n}}, by \(\seq{\mathbf{h}_j}{j}{0}{2n}\) the \trFa{\(\seqapb{2n}\)}, and by \(\seq{x_j}{j}{0}{2n}\) the \tCPa{\(\seqsb{2n}\)}{\(\seq{\mathbf{h}_j}{j}{0}{2n}\)}.
 In view of \rdefn{ab.N0940}, we have \(t_j=-\sau{0}\su{0}^\mpi  x_j\sau{0}\) for all \(j\in\mn{0}{2n}\).
 In particular,
\(
 \Huo{n}{t}
 =\sdiag{-\sau{0}\su{0}^\mpi}{n}\Huo{n}{x}\sdiag{\sau{0}}{n}
\).
 According to \rprop{P1522}, we have
\[
 \Huo{n}{x}
 =\Huo{n}{\subs}\SUuo{n}{\mathbf{h}}+\rk{\zdiag{\Oqq}{\SLuo{n-1}{\subs}}}\Huo{n}{\mathbf{h}}
 =\Hub{n}\SUuo{n}{\mathbf{h}}+\rk{\zdiag{\Oqq}{\SLub{n-1}}}\Huo{n}{\mathbf{h}}.
\]
 The application of \rthm{H.T1431} to the sequence \(\seqapb{2n}\) yields furthermore
\[
 \Huo{n}{\mathbf{h}}
 =\yuuo{0}{n}{\mathbf{h}}\vqu{n}^\ad + \vqu{n} \zuuo{0}{n}{\mathbf{h}}-\SLuo{n}{\mathbf{h}}\Huo{n}{\apbs}\SUuo{n}{\mathbf{h}}.
\]
 Since obviously \(\rk{\zdiag{\Oqq}{\SLub{n-1}}}\vqu{n}=\Ouu{\rk{n+1}q}{q}\) holds true, we get then
\beql{F.P.FTH.0}
 \Huo{n}{t}
 =\sdiag{-\sau{0}\su{0}^\mpi}{n}\ek*{\Hub{n}\SUuo{n}{\mathbf{h}}+\rk{\zdiag{\Oqq}{\SLub{n-1}}}\rk{\yuuo{0}{n}{\mathbf{h}}\vqu{n}^\ad-\SLuo{n}{\mathbf{h}}\Huo{n}{\apbs}\SUuo{n}{\mathbf{h}}}}\sdiag{\sau{0}}{n}.
\eeq
 Applying \rrem{ab.R1115b} to the sequence \(\seqsa{2n}\), we obtain
\begin{align}\label{F.P.FTH.Sg}
 \SLuo{n}{\apbs}&=-\ek{\Rqa{n}{\obg}}^\inv\SLau{n}&
&\text{and}&
 \SUuo{n}{\apbs}&=-\SUau{n}\ek{\Rqa{n}{\obg}}^\invad.
\end{align}
 Denote by \(\seqr{2n}\) the \trFa{\(\seqsa{2n}\)}.
 The application of \rlem{ab.R1835b} to the sequence \(\seqsa{2n}\) yields
\begin{align}\label{F.P.FTH.Sh}
 \SLuo{n}{\mathbf{h}}&=-\Rqa{n}{\obg}\SLuo{n}{r}&
&\text{and}&
 \SUuo{n}{\mathbf{h}}&=-\SUuo{n}{r}\ek*{\Rqa{n}{\obg}}^\ad.
\end{align}
 According to \rpropp{F.R.Fgg<D}{F.R.Fgg<D.b}, the sequence \(\seqsa{2n}\) belongs to \(\Dqqu{2n}\).
 Applying \rprop{101.S216} and \rlem{103.M39} to the sequence \(\seqsa{2n}\), we get then
\begin{align}
 \SLau{n}^\mpi&=\SLuo{n}{r},&\SUau{n}^\mpi&=\SUuo{n}{r}\label{F.P.FTH.Sr}
\shortintertext{and}
 \SLau{n}\SLau{n}^\mpi = \sdiag{\sau{0} \sau{0}^\mpi}{n}&=\SUau{n}\SUau{n}^\mpi,&\SLau{n}^\mpi\SLau{n}&= \sdiag{\sau{0}^\mpi \sau{0}}{n}=\SUau{n}^\mpi\SUau{n}.\label{F.P.FTH.aa+}
\end{align}
 In view of \rcor{C1511} it follows \(\mathbf{h}_j\apb{0}\apb{0}^\mpi=\mathbf{h}_j\) for all \(j\in\mn{0}{2n}\).
 Thus, \(\yuuo{0}{n}{\mathbf{h}}\apb{0}\apb{0}^\mpi=\yuuo{0}{n}{\mathbf{h}}\).
 Using \eqref{SR1} and \eqref{F.P.FTH.Sg}--\eqref{F.P.FTH.aa+}, we can conclude furthermore
\begin{align}\label{F.P.FTH.Sgh}
 \SLuo{n}{\apbs}\SLuo{n}{\mathbf{h}}&=\sdiag{\sau{0} \sau{0}^\mpi}{n}&
&\text{and}&
 \SUuo{n}{\apbs}\SUuo{n}{\mathbf{h}}&=\sdiag{\sau{0} \sau{0}^\mpi}{n}.
\end{align}
 By virtue of \rrem{F.R.abp=c-1} and \eqref{DA1}, hence
\beql{F.P.FTH.2}
 \yuuo{0}{n}{\mathbf{h}}\vqu{n}^\ad
 =\yuuo{0}{n}{\mathbf{h}}\apb{0}\apb{0}^\mpi\vqu{n}^\ad
 =\yuuo{0}{n}{\mathbf{h}}\sau{0}\sau{0}^\mpi\vqu{n}^\ad
 =\yuuo{0}{n}{\mathbf{h}}\vqu{n}^\ad\sdiag{\sau{0}\sau{0}^\mpi}{n}
 =\yuuo{0}{n}{\mathbf{h}}\vqu{n}^\ad\SUuo{n}{\apbs}\SUuo{n}{\mathbf{h}}
\eeq
 follows.
 In view of \rremp{R1326}{R1326.b}, \eqref{RA1}, \eqref{F.P.FTH.aa+}, and \eqref{F.P.FTH.Sg}, we obtain
\beql{F.P.FTH.RbSa+Sg}\begin{split}
 \sdiag{-\sau{0}\su{0}^\mpi}{n}
 &=-\Rqa{n}{\obg}\ek*{\Rqa{n}{\obg}}^\inv\sdiag{\sau{0}\sau{0}^\mpi \sau{0}\su{0}^\mpi}{n}
 =-\Rqa{n}{\obg}\sdiag{\sau{0}}{n}\sdiag{\sau{0}^\mpi \sau{0}}{n}\ek*{\Rqa{n}{\obg}}^\inv\sdiag{\su{0}^\mpi}{n}\\
 &=-\Rqa{n}{\obg}\sdiag{\sau{0}}{n}\SLau{n}^\mpi\SLau{n}\ek*{\Rqa{n}{\obg}}^\inv\sdiag{\su{0}^\mpi}{n}
 =\Rqa{n}{\obg}\sdiag{\sau{0}}{n}\SLau{n}^\mpi\SLuo{n}{\apbs}\sdiag{\su{0}^\mpi}{n}.
\end{split}\eeq 
 Since the combination of the second equations in \eqref{F.P.FTH.Sh} and \eqref{F.P.FTH.Sr} yields moreover \(\SUuo{n}{\mathbf{h}}=-\SUau{n}^\mpi\ek{\Rqa{n}{\obg}}^\ad\), we infer from \eqref{F.P.FTH.0}, \eqref{F.P.FTH.2}, and \eqref{F.P.FTH.RbSa+Sg}, using additionally \eqref{RA2}, then
\beql{F.P.FTH.1}
 \begin{split}
  \Huo{n}{t}
  &=\sdiag{-\sau{0}\su{0}^\mpi}{n}\ek*{\Hub{n}+\rk{\zdiag{\Oqq}{\SLub{n-1}}}\rk{\yuuo{0}{n}{\mathbf{h}}\vqu{n}^\ad\SUuo{n}{\apbs}-\SLuo{n}{\mathbf{h}}\Huo{n}{\apbs}}}\SUuo{n}{\mathbf{h}}\sdiag{\sau{0}}{n}\\
  &=\Rqa{n}{\obg}\sdiag{\sau{0}}{n}\SLau{n}^\mpi\\
  &\quad\times\ek*{\SLuo{n}{\apbs}\sdiag{\su{0}^\mpi}{n}\rk{\zdiag{\Oqq}{\SLub{n-1}}}\rk{\SLuo{n}{\mathbf{h}}\Huo{n}{\apbs}-\yuuo{0}{n}{\mathbf{h}}\vqu{n}^\ad\SUuo{n}{\apbs}}-\SLuo{n}{\apbs}\sdiag{\su{0}^\mpi}{n}\Hub{n}}\\
  &\quad\qquad\times\SUau{n}^\mpi\sdiag{\sau{0}}{n}\ek*{\Rqa{n}{\obg}}^\ad.
 \end{split}
\eeq 
 According to \rremss{H.R.Sblock}{H.R.Hblock}, we have \(\vqu{n}^\ad\SUuo{n}{\apbs}=\vqu{n}^\ad\Huo{n}{\apbs}\).
 Thus,
 \[
  \SLuo{n}{\mathbf{h}}\Huo{n}{\apbs}-\yuuo{0}{n}{\mathbf{h}}\vqu{n}^\ad\SUuo{n}{\apbs}
  =\rk{\SLuo{n}{\mathbf{h}}-\yuuo{0}{n}{\mathbf{h}}\vqu{n}^\ad}\Huo{n}{\apbs}.
 \]
 Because of \rrem{H.R.Sblock} and \rremp{H.R.updoT}{H.R.updoT.a}, furthermore
\[
 \zdiag{\Oqq}{\SLub{n-1}}
 =\SLub{n}\OIquu{1}{n}\OIquu{1}{n}^\ad
\]
 and
\[
 \SLuo{n}{\mathbf{h}}-\yuuo{0}{n}{\mathbf{h}}\vqu{n}^\ad
 =\zdiag{\Oqq}{\SLuo{n-1}{\mathbf{h}}}
 =\OIquu{1}{n}\SLuo{n-1}{\mathbf{h}}\OIquu{1}{n}^\ad.
\]
 hold true.
 By virtue of \rremp{H.R.updoT}{H.R.updoT.b} and \eqref{SD3}, we obtain
\[
  \rk{\zdiag{\Oqq}{\SLub{n-1}}}\rk{\SLuo{n}{\mathbf{h}}-\yuuo{0}{n}{\mathbf{h}}\vqu{n}^\ad}
  =\SLub{n}\OIquu{1}{n}\SLuo{n-1}{\mathbf{h}}\OIquu{1}{n}^\ad
  =\SLub{n}\SLuo{n}{\mathbf{h}}\OIquu{1}{n}\OIquu{1}{n}^\ad,
\]
 and, therefore,
\[
  \rk{\zdiag{\Oqq}{\SLub{n-1}}}\rk{\SLuo{n}{\mathbf{h}}\Huo{n}{\apbs}-\yuuo{0}{n}{\mathbf{h}}\vqu{n}^\ad\SUuo{n}{\apbs}}
  =\SLub{n}\SLuo{n}{\mathbf{h}}\OIquu{1}{n}\OIquu{1}{n}^\ad\Huo{n}{\apbs}.
 \]
 According to \rpropp{F.R.Fgg<D}{F.R.Fgg<D.a}, the sequence \(\seqs{2n+1}\) belongs to \(\Dqqu{2n+1}\).
 From \rlem{ab.L1746} and \rrem{A.R.AB=BA}, we thus conclude \(\SLau{n}\sdiag{\su{0}^\mpi}{n}\SLub{n}=\SLub{n}\sdiag{\su{0}^\mpi}{n}\SLau{n}\).
 In view of the first equations in \eqref{F.P.FTH.Sg}, \eqref{SR1}, and \eqref{RA1}, then \(\SLuo{n}{\apbs}\sdiag{\su{0}^\mpi}{n}\SLub{n}=\SLub{n}\sdiag{\su{0}^\mpi}{n}\SLuo{n}{\apbs}\) follows.
 In combination with the first equation in \eqref{F.P.FTH.Sgh} and \eqref{DA3}, we get then
\begin{multline*}
  \SLuo{n}{\apbs}\sdiag{\su{0}^\mpi}{n}\SLub{n}\SLuo{n}{\mathbf{h}}\OIquu{1}{n}\OIquu{1}{n}^\ad
  =\SLub{n}\sdiag{\su{0}^\mpi}{n}\SLuo{n}{\apbs}\SLuo{n}{\mathbf{h}}\OIquu{1}{n}\OIquu{1}{n}^\ad\\
  =\SLub{n}\sdiag{\su{0}^\mpi}{n}\sdiag{\sau{0} \sau{0}^\mpi}{n}\OIquu{1}{n}\OIquu{1}{n}^\ad
  =\SLub{n}\OIquu{1}{n}\sdiag{\su{0}^\mpi}{n-1}\OIquu{1}{n}^\ad\sdiag{\sau{0} \sau{0}^\mpi}{n}.
\end{multline*}
 Taking into account \rrem{F.R.c=ab} and \(\seqsa{2n}\in\Dqqu{2n}\), we get from \rrem{A.R.rs+} furthermore \(\ran{\sab{j}}\subseteq\ran{\sau{0}}\) and \(\nul{\sau{0}}\subseteq\nul{\sab{j}}\) for all \(j\in\mn{0}{2n-1}\).
 By virtue of \rremssp{F.R.abp=c-1}{R.AA+B}{R.AA+B.a}, then \(\sau{0}\sau{0}^\mpi\apb{j}=\apb{j}\) follows for all \(j\in\mn{0}{2n}\).
 In particular, \(\sdiag{\sau{0} \sau{0}^\mpi}{n}\Huo{n}{\apbs}=\Huo{n}{\apbs}\).
 Consequently,
\begin{multline}\label{F.P.FTH.13}
  \SLuo{n}{\apbs}\sdiag{\su{0}^\mpi}{n}\rk{\zdiag{\Oqq}{\SLub{n-1}}}\rk{\SLuo{n}{\mathbf{h}}\Huo{n}{\apbs}-\yuuo{0}{n}{\mathbf{h}}\vqu{n}^\ad\SUuo{n}{\apbs}}\\
  =\SLuo{n}{\apbs}\sdiag{\su{0}^\mpi}{n}\SLub{n}\SLuo{n}{\mathbf{h}}\OIquu{1}{n}\OIquu{1}{n}^\ad\Huo{n}{\apbs}
  =\SLub{n}\OIquu{1}{n}\sdiag{\su{0}^\mpi}{n-1}\OIquu{1}{n}^\ad\Huo{n}{\apbs}.
\end{multline}
 \rlem{ab.L0959b1} yields moreover \eqref{ab.L0959b1.A}.
 Substituting \eqref{F.P.FTH.13} into \eqref{F.P.FTH.1}, we can then use \eqref{ab.L0959b1.A} to conclude \eqref{F.P.FTH.A}. 
 From \rrem{H.R.Sblock} it is readily seen that \(\SLub{n}\vqu{n}=\yuub{0}{n}\) and \(\vqu{n}^\ad\SUub{n}=\zuub{0}{n}\) hold true.
 Taking additionally into account \eqref{DA1}, consequently
\beql{F.P.FTH.17}\begin{split}
 \yuub{0}{n}\sub{0}^\mpi\dia{1}\sub{0}^\mpi\zuub{0}{n}
 &=\SLub{n}\vqu{n}\sub{0}^\mpi\dia{1}\dia{1}^\mpi\dia{1}\sub{0}^\mpi\vqu{n}^\ad\SUub{n}\\
 &=\SLub{n}\sdiag{\sub{0}^\mpi\dia{1}}{n}\vqu{n}\dia{1}^\mpi\vqu{n}^\ad\sdiag{\dia{1}\sub{0}^\mpi}{n}\SUub{n}.
\end{split}\eeq
 Because of \rpropss{F.L.sabF}{ab.R0951}, the sequence \(\seqsb{2n}\) belongs to \(\Hggequ{2n}\).
 According to \rprop{H.R.He<D}, hence \(\seqsb{2n}\in\Dqqu{2n}\).
 Using \rrem{R.AA+B}, then
\begin{align}\label{F.P.FTH.5b}
 \sub{0}\sub{0}^\mpi\sub{j}&=\sub{j}&
&\text{and}&
 \sub{j}\sub{0}^\mpi\sub{0}&=\sub{j}&\text{for all }j&\in\mn{0}{2n}
\end{align}
 follow.
 Because of \rrem{H.R.LL}, in particular
\beql{F.P.FTH.10}
 \sdiag{\sub{0}\sub{0}^\mpi}{n-1}\LLub{n}\sdiag{\sub{0}^\mpi\sub{0}}{n-1}
 =\LLub{n}.
\eeq
 In view of \rnota{H.N.epXi}, we have furthermore \(\Defuuo{n-1}{2n}{\subs}=\NM\), by virtue of \eqref{F.P.FTH.5b}.
 Denote by \(\seq{h_j}{j}{0}{2n-2}\) the \tHTv{\(\seqsb{2n}\)}.
 Taking into account \(\seqsb{2n}\in\Dqqu{2n}\), we can conclude from \rprop{103.TS1A} then the representations
\begin{align}
 \Huo{n-1}{h}&=\sdiag{\sub{0}}{n-1}\SLub{n-1}^\mpi\LLub{n}\SUub{n-1}^\mpi\sdiag{\sub{0}}{n-1}\label{F.P.FTH.6a}
\shortintertext{and}
 \Huo{n-1}{h}&=\Dlub{n-1}\LLub{n}\Drub{n-1}.\label{F.P.FTH.6b}
\end{align}
 The application of \rlem{103.M39} to the sequence \(\seqsb{2n}\) yields moreover
\begin{align*}
 \SLub{n-1}\SLub{n-1}^\mpi&= \sdiag{\sub{0} \sub{0}^\mpi}{n-1}= \SUub{n-1}\SUub{n-1}^\mpi 
\shortintertext{and}
 \SLub{n-1}^\mpi\SLub{n-1}&= \sdiag{\sub{0}^\mpi \sub{0}}{n-1}=\SUub{n-1}^\mpi\SUub{n-1}.
\end{align*}
 In combination with \eqref{F.P.FTH.10}, \eqref{F.P.FTH.6a}, and \rremp{R1326}{R1326.b}, we infer
\beql{F.P.FTH.8}\begin{split}
 \LLub{n}
 &=\sdiag{\sub{0}\sub{0}^\mpi}{n-1}\LLub{n}\sdiag{\sub{0}^\mpi\sub{0}}{n-1}
 =\SLub{n-1}\SLub{n-1}^\mpi\LLub{n}\SUub{n-1}^\mpi\SUub{n-1}\\
 &=\SLub{n-1}\SLub{n-1}^\mpi\SLub{n-1}\SLub{n-1}^\mpi\LLub{n}\SUub{n-1}^\mpi\SUub{n-1}\SUub{n-1}^\mpi\SUub{n-1}\\
 &=\SLub{n-1}\sdiag{\sub{0}^\mpi \sub{0}}{n-1}\SLub{n-1}^\mpi\LLub{n}\SUub{n-1}^\mpi\sdiag{\sub{0} \sub{0}^\mpi}{n-1}\SUub{n-1}\\
 &=\SLub{n-1}\sdiag{\sub{0}^\mpi}{n-1}\Huo{n-1}{h}\sdiag{\sub{0}^\mpi}{n-1}\SUub{n-1}.
\end{split}\eeq 
 Denote by \(\hpseqo{2n}\) the \thpfa{\(\seqsb{2n}\)}.
 Because of \(\seqsb{2n}\in\Hggequ{2n}\), we get from \rthm{H.P.T915} in particular \(\hpu{2}=h_0\).
 Using \rremss{H.R.h2L}{F.R.ABL}, we infer then \(h_0=\Lub{1}=\osc{3}\).
 According to \rrem{F.R.ABL}, we have furthermore
\begin{align}\label{F.P.FTH.18}
 \usc{1}&=\Lau{0}=\sau{0}&
&\text{and}&
 \osc{1}&=\Lub{0}=\sub{0}.
\end{align}
 In view of \(\seqsb{2n}\in\Hggequ{2n}\), by virtue of \rprop{H.C.HTinHgge} the sequence \(\seq{h_j}{j}{0}{2n-2}\) belongs to \(\Hggequ{2n-2}\).
 Hence, \rprop{H.R.He<D} implies \(\seq{h_j}{j}{0}{2n-2}\in\Dqqu{2n-2}\).
 According to \rprop{ab.C1101}, we have \(\ran{\osc{3}}\subseteq\ran{\osc{2}}\) and \(\nul{\osc{2}}\subseteq\nul{\osc{3}}\), whereas \rprop{ab.R13371422} yields \(\ran{\osc{2}}\subseteq\ran{\dia{1}}\) and \(\nul{\dia{1}}\subseteq\nul{\osc{2}}\).
 Consequently,
\(
 \ran{h_j}
 \subseteq\ran{h_0}
 =\ran{\osc{3}}
 \subseteq\ran{\osc{2}}
 \subseteq\ran{\dia{1}}
\)
 and, analogously, \(\nul{\dia{1}}\subseteq\nul{h_j}\) follow for all \(j\in\mn{0}{2n-2}\).
 In view of \rrem{R.AA+B}, hence
\beql{F.P.FTH.9}
 \sdiag{\dia{1}\dia{1}^\mpi}{n-1}\Huo{n-1}{h}\sdiag{\dia{1}^\mpi\dia{1}}{n-1}
 =\Huo{n-1}{h}.
\eeq
 Using \rremp{R1326}{R1326.b}, \eqref{SD3}, and \eqref{DA3}, we infer from \eqref{F.P.FTH.8} then
\beql{F.P.FTH.16}\begin{split}
 \OIquu{1}{n}\LLub{n}\OIquu{1}{n}^\ad
 &=\OIquu{1}{n}\SLub{n-1}\sdiag{\sub{0}^\mpi\dia{1}\dia{1}^\mpi}{n-1}\Huo{n-1}{h}\sdiag{\dia{1}^\mpi\dia{1}\sub{0}^\mpi}{n-1}\SUub{n-1}\OIquu{1}{n}^\ad\\
 &=\SLub{n}\sdiag{\sub{0}^\mpi\dia{1}}{n}\OIquu{1}{n}\sdiag{\dia{1}^\mpi}{n-1}\Huo{n-1}{h}\sdiag{\dia{1}^\mpi}{n-1}\OIquu{1}{n}^\ad\sdiag{\dia{1}\sub{0}^\mpi}{n}\SUub{n}.
\end{split}\eeq 
 According to \rprop{ab.R13371422} and \eqref{F.P.FTH.18}, we have \(\ran{\dia{1}}=\ran{\sau{0}}\cap\ran{\sub{0}}\) and \(\nul{\sau{0}}\cup\nul{\sub{0}}=\nul{\dia{1}}\).
 By virtue of \eqref{F.P.FTH.Sr} and \rrem{ab.R1520}, thus
\begin{align*}
 \Dloou{\saus}{\subs}{n}\sdiag{\dia{1}}{n}&=\sdiag{\sau{0}}{n}\SLau{n}^\mpi\SLub{n}\sdiag{\sub{0}^\mpi\dia{1}}{n}
\shortintertext{and}
 \sdiag{\dia{1}}{n}\Droou{\saus}{\subs}{n}&=\sdiag{\dia{1}\sub{0}^\mpi}{n}\SUub{n}\SUau{n}^\mpi\sdiag{\sau{0}}{n}.
\end{align*}
 In combination with \eqref{F.P.FTH.17} and \eqref{F.P.FTH.16}, we get then
\[\begin{split}
 &\sdiag{\sau{0}}{n}\SLau{n}^\mpi\rk{\yuub{0}{n}\sub{0}^\mpi\dia{1}\sub{0}^\mpi\zuub{0}{n}+\ba\OIquu{1}{n}\LLub{n}\OIquu{1}{n}^\ad}\SUau{n}^\mpi\sdiag{\sau{0}}{n}\\
 &=\sdiag{\sau{0}}{n}\SLau{n}^\mpi\SLub{n}\sdiag{\sub{0}^\mpi\dia{1}}{n}\\
 &\qquad\times\rk*{\vqu{n}\dia{1}^\mpi\vqu{n}^\ad+\ba\OIquu{1}{n}\sdiag{\dia{1}^\mpi}{n-1}\Huo{n-1}{h}\sdiag{\dia{1}^\mpi}{n-1}\OIquu{1}{n}^\ad}\sdiag{\dia{1}\sub{0}^\mpi}{n}\SUub{n}\SUau{n}^\mpi\sdiag{\sau{0}}{n}\\
 &=\Dloou{\saus}{\subs}{n}\sdiag{\dia{1}}{n}\rk*{\vqu{n}\dia{1}^\mpi\vqu{n}^\ad+\ba\OIquu{1}{n}\sdiag{\dia{1}^\mpi}{n-1}\Huo{n-1}{h}\sdiag{\dia{1}^\mpi}{n-1}\OIquu{1}{n}^\ad}\sdiag{\dia{1}}{n}\Droou{\saus}{\subs}{n}.
\end{split}\]
 Furthermore, by virtue of \eqref{DA1}, we have
\(
 \sdiag{\dia{1}}{n}\vqu{n}\dia{1}^\mpi\vqu{n}^\ad\sdiag{\dia{1}}{n}
 =\vqu{n}\dia{1}\dia{1}^\mpi\dia{1}\vqu{n}^\ad
 =\vqu{n}\dia{1}\vqu{n}^\ad
\)
 and, because of \eqref{DA3}, \rremp{R1326}{R1326.b}, \eqref{F.P.FTH.9}, \eqref{F.P.FTH.6b}, moreover
\begin{multline*}
 \sdiag{\dia{1}}{n}\OIquu{1}{n}\sdiag{\dia{1}^\mpi}{n-1}\Huo{n-1}{h}\sdiag{\dia{1}^\mpi}{n-1}\OIquu{1}{n}^\ad\sdiag{\dia{1}}{n}
 =\OIquu{1}{n}\sdiag{\dia{1}\dia{1}^\mpi}{n-1}\Huo{n-1}{h}\sdiag{\dia{1}^\mpi\dia{1}}{n-1}\OIquu{1}{n}^\ad\\
 =\OIquu{1}{n}\Huo{n-1}{h}\OIquu{1}{n}^\ad
 =\OIquu{1}{n}\Dlub{n-1}\LLub{n}\Drub{n-1}\OIquu{1}{n}^\ad.
\end{multline*}
 Consequently,
\begin{multline*}
 \sdiag{\sau{0}}{n}\SLau{n}^\mpi\rk{\yuub{0}{n}\sub{0}^\mpi\dia{1}\sub{0}^\mpi\zuub{0}{n}+\ba\OIquu{1}{n}\LLub{n}\OIquu{1}{n}^\ad}\SUau{n}^\mpi\sdiag{\sau{0}}{n}\\
 =\Dloou{\saus}{\subs}{n}\rk{\vqu{n}\dia{1}\vqu{n}^\ad+\ba\OIquu{1}{n}\Dlub{n-1}\LLub{n}\Drub{n-1}\OIquu{1}{n}^\ad}\Droou{\saus}{\subs}{n}.
\end{multline*}
 Hence, we obtain \eqref{F.P.FTH.B} from \eqref{F.P.FTH.A}. 
 Taking into account \(\ba>0\), the formulas for \(\rank\Huo{n}{t}\) and \(\det\Huo{n}{t}\) follow, in view of \rremsss{H.R.RinLU}{F.R.detDd}{H.R.DinLU}, from \eqref{F.P.FTH.B}.
\eproof

 In the following, we use the equivalence relation ``\(\ldu\)'' introduced in \rnota{A.N.ldusim}:
 
\bcorl{F.C.FTL}
 Let \(n\in\N\) and let \(\seqs{2n+1}\in\Fggqu{2n+1}\) with \tFT{} \(\seqt{2n}\).
 Then
\beql{F.C.FTL.A}
 \LLuo{n}{t}
 =\ba\Rqa{n-1}{\obg}\Dloou{a}{b}{n-1}\Dlub{n-1}\LLub{n}\Drub{n-1}\Droou{a}{b}{n-1}\ek*{\Rqa{n-1}{\obg}}^\ad
\eeq 
 and, in particular, \(\rank\LLuo{n}{t}=\rank\LLub{n}\), \(\det\LLuo{n}{t}=\ba^{nq}\det\LLub{n}\), and \(\LLuo{n}{t}\ldu\ba\LLub{n}\).
\ecor
\bproof
 From \rnota{F.N.Dd} and \rrem{H.R.Sblock} we see that
\begin{align*}
 \Dloou{a}{b}{n}
 &=
 \bMat
  \Iq&\Ouu{q}{nq}\\
  \ub&\Dloou{a}{b}{n-1}
 \eMat&
&\text{and}&
 \Droou{a}{b}{n}
 &=
 \bMat
  \Iq&\ub\\
  \Ouu{nq}{q}&\Droou{a}{b}{n-1}
 \eMat.
\end{align*}
 In view of \rexam{H.E.ResS} and \rrem{H.R.Sblock}, we have
\begin{align*}
 \Rqa{n}{\obg}
 &=
 \bMat
  \Iq & \Ouu{q}{nq}\\
  \ub& \Rqa{n-1}{\obg}
 \eMat&
&\text{and}&
 \ek*{\Rqa{n}{\obg}}^\ad
 &=
 \bMat
  \Iq&\ub \\
  \Ouu{nq}{q}&\ek{\Rqa{n-1}{\obg}}^\ad 
 \eMat.
\end{align*}
 Let \(M\defeq\Rqa{n-1}{\obg}\Dloou{a}{b}{n-1}\), let \(U\defeq\Iq\), let \(S\defeq\Droou{a}{b}{n-1}\ek{\Rqa{n-1}{\obg}}^\ad\), and let  \(V\defeq\Iq\).
 Then there are matrices \(L\) and \(R\) such that
\begin{align*}
 \Rqa{n}{\obg}\Dloou{a}{b}{n}
 &=
 \bMat
  U&\Ouu{q}{nq}\\
  L&M
 \eMat&
&\text{and}&
 \Droou{a}{b}{n}\ek*{\Rqa{n}{\obg}}^\ad
 &=
 \bMat
  V&R\\
  \Ouu{nq}{q}&S
 \eMat.
\end{align*}
 Let \(E\defeq\diag\rk{\dia{1},\ba\Dlub{n-1}\LLub{n}\Drub{n-1}}\) and let \(F\defeq\Huo{n}{t}\).
 \rprop{F.P.FTH} yields
\[
 F
 =
 \bMat
  U&\Ouu{q}{nq}\\
  L&M
 \eMat E
 \bMat
  V&R\\
  \Ouu{nq}{q}&S
 \eMat.
 \]
 Denote by \(E=\tmat{A&B\\ C&D}\) the \tbr{} of \(E\) with \tqqa{Block} \(A\).
 Since the matrices \(U\) and  \(V\) are unitary, the application of \rlem{ab.R1038} thus yields \(F=\tmat{A&\ub\\\ub&\ub}\) and, because of \(B=\NM\) and \(C=\NM\), furthermore \(F\schca A=M(E\schca A)S=MDS\).
 From \rnota{H.N.LL} we infer \(\LLuo{n}{t}=F\schca A\).
 Consequently, \eqref{F.C.FTL.A} follows. 
 Taking into account \(\ba>0\), the formulas for rank and determinant of \(\LLuo{n}{t}\) can be obtained, in view of \rremsss{H.R.RinLU}{F.R.detDd}{H.R.DinLU}, from \eqref{F.C.FTL.A}.
 Analogously, \(\LLuo{n}{t}\ldu\ba\LLub{n}\) follows from \eqref{F.C.FTL.A}, using \rremssss{H.R.RinLU}{F.R.DdinLU}{H.R.DinLU}{R1547} and \rnota{A.N.ldusim}.
\eproof

 Now we turn our attention to the \tbHms{} built via \eqref{F.G.HaHbHK} from the \tFT{} \(\seqt{2n+1}\) of a sequence \(\seqs{2n+2}\) belonging to \(\Fggqu{2n+2}\).
 We see that the matrix \(-\ug\Huo{n}{t}+\Kuo{n}{t}\) can be basically traced back to \(\Hab{n}\):

\bpropl{F.P.FTHa}
 Let \(n\in\NO\) and let \(\seqs{2n+2}\in\Fggqu{2n+2}\) with \tFT{} \(\seqt{2n+1}\).
 Let the sequence \(\seq{u_j}{j}{0}{2n}\) be given by \(u_j\defeq-\ug t_j+t_{j+1}\).
 Then
\begin{align}
 \Huo{n}{u}
 &=\ba\Rqa{n}{\obg}\sdiag{\sau{0}}{n}\SLau{n}^\mpi\Hab{n}\SUau{n}^\mpi\sdiag{\sau{0}}{n}\ek*{\Rqa{n}{\obg}}^\ad\label{F.P.FTHa.A}
\shortintertext{and}
 \Huo{n}{u}
 &=\ba\Rqa{n}{\obg}\Dlau{n}\Hab{n}\Drau{n}\ek*{\Rqa{n}{\obg}}^\ad\label{F.P.FTHa.B}.
\end{align}
 In particular, \(\rank\Huo{n}{u}=\rank\Hab{n}\), \(\det\Huo{n}{u}=\ba^{(n+1)q}\det\Hab{n}\), and \(\Huo{n}{u}\ldu\ba\Hab{n}\).
\eprop
\bproof
 Denote by \(\seqtpa{2n+1}\) the \tamodv{\seqt{2n+1}}.
 In view of \rdefn{114.D1455}, then
 \(
  \Kuo{n}{\tpas}
  =\Huo{n}{u}
 \).
 Denote by \(\seqbpa{2n+1}\) the \tamodv{\seqsb{2n+1}}, by \(\seqapb{2n+1}\) the \tbmodv{\seqsa{2n+1}}, by \(\seq{\mathbf{h}_j}{j}{0}{2n+1}\) the \trFa{\(\seqapb{2n+1}\)}, and by \(\seq{\mathbf{x}_j}{j}{0}{2n+1}\) the \tCPa{\(\seqbpa{2n+1}\)}{\(\seq{\mathbf{h}_j}{j}{0}{2n+1}\)}.
 According to \rlemp{ab.R1110}{ab.R1110.a}, then \(\tpa{j}=-\sau{0}\su{0}^\mpi \mathbf{x}_j\sau{0}\) for all \(j\in\mn{0}{2n+1}\).
 In particular,
 \(
  \Kuo{n}{\tpas}
  =\sdiag{-\sau{0}\su{0}^\mpi}{n}\Kuo{n}{\mathbf{x}}\sdiag{\sau{0}}{n}
 \).
 Because of \rprop{P1522}, we have
 \(
  \Kuo{n}{\mathbf{x}}
  =\Kuo{n}{\bpas}\SUuo{n}{\mathbf{h}}+\SLuo{n}{\bpas}\Kuo{n}{\mathbf{h}}
 \).
 The application of \rthm{H.P0836} to the sequence \(\seqapb{2n+1}\) yields furthermore
 \(
  \Kuo{n}{\mathbf{h}}
  =-\SLuo{n}{\mathbf{h}}\Kuo{n}{\apbs}\SUuo{n}{\mathbf{h}}
 \).
 From \rrem{F.R.abp=c-1} one can easily see \(\Kuo{n}{\bpas}=\Hab{n}=\Kuo{n}{\apbs}\).
 Consequently,
\beql{F.P.FTHa.0}
 \Huo{n}{u}
 =\Kuo{n}{\tpas}
 =\sdiag{-\sau{0}\su{0}^\mpi}{n}\Kuo{n}{\mathbf{x}}\sdiag{\sau{0}}{n}
 =\sdiag{-\sau{0}\su{0}^\mpi}{n}\rk{\Hab{n}\SUuo{n}{\mathbf{h}}-\SLuo{n}{\bpas}\SLuo{n}{\mathbf{h}}\Hab{n}\SUuo{n}{\mathbf{h}}}\sdiag{\sau{0}}{n}.
\eeq
 Applying \rlem{ab.R1835b} to the sequence \(\seqsa{2n+1}\), we obtain \eqref{F.P.FTH.Sh}, where \(\seqr{2n+1}\) is the \trFa{\(\seqsa{2n+1}\)}.
 The application of \rrem{ab.R1115a} to \(\seqsb{2n+1}\) and of \rrem{ab.R1115b} to \(\seqsa{2n+1}\) provides us
 \begin{equation}\label{F.P.FTHa.Sf}
  \SLuo{n}{\bpas}
  =\ek{\Rqa{n}{\ug}}^\inv\SLub{n}
 \end{equation}
 and \eqref{F.P.FTH.Sg}.
 According to \rprop{F.R.Fgg<D}, we have \(\seqsa{2n+1}\in\Dqqu{2n+1}\).
 The application of \rprop{101.S216} and \rlem{103.M39} to the sequence \(\seqsa{2n+1}\) yields then \eqref{F.P.FTH.Sr} and \eqref{F.P.FTH.aa+}, \tresp{}
 By virtue of \rremp{R1326}{R1326.b}, \eqref{RA1}, \eqref{F.P.FTH.aa+}, and \eqref{F.P.FTH.Sg}, we get \eqref{F.P.FTH.RbSa+Sg}.
 Combining the second equations in \eqref{F.P.FTH.Sh} and \eqref{F.P.FTH.Sr}, we obtain furthermore \(\SUuo{n}{\mathbf{h}}=-\SUau{n}^\mpi\ek{\Rqa{n}{\obg}}^\ad\).
 Consequently, using \eqref{F.P.FTH.RbSa+Sg} and \eqref{RA2}, we infer from \eqref{F.P.FTHa.0} that
 \begin{equation}\label{F.P.FTHa.1}
  \begin{split}
  \Huo{n}{u}
  &=\Rqa{n}{\obg}\sdiag{\sau{0}}{n}\SLau{n}^\mpi\rk{\SLuo{n}{\apbs}\sdiag{\su{0}^\mpi}{n}-\SLuo{n}{\apbs}\sdiag{\su{0}^\mpi}{n}\SLuo{n}{\bpas}\SLuo{n}{\mathbf{h}}}\Hab{n}\SUuo{n}{\mathbf{h}}\sdiag{\sau{0}}{n}\\
  &=\Rqa{n}{\obg}\sdiag{\sau{0}}{n}\SLau{n}^\mpi\rk{\SLuo{n}{\apbs}\sdiag{\su{0}^\mpi}{n}\SLuo{n}{\bpas}\SLuo{n}{\mathbf{h}}-\SLuo{n}{\apbs}\sdiag{\su{0}^\mpi}{n}}\Hab{n}\SUau{n}^\mpi\sdiag{\sau{0}}{n}\ek*{\Rqa{n}{\obg}}^\ad.
 \end{split}
 \end{equation}
 According to \rprop{F.R.Fgg<D}, we have \(\seqs{2n+2}\in\Dqqu{2n+2}\).
 Using \rlem{ab.L1746} and \rrem{A.R.AB=BA}, we get then \(\SLau{n}\sdiag{\su{0}^\mpi}{n}\SLub{n}=\SLub{n}\sdiag{\su{0}^\mpi}{n}\SLau{n}\).
 By virtue of \eqref{F.P.FTHa.Sf}, the first equations in \eqref{F.P.FTH.Sg}, \eqref{SR1}, and \eqref{RA1}, and formula~\eqref{SR4}, thus \(\SLuo{n}{\apbs}\sdiag{\su{0}^\mpi}{n}\SLuo{n}{\bpas}=\SLuo{n}{\bpas}\sdiag{\su{0}^\mpi}{n}\SLuo{n}{\apbs}\) holds true.
 From \eqref{SR1} and \eqref{F.P.FTH.Sg}--\eqref{F.P.FTH.aa+}, furthermore \eqref{F.P.FTH.Sgh} follows.
 Because of \rrem{A.R.rs+}, we have, in view of \rrem{F.R.c=ab} and \(\seqsa{2n+1}\in\Dqqu{2n+1}\), obviously \(\ran{\sab{j}}\subseteq\ran{\sau{0}}\) and \(\nul{\sau{0}}\subseteq\nul{\sab{j}}\) for all \(j\in\mn{0}{2n}\).
 \rrem{R.AA+B} yields then
\begin{align}\label{F.P.FTHa.5}
 \sau{0}\sau{0}^\mpi\sab{j}&=\sab{j}&
&\text{and}&
 \sab{j}\sau{0}^\mpi\sau{0}&=\sab{j}&\text{for all }j&\in\mn{0}{2n}.
\end{align}
 In particular, \(\sdiag{\sau{0} \sau{0}^\mpi}{n}\Hab{n}=\Hab{n}\).
 Thus, from the first equation in \eqref{F.P.FTH.Sgh}, we get
 \beql{F.P.FTHa.13}
  \SLuo{n}{\apbs}\sdiag{\su{0}^\mpi}{n}\SLuo{n}{\bpas}\SLuo{n}{\mathbf{h}}\Hab{n}
  =\SLuo{n}{\bpas}\sdiag{\su{0}^\mpi}{n}\SLuo{n}{\apbs}\SLuo{n}{\mathbf{h}}\Hab{n}
  =\SLuo{n}{\bpas}\sdiag{\su{0}^\mpi}{n}\Hab{n}.
 \eeq
 From \rrem{ab.L1118}, we know that \(\SLuo{n}{\bpas}-\SLuo{n}{\apbs}=\ba\sdiag{\su{0}}{n}\).
 Because of \rrem{F.R.Fgg-r}, we have \(\ran{\sau{0}}\subseteq\ran{\su{0}}\) and \(\nul{\su{0}}\subseteq\nul{\sau{0}}\).
 According to \rremp{R.AA+B}{R.AA+B.a}, then \(\su{0}\su{0}^\mpi\sau{0}=\sau{0}\).
 In view of \eqref{F.P.FTHa.5}, therefore \(\su{0}\su{0}^\mpi\sab{j}=\sab{j}\) holds true for all \(j\in\mn{0}{2n}\), implying \(\sdiag{\su{0} \su{0}^\mpi}{n}\Hab{n}=\Hab{n}\).
 By virtue of \eqref{F.P.FTHa.13} and \rremp{R1326}{R1326.b}, we thus obtain
 \[
  \rk{\SLuo{n}{\apbs}\sdiag{\su{0}^\mpi}{n}\SLuo{n}{\bpas}\SLuo{n}{\mathbf{h}}-\SLuo{n}{\apbs}\sdiag{\su{0}^\mpi}{n}}\Hab{n}
  =\rk{\SLuo{n}{\bpas}-\SLuo{n}{\apbs}}\sdiag{\su{0}^\mpi}{n}\Hab{n}
  =\ba\Hab{n}.
 \]
 Substituting this into \eqref{F.P.FTHa.1}, we get \eqref{F.P.FTHa.A}. 
 According to \eqref{F.P.FTHa.5}, we have \(\sdiag{\sau{0}\sau{0}^\mpi}{n}\Hab{n}\sdiag{\sau{0}^\mpi\sau{0}}{n}=\Hab{n}\).
 In view of \eqref{F.P.FTH.Sr}, the application of \rrem{ab.R1602} to the sequence \(\seqsa{2n+1}\) provides us
\begin{align}\label{F.P.FTHa.Daa}
 \Dlau{n}\sdiag{\sau{0}}{n}&=\sdiag{\sau{0}}{n}\SLau{n}^\mpi\sdiag{\sau{0}}{n}&
&\text{and}&
 \sdiag{\sau{0}}{n}\Drau{n}&=\sdiag{\sau{0}}{n}\SUau{n}^\mpi\sdiag{\sau{0}}{n}. 
\end{align}
 Taking additionally into account \rremp{R1326}{R1326.b}, consequently
\[\begin{split}
 \sdiag{\sau{0}}{n}\SLau{n}^\mpi\Hab{n}\SUau{n}^\mpi\sdiag{\sau{0}}{n}
 &=\sdiag{\sau{0}}{n}\SLau{n}^\mpi\sdiag{\sau{0}\sau{0}^\mpi}{n}\Hab{n}\sdiag{\sau{0}^\mpi\sau{0}}{n}\SUau{n}^\mpi\sdiag{\sau{0}}{n}\\
 &=\Dlau{n}\sdiag{\sau{0}\sau{0}^\mpi}{n}\Hab{n}\sdiag{\sau{0}^\mpi\sau{0}}{n}\Drau{n}
 =\Dlau{n}\Hab{n}\Drau{n}.
\end{split}\]
 Hence, \eqref{F.P.FTHa.B} follows from \eqref{F.P.FTHa.A}. 
 The formulas for rank and determinant of \(\Huo{n}{u}\) can be seen, in view of \rremss{H.R.RinLU}{H.R.DinLU} and \(\ba>0\), from \eqref{F.P.FTHa.B}.
 Analogously, \(\Huo{n}{u}\ldu\ba\Hab{n}\) follows from \eqref{F.P.FTHa.B}, using \rremsss{H.R.RinLU}{H.R.DinLU}{R1547} and \rnota{A.N.ldusim}.
\eproof

 To obtain in \rprop{F.P.FTHb} a convenient form of the \tbHm{} \(\obg\Huo{n}{t}-\Kuo{n}{t}\), we need some further technical results:

\bleml{ab.L1539a}
 Suppose \(\kappa\geq1\).
 Let \(\seqska\) be a sequence of complex \tpqa{matrices} and let \(n\in\mn{0}{\kappa-1}\).
 Let the sequence \(\seq{h_j}{j}{0}{\kappa-1}\) be given by \(h_j\defeq\su{j+1}\).
 For all \(k\in\mn{0}{\kappa-n-1}\), then
\[
 \SLu{n}\sdiag{\su{0}^\mpi}{n}\yuu{k+1}{k+n+1}-\SLuo{n}{h}\sdiag{\su{0}^\mpi}{n}\yuu{k}{k+n}
 =\sdiag{\su{0}\su{0}^\mpi}{n}\yuu{k+1}{k+n+1}-\yuu{1}{n+1}\su{0}^\mpi\su{k}.
\]
\elem
\bproof
 Consider an arbitrary \(k\in\mn{0}{\kappa-n-1}\).
 In view of \rnotass{H.N.T}{H.N.updo}, we have
\(
 \Tp{n}\yuu{k+1}{k+n+1}-\yuu{k}{k+n}
 =-\vpu{n}\su{k}
\).
 Taking into account \rrem{H.R.Sblock}, consequently
\[
 \SLuo{n}{h}\sdiag{\su{0}^\mpi}{n}\rk{\Tp{n}\yuu{k+1}{k+n+1}-\yuu{k}{k+n}}\\
 =-\SLuo{n}{h}\vqu{n}\su{0}^\mpi\su{k}
 =-\yuu{1}{n+1}\su{0}^\mpi\su{k}.
\]
 Because of \rnota{M.N.S}, furthermore \(\SLuo{n}{h}\Tq{n}+\sdiag{\su{0}}{n}=\SLu{n}\).
 Thus, using \rremp{R1326}{R1326.b} and \rremp{R1317}{R1317.b}, we conclude
\[\begin{split}
 &\SLu{n}\sdiag{\su{0}^\mpi}{n}\yuu{k+1}{k+n+1}-\SLuo{n}{h}\sdiag{\su{0}^\mpi}{n}\yuu{k}{k+n}\\
 &=\SLuo{n}{h}\Tq{n}\sdiag{\su{0}^\mpi}{n}\yuu{k+1}{k+n+1}+\sdiag{\su{0}\su{0}^\mpi}{n}\yuu{k+1}{k+n+1}-\SLuo{n}{h}\sdiag{\su{0}^\mpi}{n}\yuu{k}{k+n}\\
 &=\sdiag{\su{0}\su{0}^\mpi}{n}\yuu{k+1}{k+n+1}+\SLuo{n}{h}\sdiag{\su{0}^\mpi}{n}\Tp{n}\yuu{k+1}{k+n+1}-\SLuo{n}{h}\sdiag{\su{0}^\mpi}{n}\yuu{k}{k+n}\\
 &=\sdiag{\su{0}\su{0}^\mpi}{n}\yuu{k+1}{k+n+1}-\yuu{1}{n+1}\su{0}^\mpi\su{k}.\qedhere
\end{split}\]
\eproof

\bleml{ab.L1534a}
 Suppose \(\kappa\geq1\).
 Let \(\seqska\) be a sequence of complex \tpqa{matrices} and let \(n\in\mn{0}{\kappa-1}\).
 For all \(k\in\mn{0}{\kappa-n-1}\), then
\beql{ab.L1534a.1}
 \SLub{n}\sdiag{\su{0}^\mpi}{n}\yauu{k}{k+n}-\SLau{n}\sdiag{\su{0}^\mpi}{n}\yuub{k}{k+n}
 =\ba\rk*{\sdiag{\su{0}\su{0}^\mpi}{n}\yuu{k+1}{k+n+1}-\yuu{1}{n+1}\su{0}^\mpi\su{k}}.
\eeq
\elem
\bproof
 Let \(k\in\mn{0}{\kappa-n-1}\).
 Let the sequence \(\seq{h_j}{j}{0}{\kappa-1}\) be given by \(h_j\defeq\su{j+1}\).
 Obviously, we have \(\SLub{n}=\obg\SLu{n}-\SLuo{n}{h}\) and \(\SLau{n}=-\ug\SLu{n}+\SLuo{n}{h}\) as well as \(\yauu{k}{k+n}=-\ug\yuu{k}{k+n}+\yuu{k+1}{k+n+1}\) and \(\yuub{k}{k+n}=\obg\yuu{k}{k+n}-\yuu{k+1}{k+n+1}\).
 Consequently, a straightforward calculation yields
\[
 \SLub{n}\sdiag{\su{0}^\mpi}{n}\yauu{k}{k+n}-\SLau{n}\sdiag{\su{0}^\mpi}{n}\yuub{k}{k+n}
 =\ba\SLu{n}\sdiag{\su{0}^\mpi}{n}\yuu{k+1}{k+n+1}-\ba\SLuo{n}{h}\sdiag{\su{0}^\mpi}{n}\yuu{k}{k+n}.
\]
 Using \rlem{ab.L1539a}, then \eqref{ab.L1534a.1} follows.
\eproof

 We are now able to basically reduce the \tbHm{} \(\obg\Huo{n}{t}-\Kuo{n}{t}\) to the Schur complement \(\LLu{n+1}\) of \(\su{0}\) in \(\Hu{n+1}\):

\bpropl{F.P.FTHb}
 Let \(n\in\NO\) and let \(\seqs{2n+2}\in\Fggqu{2n+2}\) with \tFT{} \(\seqt{2n+1}\).
 Let the sequence \(\seq{v_j}{j}{0}{2n}\) be given by \(v_j\defeq\obg t_j-t_{j+1}\).
 Then
\begin{align}\label{F.P.FTHb.A}
 \Huo{n}{v}&=\ba\sdiag{\sau{0}}{n}\SLau{n}^\mpi\LLu{n+1}\SUau{n}^\mpi\sdiag{\sau{0}}{n}&
&\text{and}&
 \Huo{n}{v}&=\ba\Dloou{a}{s}{n}\Dlu{n}\LLu{n+1}\Dru{n}\Droou{a}{s}{n}.
\end{align}
 In particular, \(\rank\Huo{n}{v}=\rank\LLu{n+1}\), \(\det\Huo{n}{v}=\ba^{(n+1)q}\det\LLu{n+1}\), and \(\Huo{n}{v}\ldu\ba\LLu{n+1}\).
\eprop
\bproof
 Denote by \(\seqtpb{2n+1}\) the \tbmodv{\seqt{2n+1}}.
 In view of \rdefn{ab.N1137b}, then
\(
 \Kuo{n}{\mathbf{v}}
 =\Huo{n}{v}
\).
 Denote by \(\seqr{2n+1}\) the \trFa{\(\seqsa{2n+1}\)} and by \(\seq{\mathbf{y}_j}{j}{0}{2n+1}\) the \tCPa{\(\seqsb{2n+1}\)}{\(\seq{r_j}{j}{0}{2n+1}\)}.
 According to \rlemp{ab.R1110}{ab.R1110.b}, we have \(\tpb{j}=-\sau{0}\su{0}^\mpi \mathbf{y}_j\sau{0}\) for all \(j\in\mn{0}{2n+1}\).
 In particular,
\(
 \Kuo{n}{\mathbf{v}}
 =\sdiag{-\sau{0}\su{0}^\mpi}{n}\Kuo{n}{\mathbf{y}}\sdiag{\sau{0}}{n}
\).
 From \rprop{P1522} we can conclude
\[
 \Kuo{n}{\mathbf{y}}
 =\Kuo{n}{b}\SUuo{n}{r}+\SLuo{n}{b}\Kuo{n}{r}
 =\Kub{n}\SUuo{n}{r}+\SLub{n}\Kuo{n}{r}.
\]
 The application of \rthm{H.P0836} to the sequence \(\seqsa{2n+1}\) yields furthermore
\[
 \Kuo{n}{r}
 =-\SLuo{n}{r}\Kuo{n}{a}\SUuo{n}{r}
 =-\SLuo{n}{r}\Kau{n}\SUuo{n}{r}.
\]
 Consequently,
\beql{F.P.FTHb.0}
 \Huo{n}{v}
 =\sdiag{-\sau{0}\su{0}^\mpi}{n}\Kuo{n}{\mathbf{y}}\sdiag{\sau{0}}{n}
 =\sdiag{-\sau{0}\su{0}^\mpi}{n}\rk{\Kub{n}\SUuo{n}{r}-\SLub{n}\SLuo{n}{r}\Kau{n}\SUuo{n}{r}}\sdiag{\sau{0}}{n}
\eeq
 follows.
 According to \rprop{F.R.Fgg<D}, we have \(\seqsa{2n+1}\in\Dqqu{2n+1}\).
 Applying \rprop{101.S216} and \rlem{103.M39} to the sequence \(\seqsa{2n+1}\), we obtain then \eqref{F.P.FTH.Sr} and \eqref{F.P.FTH.aa+}, \tresp{}
 From \rremp{R1326}{R1326.b} and \eqref{F.P.FTH.aa+} we infer
\beql{F.P.FTHb.Sa+Sa}
 \sdiag{-\sau{0}\su{0}^\mpi}{n}
 =-\sdiag{\sau{0}\sau{0}^\mpi \sau{0}\su{0}^\mpi}{n}
 =-\sdiag{\sau{0}}{n}\sdiag{\sau{0}^\mpi \sau{0}}{n}\sdiag{\su{0}^\mpi}{n}
 =-\sdiag{\sau{0}}{n}\SLau{n}^\mpi\SLau{n}\sdiag{\su{0}^\mpi}{n}.
\eeq
 Taking additionally into account \eqref{F.P.FTHb.0} and \eqref{F.P.FTH.Sr}, we get
\beql{F.P.FTHb.1}\begin{split}
 \Huo{n}{v}
 &=-\sdiag{\sau{0}}{n}\SLau{n}^\mpi\SLau{n}\sdiag{\su{0}^\mpi}{n}\rk{\Kub{n}\SUau{n}^\mpi-\SLub{n}\SLau{n}^\mpi\Kau{n}\SUau{n}^\mpi}\sdiag{\sau{0}}{n}\\
 &=\sdiag{\sau{0}}{n}\SLau{n}^\mpi\rk*{\SLau{n}\sdiag{\su{0}^\mpi}{n}\SLub{n}\SLau{n}^\mpi\Kau{n}-\SLau{n}\sdiag{\su{0}^\mpi}{n}\Kub{n}}\SUau{n}^\mpi\sdiag{\sau{0}}{n}.
\end{split}\eeq 
 Because of \rprop{ab.R0951}, the sequence \(\seqs{2n+2}\) belongs to \(\Hggequ{2n+2}\) and, according to \rprop{H.R.He<D}, hence it belongs to \(\Dqqu{2n+2}\).
 By virtue of \rlem{ab.L1746} and \rrem{A.R.AB=BA}, we can conclude then \(\SLau{n}\sdiag{\su{0}^\mpi}{n}\SLub{n}=\SLub{n}\sdiag{\su{0}^\mpi}{n}\SLau{n}\).
 In view of \(\seqsa{2n+1}\in\Dqqu{2n+1}\), furthermore \rremp{R.AA+B}{R.AA+B.a} yields \(\sau{0}\sau{0}^\mpi\sau{j}=\sau{j}\) for all \(j\in\mn{0}{2n+1}\) and, consequently, \(\sdiag{\sau{0} \sau{0}^\mpi}{n}\Kau{n}=\Kau{n}\).
 Using \eqref{F.P.FTH.aa+}, we infer then
\beql{F.P.FTHb.4}
 \SLau{n}\sdiag{\su{0}^\mpi}{n}\SLub{n}\SLau{n}^\mpi\Kau{n}
 =\SLub{n}\sdiag{\su{0}^\mpi}{n}\sdiag{\sau{0} \sau{0}^\mpi}{n}\Kau{n}
 =\SLub{n}\sdiag{\su{0}^\mpi}{n}\Kau{n}.
\eeq
 Because of \rlem{ab.L1534a}, we have \eqref{ab.L1534a.1} for all \(k\in\mn{0}{n+1}\).
 In view of \(\Kau{n}=\mat{\yauu{1}{n+1},\yauu{2}{n+2},\dotsc,\yauu{n+1}{2n+1}}\), \(\Kub{n}=\mat{\yuub{1}{n+1},\yuub{2}{n+2},\dotsc,\yuub{n+1}{2n+1}}\), and \(\Gu{n}=\mat{\yuu{2}{n+2},\yuu{3}{n+3},\dotsc,\yuu{n+2}{2n+2}}\), the equations \eqref{ab.L1534a.1} for \(k\in\mn{1}{n+1}\) can be summarized to the single equation
\beql{F.P.FTHb.7}
 \SLub{n}\sdiag{\su{0}^\mpi}{n}\Kau{n}-\SLau{n}\sdiag{\su{0}^\mpi}{n}\Kub{n}
 =\ba\rk*{\sdiag{\su{0}\su{0}^\mpi}{n}\Gu{n}-\yuu{1}{n+1}\su{0}^\mpi\zuu{1}{n+1}}.
\eeq
 By virtue of \(\seqs{2n+2}\in\Dqqu{2n+2}\), we obtain from \rrem{R.AA+B} furthermore
\begin{align}\label{F.P.FTHb.5}
 \su{0}\su{0}^\mpi\su{j}&=\su{j}&
&\text{and}&
 \su{j}\su{0}^\mpi\su{0}&=\su{j}&\text{for all }j&\in\mn{0}{2n+2}.
\end{align}
 In particular, \(\sdiag{\su{0}\su{0}^\mpi}{n}\Gu{n}=\Gu{n}\).
 Because of \rrem{H.R.LL}, we get from \eqref{F.P.FTHb.7} then
\[
 \SLub{n}\sdiag{\su{0}^\mpi}{n}\Kau{n}-\SLau{n}\sdiag{\su{0}^\mpi}{n}\Kub{n}
 =\ba\LLu{n+1}.
\]
 Combining the last identity with \eqref{F.P.FTHb.4} and  \eqref{F.P.FTHb.1}, we infer the first equation in \eqref{F.P.FTHb.A}. 
 According to \rrem{H.R.LL} and \eqref{F.P.FTHb.5}, we have
\beql{F.P.FTHb.10}
 \sdiag{\su{0}\su{0}^\mpi}{n}\LLu{n+1}\sdiag{\su{0}^\mpi\su{0}}{n}
 =\LLu{n+1}.
\eeq
 Denote by \(\seq{h_j}{j}{0}{2n}\) the \tHTv{\(\seqs{2n+2}\)}.
 In view of \rnota{H.N.epXi}, we get \(\Defuu{n}{2n+2}=\NM\) from \eqref{F.P.FTHb.5}.
 Taking into account \(\seqs{2n+2}\in\Dqqu{2n+2}\), the application of \rprop{103.TS1A} thus yields the representations
\begin{align}\label{F.P.FTHb.6}
 \Huo{n}{h}&=\sdiag{\su{0}}{n}\SLu{n}^\mpi\LLu{n+1}\SUu{n}^\mpi\sdiag{\su{0}}{n}&
&\text{and}&
 \Huo{n}{h}&=\Dlu{n}\LLu{n+1}\Dru{n}.
\end{align}
 By virtue of \rlem{103.M39}, we have furthermore \(\SLu{n}\SLu{n}^\mpi=\sdiag{\su{0} \su{0}^\mpi}{n}= \SUu{n}\SUu{n}^\mpi\) and \(\SLu{n}^\mpi\SLu{n}= \sdiag{\su{0}^\mpi \su{0}}{n}=\SUu{n}^\mpi\SUu{n}\).
 In combination with \eqref{F.P.FTHb.10}, \rremp{R1326}{R1326.b}, and the first equation in \eqref{F.P.FTHb.6}, then
\beql{F.P.FTHb.8}\begin{split}
 \LLu{n+1}
 =\SLu{n}\SLu{n}^\mpi\LLu{n+1}\SUu{n}^\mpi\SUu{n}
 &=\SLu{n}\SLu{n}^\mpi\SLu{n}\SLu{n}^\mpi\LLu{n+1}\SUu{n}^\mpi\SUu{n}\SUu{n}^\mpi\SUu{n}\\
 &=\SLu{n}\sdiag{\su{0}^\mpi \su{0}}{n}\SLu{n}^\mpi\LLu{n+1}\SUu{n}^\mpi\sdiag{\su{0} \su{0}^\mpi}{n}\SUu{n}
 =\SLu{n}\sdiag{\su{0}^\mpi}{n}\Huo{n}{h}\sdiag{\su{0}^\mpi}{n}\SUu{n}
\end{split}\eeq 
 follows.
 Denote by \(\hpseqo{2n+2}\) the \thpfa{\(\seqs{2n+2}\)}.
 Because of \(\seqs{2n+2}\in\Hggequ{2n+2}\), we get from \rthm{H.P.T915} in particular \(\hpu{2}=h_0\).
 Using \rremss{H.R.h2L}{F.R.ABL}, we infer then \(h_0=\Lu{1}=\usc{2}\).
 According to \rrem{F.R.ABL}, we have \(\usc{1}=\Lau{0}=\sau{0}\).
 In view of \(\seqs{2n+2}\in\Hggequ{2n+2}\), the sequence \(\seq{h_j}{j}{0}{2n}\) belongs to \(\Hggequ{2n}\), by virtue of \rprop{H.C.HTinHgge}.
 Hence, \rprop{H.R.He<D} implies \(\seq{h_j}{j}{0}{2n}\in\Dqqu{2n}\).
 \rprop{ab.C1101} shows \(\ran{\usc{2}}\subseteq\ran{\usc{1}}\) and \(\nul{\usc{1}}\subseteq\nul{\usc{2}}\).
 Consequently,
\(
 \ran{h_j}
 \subseteq\ran{h_0}
 =\ran{\usc{2}}
 \subseteq\ran{\usc{1}}
 =\ran{\sau{0}}
\)
 and, analogously, \(\nul{\sau{0}}\subseteq\nul{h_j}\) follow for all \(j\in\mn{0}{2n}\).
 In view of \rrem{R.AA+B}, hence
\beql{F.P.FTHb.9}
 \sdiag{\sau{0}\sau{0}^\mpi}{n}\Huo{n}{h}\sdiag{\sau{0}^\mpi\sau{0}}{n}
 =\Huo{n}{h}.
\eeq
 According to \rrem{F.R.Fgg-r}, furthermore \(\ran{\sau{0}}\subseteq\ran{\su{0}}\) and \(\nul{\su{0}}\subseteq\nul{\sau{0}}\) hold true.
 By virtue of \eqref{F.P.FTH.Sr} and \rrem{ab.R1520}, thus
\[
 \Dloou{a}{s}{n}\sdiag{\sau{0}}{n}=\sdiag{\sau{0}}{n}\SLau{n}^\mpi\SLu{n}\sdiag{\su{0}^\mpi\sau{0}}{n}
\text{ and }
 \sdiag{\sau{0}}{n}\Droou{a}{s}{n}=\sdiag{\sau{0}\su{0}^\mpi}{n}\SUu{n}\SUau{n}^\mpi\sdiag{\sau{0}}{n}
\]
 are valid.
 In combination with \eqref{F.P.FTHb.8}, \eqref{F.P.FTHb.9}, \rremp{R1326}{R1326.b}, and the second equation in \eqref{F.P.FTHb.6}, then
\[\begin{split}
 \sdiag{\sau{0}}{n}\SLau{n}^\mpi\LLu{n+1}\SUau{n}^\mpi\sdiag{\sau{0}}{n}
 &=\sdiag{\sau{0}}{n}\SLau{n}^\mpi\SLu{n}\sdiag{\su{0}^\mpi}{n}\Huo{n}{h}\sdiag{\su{0}^\mpi}{n}\SUu{n}\SUau{n}^\mpi\sdiag{\sau{0}}{n}\\
 &=\sdiag{\sau{0}}{n}\SLau{n}^\mpi\SLu{n}\sdiag{\su{0}^\mpi\sau{0}\sau{0}^\mpi}{n}\Huo{n}{h}\sdiag{\sau{0}^\mpi\sau{0}\su{0}^\mpi}{n}\SUu{n}\SUau{n}^\mpi\sdiag{\sau{0}}{n}\\
 &=\Dloou{a}{s}{n}\sdiag{\sau{0}\sau{0}^\mpi}{n}\Huo{n}{h}\sdiag{\sau{0}^\mpi\sau{0}}{n}\Droou{a}{s}{n}\\
 &=\Dloou{a}{s}{n}\Huo{n}{h}\Droou{a}{s}{n}
 =\Dloou{a}{s}{n}\Dlu{n}\LLu{n+1}\Dru{n}\Droou{a}{s}{n}
\end{split}\]
 follows.
 Hence, we infer the second equation in \eqref{F.P.FTHb.A} from the first one. 
 Taking into account \(\ba>0\), the formulas for rank and determinant of \(\Huo{n}{v}\) follow, in view of \rremss{F.R.detDd}{H.R.DinLU}, from the second equation in \eqref{F.P.FTHb.A}.
 Analogously, \(\Huo{n}{v}\ldu\ba\LLu{n+1}\) follows from \eqref{F.P.FTHb.A}, using \rremsss{F.R.DdinLU}{H.R.DinLU}{R1547}, and \rnota{A.N.ldusim}.
\eproof

 To obtain in \rprop{F.P.FTHab} a convenient representation of the \tbHm{} built via \eqref{F.G.HabHKG} from the \tFT{} of a sequence, we still need one further technical result:

\bleml{ab.L1811a}
 Suppose \(\kappa\geq3\).
 Let \(\seqska\in\Kggeqka\).
 Denote by \(\seqbpa{\kappa-1}\) the \tamodv{\seqsb{\kappa-1}}.
 For all \(n\in\NO\) with \(2n+3\leq\kappa\), then
\beql{ab.L1811a.A}
 \SLuo{n}{\bpas}\sdiag{\su{0}^\mpi}{n}\LLau{n+1}-\SLau{n}\sdiag{\su{0}^\mpi}{n}\rk{\Guo{n}{\bpas}-\yuuo{1}{n+1}{\bpas}\sau{0}^\mpi\zauu{1}{n+1}}
 =\ba\LLau{n+1}.
\eeq
\elem
\bproof
 Consider an arbitrary \(n\in\NO\) with \(2n+3\leq\kappa\).
 Because of \rlem{ab.L1846b}, we have \eqref{ab.L1846b.B} for all \(k\in\mn{1}{n+2}\).
 Consequently,
\begin{multline*}
 \rk*{\SLau{n}\sdiag{\su{0}^\mpi}{n}\yuuo{1}{n+1}{\bpas}-\SLuo{n}{\bpas}\sdiag{\su{0}^\mpi}{n}\yauu{1}{n+1}}\sau{0}^\mpi\zauu{1}{n+1}\\
 =\rk*{\obg\yauu{0}{n}\su{0}^\mpi\sau{0}-\ba\sdiag{\su{0}\su{0}^\mpi}{n}\yauu{1}{n+1}}\sau{0}^\mpi\zauu{1}{n+1}
\end{multline*}
 and, by virtue of the \tbr{s} \(\Guo{n}{\bpas}=\mat{\yuuo{2}{n+2}{\bpas},\yuuo{3}{n+3}{\bpas},\dotsc,\yuuo{n+2}{2n+2}{\bpas}}\) and \(\Gau{n}=\mat{\yauu{2}{n+2},\yauu{3}{n+3},\dotsc,\yauu{n+2}{2n+2}}\), furthermore
\[
 \SLau{n}\sdiag{\su{0}^\mpi}{n}\Guo{n}{\bpas}-\SLuo{n}{\bpas}\sdiag{\su{0}^\mpi}{n}\Gau{n}
 =\obg\yauu{0}{n}\su{0}^\mpi\zauu{1}{n+1}-\ba\sdiag{\su{0}\su{0}^\mpi}{n}\Gau{n}.
\]
 Hence, we obtain
\begin{multline}\label{ab.L1811a.1}
  \SLau{n}\sdiag{\su{0}^\mpi}{n}\rk{\Guo{n}{\bpas}-\yuuo{1}{n+1}{\bpas}\sau{0}^\mpi\zauu{1}{n+1}}\\
  =\rk*{\SLuo{n}{\bpas}\sdiag{\su{0}^\mpi}{n}-\ba\sdiag{\su{0}\su{0}^\mpi}{n}}\rk{\Gau{n}-\yauu{1}{n+1}\sau{0}^\mpi\zauu{1}{n+1}}\\
  +\obg\yauu{0}{n}\su{0}^\mpi\rk{\Iq-\sau{0}\sau{0}^\mpi}\zauu{1}{n+1}.
\end{multline}
 From \rrem{K.R.Ketr} we infer that the sequence \(\seqs{2n+3}\) belongs to \(\Kggequ{2n+3}\) and, because of \eqref{Kgg2n-1e}, therefore it belongs to \(\Hggequ{2n+3}\).
 According to \rprop{H.R.He<D}, consequently \(\seqs{2n+3}\in\Dqqu{2n+3}\).
 \rrem{A.R.rs+} yields then \(\ran{\sau{j}}\subseteq\ran{\su{0}}\) for all \(j\in\mn{0}{2n+2}\).
 Using \rremp{R.AA+B}{R.AA+B.a}, we hence get \(\sdiag{\su{0}\su{0}^\mpi}{n}\yauu{1}{n+1}=\yauu{1}{n+1}\) and \(\sdiag{\su{0}\su{0}^\mpi}{n}\Gau{n}=\Gau{n}\).
 From \rrem{H.R.LL} thus
\beql{ab.L1811a.2}
 \rk*{\SLuo{n}{\bpas}\sdiag{\su{0}^\mpi}{n}-\ba\sdiag{\su{0}\su{0}^\mpi}{n}}\rk{\Gau{n}-\yauu{1}{n+1}\sau{0}^\mpi\zauu{1}{n+1}}
 =\SLuo{n}{\bpas}\sdiag{\su{0}^\mpi}{n}\LLau{n+1}-\ba\LLau{n+1}
\eeq
 follows.
 Because of \(\seqs{2n+3}\in\Kggequ{2n+3}\) and \eqref{Kgg2n-1e}, the sequence \(\seqsa{2n+2}\) belongs to \(\Hggqu{2n+2}\) and, according to \rprop{P4-25}, hence it belongs to \(\Dtqqu{2n+2}\).
 From \rremp{R.AA+B}{R.AA+B.a} we conclude \(\sau{0}\sau{0}^\mpi\zauu{1}{n+1}=\zauu{1}{n+1}\).
 Substituting this together with \eqref{ab.L1811a.2} into \eqref{ab.L1811a.1}, we get \(\SLau{n}\sdiag{\su{0}^\mpi}{n}\rk{\Guo{n}{\bpas}-\yuuo{1}{n+1}{\bpas}\sau{0}^\mpi\zauu{1}{n+1}}=\SLuo{n}{\bpas}\sdiag{\su{0}^\mpi}{n}\LLau{n+1}-\ba\LLau{n+1}\), implying \eqref{ab.L1811a.A}.
\eproof

 The \tbHm{} \(-\ug\obg\Huo{n}{t}+(\ug+\obg)\Kuo{n}{t}-\Guo{n}{t}\) can be basically reduced to the Schur complement \(\LLau{n+1}\):

\bpropl{F.P.FTHab}
 Let \(n\in\NO\) and let \(\seqs{2n+3}\in\Fggqu{2n+3}\) with \tFT{} \(\seqt{2n+2}\).
 Let the sequence \(\seq{w_j}{j}{0}{2n}\) be given by \(w_j\defeq-\ug\obg t_j+(\ug+\obg)t_{j+1}-t_{j+2}\).
 Then
\begin{align}\label{F.P.FTHab.A}
 \Huo{n}{w}&=\ba\sdiag{\sau{0}}{n}\SLau{n}^\mpi\LLau{n+1}\SUau{n}^\mpi\sdiag{\sau{0}}{n}&
&\text{and}&
 \Huo{n}{w}&=\ba\Dlau{n}\LLau{n+1}\Drau{n}.
\end{align}
 In particular, \(\rank\Huo{n}{w}=\rank\LLau{n+1}\), \(\det\Huo{n}{w}=\ba^{(n+1)q}\det\LLau{n+1}\), and \(\Huo{n}{w}\ldu\ba\LLau{n+1}\).
\eprop
\bproof
 Let \(\seqtpc{2n+2}\) be the \tabmodv{\seqt{2n+2}}.
 In view of \rdefn{ab.N1137c}, then
\(
 \Guo{n}{\mathbf{w}}
 =\Huo{n}{w}
\)
 Denote by \(\seqbpa{2n+2}\) the \tamodv{\seqsb{2n+2}}, by \(\seqr{2n+2}\) the \trFa{\(\seqsa{2n+2}\)}, and by \(\seq{\mathbf{z}_j}{j}{0}{2n+2}\) the \tCPa{\(\seqbpa{2n+2}\)}{\(\seqr{2n+2}\)}.
 According to \rlemp{ab.R1110}{ab.R1110.c}, then \(\tpc{j}=-\sau{0}\su{0}^\mpi \mathbf{z}_j\sau{0}\) for all \(j\in\mn{0}{2n+2}\).
 In particular,
\(
 \Guo{n}{\mathbf{w}}
 =\sdiag{-\sau{0}\su{0}^\mpi}{n}\Guo{n}{\mathbf{z}}\sdiag{\sau{0}}{n}
\).
 Because of \rprop{P1522}, we have
\[
 \Guo{n}{\mathbf{z}}
 =\Guo{n}{\bpas}\SUuo{n}{r}+\yuuo{1}{n+1}{\bpas}\zuuo{1}{n+1}{r}+\SLuo{n}{\bpas}\Guo{n}{r}.
\]
 The application of \rcor{ab.R1112} and \rthm{H.P0938} to the sequence \(\seqsa{2n+2}\) yields \(\zuuo{1}{n+1}{r}=-\sau{0}^\mpi\zauu{1}{n+1}\SUuo{n}{r}\) and \(\Guo{n}{r}=-\SLuo{n}{r}\LLau{n+1}\SUuo{n}{r}\),
 \tresp{}
 Consequently,
\beql{F.P.FTHab.0}
 \Huo{n}{w}
 =\sdiag{-\sau{0}\su{0}^\mpi}{n}\rk{\Guo{n}{\bpas}\SUuo{n}{r}-\yuuo{1}{n+1}{\bpas}\sau{0}^\mpi\zauu{1}{n+1}\SUuo{n}{r}-\SLuo{n}{\bpas}\SLuo{n}{r}\LLau{n+1}\SUuo{n}{r}}\sdiag{\sau{0}}{n}.
\eeq
 According to \rprop{F.R.Fgg<D}, we have \(\seqsa{2n+2}\in\Dqqu{2n+2}\).
 The application of \rprop{101.S216} and \rlem{103.M39} to the sequence \(\seqsa{2n+2}\) yields then \eqref{F.P.FTH.Sr} and \eqref{F.P.FTH.aa+}, \tresp{}
 From \rremp{R1326}{R1326.b} and \eqref{F.P.FTH.aa+} we infer \eqref{F.P.FTHb.Sa+Sa}.
 Taking additionally into account \eqref{F.P.FTHab.0} and \eqref{F.P.FTH.Sr}, we get
\beql{F.P.FTHab.1}\begin{split}
 \Huo{n}{w}
 &=-\sdiag{\sau{0}}{n}\SLau{n}^\mpi\SLau{n}\sdiag{\su{0}^\mpi}{n}\rk{\Guo{n}{\bpas}-\yuuo{1}{n+1}{\bpas}\sau{0}^\mpi\zauu{1}{n+1}-\SLuo{n}{\bpas}\SLuo{n}{r}\LLau{n+1}}\SUuo{n}{r}\sdiag{\sau{0}}{n}\\
 &=\sdiag{\sau{0}}{n}\SLau{n}^\mpi\\
 &\quad\times\ek*{\SLau{n}\sdiag{\su{0}^\mpi}{n}\SLuo{n}{\bpas}\SLau{n}^\mpi\LLau{n+1}-\SLau{n}\sdiag{\su{0}^\mpi}{n}\rk{\Guo{n}{\bpas}-\yuuo{1}{n+1}{\bpas}\sau{0}^\mpi\zauu{1}{n+1}}}\SUau{n}^\mpi\sdiag{\sau{0}}{n}.
\end{split}\eeq 
 Because of \rprop{F.R.Fgg<D}, the sequence \(\seqs{2n+3}\) belongs to \(\Dqqu{2n+3}\).
 Using \rlem{ab.L1746} and \rrem{A.R.AB=BA}, we can conclude then \(\SLau{n}\sdiag{\su{0}^\mpi}{n}\SLub{n}=\SLub{n}\sdiag{\su{0}^\mpi}{n}\SLau{n}\).
 The application of \rrem{ab.R1115a} to the sequence \(\seqsb{2n+2}\) yields furthermore \(\SLuo{n}{\bpas}=\ek{\Rqa{n}{\ug}}^\inv\SLub{n}\).
 Using additionally the first equations in \eqref{SR1} and \eqref{RA1}, consequently \(\SLau{n}\sdiag{\su{0}^\mpi}{n}\SLuo{n}{\bpas}=\SLuo{n}{\bpas}\sdiag{\su{0}^\mpi}{n}\SLau{n}\) follows.
 In view of \(\seqsa{2n+2}\in\Dqqu{2n+2}\), we infer from \rrem{R.AA+B} moreover
\begin{align}\label{F.P.FTHab.5}
 \sau{0}\sau{0}^\mpi\sau{j}&=\sau{j}&
&\text{and}&
 \sau{j}\sau{0}^\mpi\sau{0}&=\sau{j}&\text{for all }j&\in\mn{0}{2n+2}.
\end{align}
 By virtue of \rrem{H.R.LL}, thus \(\sdiag{\sau{0} \sau{0}^\mpi}{n}\LLau{n+1}=\LLau{n+1}\).
 Using \eqref{F.P.FTH.aa+}, then we obtain
\beql{F.P.FTHab.4}
 \SLau{n}\sdiag{\su{0}^\mpi}{n}\SLuo{n}{\bpas}\SLau{n}^\mpi\LLau{n+1}
 =\SLuo{n}{\bpas}\sdiag{\su{0}^\mpi}{n}\SLau{n}\SLau{n}^\mpi\LLau{n+1}
 =\SLuo{n}{\bpas}\sdiag{\su{0}^\mpi}{n}\LLau{n+1}.
\eeq
 According to~\zitaa{MR3775449}{\cprop{9.3}{28}}, we have \(\seqs{2n+3}\in\Kggequ{2n+3}\).
 Hence, \rlem{ab.L1811a} yields \eqref{ab.L1811a.A}.
 Substituting \eqref{F.P.FTHab.4} into \eqref{F.P.FTHab.1}, we can then use \eqref{ab.L1811a.A} to conclude the first equation in \eqref{F.P.FTHab.A}. 
 From \rrem{H.R.LL} and \eqref{F.P.FTHab.5} we get \(\sdiag{\sau{0}\sau{0}^\mpi}{n}\LLau{n+1}\sdiag{\sau{0}^\mpi\sau{0}}{n}=\LLau{n+1}\).
 By virtue of \eqref{F.P.FTH.Sr}, the application of \rrem{ab.R1602} to the sequence \(\seqsa{2n+2}\) yields \eqref{F.P.FTHa.Daa}.
 Taking additionally into account \rremp{R1326}{R1326.b}, consequently
\[\begin{split}
 \sdiag{\sau{0}}{n}\SLau{n}^\mpi\LLau{n+1}\SUau{n}^\mpi\sdiag{\sau{0}}{n}
 &=\sdiag{\sau{0}}{n}\SLau{n}^\mpi\sdiag{\sau{0}\sau{0}^\mpi}{n}\LLau{n+1}\sdiag{\sau{0}^\mpi\sau{0}}{n}\SUau{n}^\mpi\sdiag{\sau{0}}{n}\\
 &=\Dlau{n}\sdiag{\sau{0}\sau{0}^\mpi}{n}\LLau{n+1}\sdiag{\sau{0}^\mpi\sau{0}}{n}\Drau{n}
 =\Dlau{n}\LLau{n+1}\Drau{n}
\end{split}\]
 follows.
 Hence, we infer the second equation in \eqref{F.P.FTHab.A} from the first one. 
 Taking into account \(\ba>0\), the formulas for rank and determinant of \(\Huo{n}{w}\) follow, in view of \rrem{H.R.DinLU}, from the second equation in \eqref{F.P.FTHab.A}.
 Analogously, \(\Huo{n}{w}\ldu\ba\LLau{n+1}\) follows from \eqref{F.P.FTHab.A}, using \rrem{H.R.DinLU} and \rnota{A.N.ldusim}.
\eproof
 
 In view of \rrem{F.R.012}, representations for the \tbHms{} \(\Huo{n}{t}\), \(\Kuo{n}{t}\), and \(\Guo{n}{t}\) built via \rnota{N.HKG} from the \tFT{} \(\seqt{\kappa-1}\) of a sequence \(\seqska\in\Fggqka\) can be obtained by appropriate combinations of \rpropsss{F.P.FTHa}{F.P.FTHb}{F.P.FTHab}.

\section{A Schur type algorithm in the class \(\Fggqka\)}\label{S.F.kFT}

 We are now going to iterate the \tFT{ation} introduced in \rdefn{ab.N0940}:
 
\bdefnl{ab.N1020}
 Let \(\seqska \) be a sequence of complex \tpqa{matrices}.
 Let the sequence \(\seq{\su{j}^\FTa{0}}{j}{0}{\kappa}\) be given by
 \(
  \su{j}^\FTa{0}
  \defeq\su{j}
 \).
 If \(\kappa\geq1\), then, for all \(k\in\mn{1}{\kappa}\), let the sequence \(\seq{\su{j}^\FTa{k}}{j}{0}{\kappa-k}\) be recursively defined to be the \tFT{} of the sequence \(\seq{\su{j}^\FTa{k-1}}{j}{0}{\kappa-(k-1)}\).
 For all \(k\in\mn{0}{\kappa}\), then we call the sequence \(\seq{\su{j}^\FTa{k}}{j}{0}{\kappa-k}\) the \emph{\hnFTv{k}{\(\seqska \)}}.
\edefn

\breml{F.R.FT1}
 Suppose \(\kappa\geq1\).
 Let \(\seqska\) be a sequence of complex \tpqa{matrices}.
 Then \(\seq{\su{j}^\FT}{j}{0}{\kappa-1}\) is exactly the \tFTv{\(\seqska \)} from \rdefn{ab.N0940}.
\erem

 Since we consider the \tFT{ation} only for a fixed interval \(\ab\), we omit to indicate in our notation the dependence of the sequence \(\seq{\su{j}^\FTa{k}}{j}{0}{\kappa-k}\) on the real numbers \(\ug\) and \(\obg\).
 
\breml{F.R.kFTtr}
 Let \(k\in\mn{0}{\kappa}\) and let \(\seqska \) be a sequence of complex \tpqa{matrices} with \tnFT{k} \(\seq{u_j}{j}{0}{\kappa-k}\).
 In view of \rrem{F.R.FTtr}, we see that, for each \(\ell\in\mn{0}{\kappa-k}\), the matrix \(u_{\ell}\) is built from the matrices \(\su{0},\su{1},\dotsc,\su{\ell+k}\).
 In particular, for each \(m\in\mn{k}{\kappa}\), the \tnFTv{k}{\(\seqs{m}\)} coincides with \(\seq{u_j}{j}{0}{m-k}\).
\erem

 Now we are going to prove a result, which plays a key role in our following considerations.

\bthml{ab.P1030}
 If \(\seqska\in\Fggqka\), then \(\seq{\su{j}^\FTa{k}}{j}{0}{\kappa-k}\in\Fggqu{\kappa-k}\) for all \(k\in\mn{0}{\kappa}\).
\ethm
\bproof
 Let \(\seqska\in\Fggqka\).
 By virtue of \rdefn{ab.N1020}, then \(\seq{\su{j}^\FTa{0}}{j}{0}{\kappa}\in\Fggqu{\kappa}\).
 
 Now assume \(\kappa\geq1\) and let \(k\in\mn{1}{\kappa}\).
 Denote by \(\seqt{\kappa-1}\) the \tFTv{\(\seqska\)}.
 According to \rprop{ab.C0929}, we have \(\dia{1}\in\Cggq\).
 
 First consider the case \(k=1\).
 \rprop{F.R.Fgg<D} yields \(\seqska\in\Dqqu{\kappa}\).
 From \rlem{ab.P1728a} we get then \(\Huo{0}{t}=\dia{1}\).
 Hence, \(\Huo{0}{t}\in\Cggq\), \tie{}, \(\seqt{0}\in\Fggqu{0}\).
  
 Now assume \(\kappa\geq2\) and consider an arbitrary \(n\in\NO\) with \(2n+2\leq\kappa\).
 Because of \rprop{ab.R0933}, then \(\seqs{2n+2}\in\Fggqu{2n+2}\).
 In view of \eqref{Fgg2n}, thus \(\seqsab{2n}\in\Hggqu{2n}\) and \(\seqs{2n+2}\in\Hggqu{2n+2}\).
 Hence, the matrices \(\Hab{n}\) and \(\Hu{n+1}\) are both \tnnH{}.
 In particular, the Schur complement \(\LLu{n+1}\) of \(\su{0}\) in \(\Hu{n+1}\) is \tnnH{}, according to \rrem{L.AEP}.
 Because of \rlem{F.R.Fgg-s}, the matrices \(\su{0},\su{1},\dotsc,\su{2n+2}\) are \tH{}.
 Consequently, the matrices \(\sau{0},\sau{1},\dotsc,\sau{2n+1}\) are \tH{} as well.
 Using \rremss{H.R.D*}{F.R.Dd*}, then \(\Dlu{n}^\ad=\Dru{n}\), \(\Dlau{n}^\ad=\Drau{n}\), and \(\rk{\Dloou{\saus}{s}{n}}^\ad=\Droou{\saus}{s}{n}\) follow.
 Taking additionally into account \(\ba>0\), we can conclude, in view of \rremss{A.R.XAX}{A.R.kK}, from \rpropss{F.P.FTHa}{F.P.FTHb} that the matrices \(-\ug\Huo{n}{t}+\Kuo{n}{t}\) and \(\obg\Huo{n}{t}-\Kuo{n}{t}\) are both \tnnH{}, \tie{}, \(\seqt{2n+1}\) belongs to \(\Fggqu{2n+1}\).
  
 Now assume \(\kappa\geq3\) and consider an arbitrary \(n\in\NO\) with \(2n+3\leq\kappa\).
 By virtue of  \rprop{ab.R0933}, then \(\seqs{2n+3}\in\Fggqu{2n+3}\).
 In view of \eqref{Fgg2n+1}, thus \(\set{\seqsa{2n+2},\seqsb{2n+2}}\subseteq\Hggqu{2n+2}\).
 Hence, the matrices \(\Hau{n+1}\) and \(\Hub{n+1}\) are both \tnnH{}.
 In particular, the Schur complements \(\LLau{n+1}\) and \(\LLub{n+1}\) are both \tnnH{}, according to \rrem{L.AEP}.
 Because of \rlem{F.R.Fgg-s}, the matrices \(\su{0},\su{1},\dotsc,\su{2n+3}\) are \tH{}.
 Consequently, the matrices \(\sau{0},\sau{1},\dotsc,\sau{2n+2}\) and \(\sub{0},\sub{1},\dotsc,\sub{2n+2}\) are \tH{} as well.
 \rremss{H.R.D*}{F.R.Dd*} yield \(\Dlau{n}^\ad=\Drau{n}\), \(\Dlub{n}^\ad=\Drub{n}\), and \(\rk{\Dloou{\saus}{\subs}{n+1}}^\ad=\Droou{\saus}{\subs}{n+1}\).
 Taking additionally into account \(\dia{1}\in\Cggq\) and \(\ba>0\), we conclude, in view of \rremss{A.R.XAX}{A.R.kK}, form \rpropss{F.P.FTH}{F.P.FTHab} that the matrices \(\Huo{n+1}{t}\) and \(-\ug\obg\Huo{n}{t}+(\ug+\obg)\Kuo{n}{t}-\Guo{n}{t}\) are both \tnnH{}, \tie{}, \(\seqt{2n+2}\) belongs to \(\Fggqu{2n+2}\).
  
 We thus have proved \(\seqt{m}\in\Fggqu{m}\) for all \(m\in\mn{0}{\kappa-1}\).
 In particular, \(\seqt{\kappa-1}\in\Fggqu{\kappa-1}\) holds true.
 Because of \rrem{F.R.FT1}, hence \(\seq{\su{j}^\FT}{j}{0}{\kappa-1}\in\Fggqu{\kappa-1}\).
 In view of the recursive construction of the \tnFT{k}, we now can proceed our considerations by mathematical induction to complete the proof.
\eproof

 We are now going to prove a relation similar to \rthm{H.P.T915}  between the \tfp{s} of an \tFnnd{} sequence \(\seqska\) and its associated \tFT{s} \(\seq{\su{j}^\FTa{0}}{j}{0}{\kappa},\seq{\su{j}^\FTa{1}}{j}{0}{\kappa-1},\seq{\su{j}^\FTa{2}}{j}{0}{\kappa-2},\dotsc\)
 To that end, block \(LDU\)~factorizations of the corresponding \tbHms{} \(\Hu{n}\), \(\Hau{n}\), \(\Hub{n}\), and \(\Hab{n}\) are needed.
 By a block \(LDU\)~factorization of a given complex \taaa{(n+1)p}{(n+1)q}{block} matrix \(M\), we mean a representation \(M=LDU\) where \(L\in\nudpu{n}\) (\tresp{}\ \(U\in\nodqu{n}\)) is a lower (\tresp{}\ upper) block triangular factor belonging to the corresponding class introduced in \rnota{ab.N1308} and \(D\) is a complex \taaa{(n+1)p}{(n+1)q}{block} diagonal matrix with \tpqa{blocks}.
 To avoid unnecessary explicit computations, we will use the equivalence relation ``\(\ldu\)'' introduced in \rnota{A.N.ldusim}, omitting the left and right factors belonging to \(\nudqu{n}\) and \(\nodqu{n}\), \tresp{} 
 For each matrix \(X=X^{\ok{s}}\) built from the sequence \(\seqska\), we denote (if possible) by \(X^\FTa{k}\defeq X^{\ok{u}}\) the corresponding matrix built from the \tnFT{k} \(\seq{u_j}{j}{0}{\kappa-k}\) of \(\seqska\) instead of \(\seqska\).

\bleml{F.L.Hsim}
 Let \(\seqska\in\Fggqka\).
\benui
 \il{F.L.Hsim.a} If \(\kappa\geq3\), for all \(n\in\N\) with \(2n+1\leq\kappa\), then
\begin{align}\label{F.L.Hsim.A}
 \ba^2\Hau{n}
 &\ldu\ba^2\diag\rk{\sau{0},\LLau{n}}
 \ldu\ba\diag\rk{\ba\sau{0},\Hab{n-1}^\FTa{1}}
 \ldu\diag\rk{\ba^2\sau{0},\Hau{n-1}^\FTa{2}}
\shortintertext{and}\label{F.L.Hsim.B}
 \ba^2\Hub{n}
 &\ldu\ba^2\diag\rk{\sub{0},\LLub{n}}
 \ldu\ba\diag\rk{\ba\sub{0},\LLu{n}^\FTa{1}}
 \ldu\diag\rk{\ba^2\sub{0},\Hub{n-1}^\FTa{2}}.
\end{align}
 \il{F.L.Hsim.b} If \(\kappa\geq4\), for all \(n\in\N\) with \(2n+2\leq\kappa\), then
\begin{align}
 \ba^2\Hab{n}
 &\ldu\ba\Hau{n}^\FTa{1}
 \ldu\ba\diag\rk{\sau{0}^\FTa{1},\LLau{n}^\FTa{1}}
 \ldu\diag\rk{\ba\sau{0}^\FTa{1},\Hab{n-1}^\FTa{2}}\label{F.L.Hsim.C}
\shortintertext{and}
 \ba^2\LLu{n+1}
 &\ldu\ba\Hub{n}^\FTa{1}
 \ldu\ba\diag\rk{\sub{0}^\FTa{1},\LLub{n}^\FTa{1}}
 \ldu\diag\rk{\ba\sub{0}^\FTa{1},\LLu{n}^\FTa{2}}.\label{F.L.Hsim.D}
\end{align}
\eenui
\elem
\bproof
 \eqref{F.L.Hsim.a} Assume \(\kappa\geq3\) and consider an arbitrary \(n\in\N\) with \(2n+1\leq\kappa\).
 According to \rprop{ab.R0933}, we have \(\seqs{2n+1}\in\Fggqu{2n+1}\).
 In view of \rprop{F.R.Fgg<D}, in particular \(\set{\seqsa{2n},\seqsb{2n}}\subseteq\Dqqu{2n}\) holds true.
 Applying \rlem{H.L.HsL} to \(\seqsa{2n}\) and \(\seqsb{2n}\), we can conclude then \(\Hau{n}\ldu\diag\rk{\sau{0},\LLau{n}}\) and \(\Hub{n}\ldu\diag\rk{\sub{0},\LLub{n}}\), \tresp{}
 Using \rprop{F.P.FTHab} and \rcor{F.C.FTL}, we obtain furthermore \(\ba\LLau{n}\ldu\Hab{n-1}^\FTa{1}\) and \(\ba\LLub{n}\ldu\LLu{n}^\FTa{1}\), \tresp{}
 Because of \rthm{ab.P1030}, the sequence \(\seq{\su{j}^\FT}{j}{0}{2n}\) belongs to \(\Fggqu{2n}\).
 The application of \rpropss{F.P.FTHa}{F.P.FTHb} to \(\seq{\su{j}^\FT}{j}{0}{2n}\) thus yields \(\ba\Hab{n-1}^\FTa{1}\ldu\Hau{n-1}^\FTa{2}\) and \(\ba\LLu{n}^\FTa{1}\ldu\Hub{n-1}^\FTa{2}\), \tresp{}
 By virtue of \rrem{A.R.ldudiag}, then \eqref{F.L.Hsim.A} and \eqref{F.L.Hsim.B} follow.
 
 \eqref{F.L.Hsim.b} Assume \(\kappa\geq4\) and consider an arbitrary \(n\in\N\) with \(2n+2\leq\kappa\).
 According to \rprop{ab.R0933}, we have \(\seqs{2n+2}\in\Fggqu{2n+2}\).
 Using \rpropss{F.P.FTHa}{F.P.FTHb}, we obtain thus \(\ba\Hab{n}\ldu\Hau{n}^\FTa{1}\) and \(\ba\LLu{n+1}\ldu\Hub{n}^\FTa{1}\), \tresp{}
 Because of \rthm{ab.P1030}, the sequence \(\seq{\su{j}^\FT}{j}{0}{2n+1}\) belongs to \(\Fggqu{2n+1}\).
 In view of \rprop{F.R.Fgg<D}, in particular \(\set{\seq{\sau{j}^\FT}{j}{0}{2n},\seq{\sub{j}^\FT}{j}{0}{2n}}\subseteq\Dqqu{2n}\) holds true.
 Applying \rlem{H.L.HsL} to \(\seq{\sau{j}^\FT}{j}{0}{2n}\) and \(\seq{\sub{j}^\FT}{j}{0}{2n}\), we can conclude then \(\Hau{n}^\FT\ldu\diag\rk{\sau{0}^\FT,\LLau{n}^\FT}\) and \(\Hub{n}^\FT\ldu\diag\rk{\sub{0}^\FT,\LLub{n}^\FT}\), \tresp{}
 The application of \rprop{F.P.FTHab} and \rcor{F.C.FTL} to \(\seq{\su{j}^\FT}{j}{0}{2n+1}\) yields furthermore \(\ba\LLau{n}^\FTa{1}\ldu\Hab{n-1}^\FTa{2}\) and \(\ba\LLub{n}^\FT\ldu\LLu{n}^\FTa{2}\), \tresp{}
 By virtue of \rrem{A.R.ldudiag}, then \eqref{F.L.Hsim.C} and \eqref{F.L.Hsim.D} follow.
\eproof

 By repeated application of \rlem{F.L.Hsim}, we are able now to reduce each of the four \tbHms{} \(\Hu{n}\), \(\Hau{n}\), \(\Hub{n}\), and \(\Hab{n}\) to block diagonal form, up to equivalence: 

\bleml{F.L.HsimFT0}
 Let \(n\in\N\) and let \(\seqs{2n}\in\Fggqu{2n}\).
 Then
\[
 \Hu{n}
 \ldu\diag\rk{\su{0},\ba^{-1}\sub{0}^\FTa{1},\ba^{-3}\sub{0}^\FTa{3},\dotsc,\ba^{-(2n-1)}\sub{0}^\FTa{2n-1}}
 \;\text{and}\;
 \Hab{n-1}\ldu\diag\seq{\ba^{-(2j+1)}\sau{0}^\FTa{2j+1}}{j}{0}{n-1}.
 \]
\elem
\bproof
 In view of \rprop{F.R.Fgg<D}, the sequence \(\seqs{2n}\) belongs to \(\Dqqu{2n}\).
 Hence, using \rlem{H.L.HsL}, we infer
\[
 \Hu{n}
 \ldu\diag\rk{\su{0},\LLu{n}}.
\]
 According to \rthm{ab.P1030}, we have, for all \(m\in\mn{0}{2n}\), furthermore \(\seq{\su{j}^\FTa{m}}{j}{0}{2n-m}\in\Fggqu{2n-m}\).
 Applying \rpropss{F.P.FTHb}{F.P.FTHa} to the sequence \(\seq{\su{j}^\FTa{2(n-1)}}{j}{0}{2}\), we obtain \(\Hub{0}^\FTa{2n-1}\ldu\ba\LLu{1}^\FTa{2n-2}\) and \(\Hau{0}^\FTa{2n-1}\ldu\ba\Hab{0}^\FTa{2n-2}\), \tresp{}
 Hence,
\begin{align*}
 \ba^{-2(n-1)}\LLu{1}^\FTa{2(n-1)}&\ldu\ba^{-(2n-1)}\sub{0}^\FTa{2n-1}&
&\text{and}&
 \ba^{-2(n-1)}\Hab{0}^\FTa{2(n-1)}&\ldu\ba^{-(2n-1)}\sau{0}^\FTa{2n-1}
\end{align*}
 hold true.
 In view of \rrem{A.R.ldudiag}, if \(n=1\), the proof is complete.
 Now assume \(n\geq2\).
 For all \(k\in\mn{0}{n-2}\), the application of \rlemp{F.L.Hsim}{F.L.Hsim.b} to the sequence \(\seq{\su{j}^\FTa{2k}}{j}{0}{2(n-k)}\) then yields \(\ba^2\LLu{n-k}^\FTa{2k}\ldu\diag\rk{\ba\sub{0}^\FTa{2k+1},\LLu{n-k-1}^\FTa{2k+2}}\)
 and \(\ba^2\Hab{n-k-1}^\FTa{2k}\ldu\diag\rk{\ba\sau{0}^\FTa{2k+1},\Hab{n-k-2}^\FTa{2k+2}}\).
 Taking into account \(\ba>0\), consequently
\[
 \ba^{-2k}\LLu{n-k}^\FTa{2k}
 \ldu\diag\rk{\ba^{-(2k+1)}\sub{0}^\FTa{2k+1},\ba^{-2(k+1)}\LLu{n-(k+1)}^\FTa{2(k+1)}}
\]
 and
\[
 \ba^{-2k}\Hab{n-k-1}^\FTa{2k}
 \ldu\diag\rk{\ba^{-(2k+1)}\sau{0}^\FTa{2k+1},\ba^{-2(k+1)}\Hab{n-(k+1)-1}^\FTa{2(k+1)}}
\]
 follow.
 Using \rrem{A.R.ldudiag}, we thus can conclude the assertion.
\eproof

\bleml{F.L.HsimFT1}
 If \(n\in\NO\) and \(\seqs{2n+1}\in\Fggqu{2n+1}\), then
\begin{align*}
 \Hau{n}&\ldu\diag\seq{\ba^{-2j}\sau{0}^\FTa{2j}}{j}{0}{n}&
&\text{and}&
 \Hub{n}&\ldu\diag\seq{\ba^{-2j}\sub{0}^\FTa{2j}}{j}{0}{n}.
\end{align*}
\elem
\bproof
 Obviously, the assertion is valid for \(n=0\).
 Now assume \(n\geq1\).
 According to \rthm{ab.P1030}, we have \(\seq{\su{j}^\FTa{m}}{j}{0}{2n-m+1}\in\Fggqu{2n-m+1}\) for all \(m\in\mn{0}{2n+1}\).
 For all \(k\in\mn{0}{n-1}\), the application of \rlemp{F.L.Hsim}{F.L.Hsim.a} to the sequence \(\seq{\su{j}^\FTa{2k}}{j}{0}{2(n-k)+1}\) then yields \(\ba^2\Hau{n-k}^\FTa{2k}\ldu\diag\rk{\ba^2\sau{0}^\FTa{2k},\Hau{n-k-1}^\FTa{2k+2}}\) and
 \(\ba^2\Hub{n-k}^\FTa{2k}\ldu\diag\rk{\ba^2\sub{0}^\FTa{2k},\Hub{n-k-1}^\FTa{2k+2}}\).
 Taking into account \(\ba>0\), consequently
\[
 \ba^{-2k}\Hau{n-k}^\FTa{2k}
 \ldu\diag\rk{\ba^{-2k}\sau{0}^\FTa{2k},\ba^{-2(k+1)}\Hau{n-(k+1)}^\FTa{2(k+1)}}
\]
 and
\[
 \ba^{-2k}\Hub{n-k}^\FTa{2k}
 \ldu\diag\rk{\ba^{-2k}\sub{0}^\FTa{2k},\ba^{-2(k+1)}\Hub{n-(k+1)}^\FTa{2(k+1)}}.
\]
 Because of \(\Hau{0}^\FTa{2n}=\sau{0}^\FTa{2n}\) and \(\Hub{0}^\FTa{2n}=\sub{0}^\FTa{2n}\), the assertion follows from \rrem{A.R.ldudiag}.
\eproof

 Using \rlem{L3219}, we now obtain from \rlemss{F.L.HsimFT0}{F.L.HsimFT1} a relation between the \tfp{s} \(\fpseqka\) of a \tFnnd{} sequence \(\seqska\) and its \tnFT{k} \(\seq{\su{j}^\FTa{k}}{j}{0}{\kappa-k}\):

\bpropl{F.P.FPFT}
 Let \(\seqska\in\Fggqka\) with \tfpf{} \(\fpseqka\).
 Then \(\fpu{0}=\su{0}\) and furthermore
\begin{align*}
 \fpu{4k+1}&=\ba^{-2k}\sau{0}^\FTa{2k}&
&\text{and}&
 \fpu{4k+2}&=\ba^{-2k}\sub{0}^\FTa{2k}
\end{align*}
 for all \(k\in\NO\) with \(2k+1\leq\kappa\) and
\begin{align*}
 \fpu{4k+3}&=\ba^{-(2k+1)}\sau{0}^\FTa{2k+1}&
&\text{and}&
 \fpu{4k+4}&=\ba^{-(2k+1)}\sub{0}^\FTa{2k+1}
\end{align*}
 for all \(k\in\NO\) with \(2k+2\leq\kappa\).
\eprop
\bproof
 In view of \eqref{F.G.f012}, we have \(\fpu{0}=\su{0}\).
 
 Now assume \(\kappa\geq1\) and consider an arbitrary \(k\in\NO\) with \(2k+1\leq\kappa\).
 According to \rprop{ab.R0933}, the sequence \(\seqs{2k+1}\) belongs to \(\Fggqu{2k+1}\).
 In view of \eqref{Fgg2n+1}, thus \(\set{\seqsa{2k},\seqsb{2k}}\subseteq\Hggqu{2k}\).
 \rlem{F.L.HsimFT1} yields furthermore \(\Hau{k}\ldu\diag\seq{\ba^{-2j}\sau{0}^\FTa{2j}}{j}{0}{k}\) and \(\Hub{k}\ldu\diag\seq{\ba^{-2j}\sub{0}^\FTa{2j}}{j}{0}{k}\).
 Hence, we infer from \rlem{L3219} in particular \(\Lau{k}=\ba^{-2k}\sau{0}^\FTa{2k}\) and \(\Lub{k}=\ba^{-2k}\sub{0}^\FTa{2k}\).
 By virtue of \rrem{F.R.ABL}, then \(\fpu{4k+1}=\usc{2k+1}=\ba^{-2k}\sau{0}^\FTa{2k}\) and \(\fpu{4k+2}=\osc{2k+1}=\ba^{-2k}\sub{0}^\FTa{2k}\) follow.
 
 Now assume \(\kappa\geq2\) and consider an arbitrary \(k\in\NO\) with \(2k+2\leq\kappa\).
 According to \rprop{ab.R0933}, the sequence \(\seqs{2(k+1)}\) belongs to \(\Fggqu{2(k+1)}\).
 In view of \eqref{Fgg2n}, thus \(\seqsab{2k}\in\Hggqu{2k}\) and \(\seqs{2(k+1)}\in\Hggqu{2(k+1)}\).
 \rlem{F.L.HsimFT0} yields furthermore \(\Hab{k}\ldu\diag\seq{\ba^{-(2j+1)}\sau{0}^\FTa{2j+1}}{j}{0}{k}\) and \(\Hu{k+1}\ldu\diag\rk{\su{0},\ba^{-1}\sub{0}^\FTa{1},\ba^{-3}\sub{0}^\FTa{3},\dotsc,\ba^{-(2k+1)}\sub{0}^\FTa{2k+1}}\).
 Hence, from \rlem{L3219} we infer in particular \(\Lab{k}=\ba^{-(2k+1)}\sau{0}^\FTa{2k+1}\) and \(\Lu{k+1}=\ba^{-(2k+1)}\sub{0}^\FTa{2k+1}\).
 By virtue of \rrem{F.R.ABL}, then \(\fpu{4k+3}=\osc{2k+2}=\ba^{-(2k+1)}\sau{0}^\FTa{2k+1}\) and \(\fpu{4k+4}=\usc{2k+2}=\ba^{-(2k+1)}\sub{0}^\FTa{2k+1}\) follow.
\eproof

 Using \rprop{F.P.FPFT}, we are now able to express the \tfdf{} by the \tFT{s} of an \tFnnd{} sequence:
 
\bcorl{F.C.diaFT}
 If \(\seqska\in \Fggqka \), then \(\dia{j}=\ba^{-(j-1)}\su{0}^\FTa{j}\) for all \(j\in\mn{0}{\kappa}\).
\ecor
\bproof
 In the case \(\kappa<\infp\), we first extend the sequence \(\seqska\), by virtue of \rthm{165.T112}, to a sequence \(\seqs{\kappa+1}\in\Fggqu{\kappa+1}\).
 In view of \rremss{F.R.diatr}{F.R.kFTtr}, then the assertion follows from \rremss{F.R.f2n-1}{F.R.012} and \rprop{F.P.FPFT}.
\eproof

 Using the parallel sum given via \eqref{ps}, the effect caused by \tFT{ation} on the \tfp{s} can be completely described:

\bcorl{F.C.Salg1}
 Assume \(\kappa\geq1\) and let \(\seqska\in\Fggqka\) with \tFT{} \(\seqt{\kappa-1}\).
 Denote by \(\fgpseq{2(\kappa-1)}\) the \tfpfa{\(\seqt{\kappa-1}\)}.
 Then \(\fgpu{0}=\ba\rk{\fpu{1}\ps\fpu{2}}\) and \(\fgpu{j}=\ba\fpu{j+2}\) for all \(j\in\mn{1}{2(\kappa-1)}\).
\ecor
\bproof
 According to \eqref{F.G.f012}, we have \(\fgpu{0}=t_0\).
 Because of \rprop{F.R.Fgg<D}, the sequence \(\seqska\) belongs to \(\Dqqu{\kappa}\).
 Hence, we infer \(t_0=\dia{1}\) from \rlem{ab.P1728a}.
 By virtue of \rprop{F.R.d+=f}, furthermore \(\dia{1}=\ba\rk{\fpu{1}\ps\fpu{2}}\) holds true.
 Consequently, we get \(\fgpu{0}=\ba\rk{\fpu{1}\ps\fpu{2}}\).
 \rthm{ab.P1030} yields \(\seqt{\kappa-1}\in\Fggqu{\kappa-1}\).
 Taking into account the recursive construction of the \tnFT{k}, we can then conclude \(\fgpu{j}=\ba\fpu{j+2}\) for all \(j\in\mn{1}{2(\kappa-1)}\) from \rprop{F.P.FPFT}.
\eproof

 We are now going to show that the \tFT{ation} of an \tFnnd{} sequence essentially coincides with left shifting of the corresponding \tfpf{} or the \tfcf{}.
 By virtue of the recursive construction of the \tnFT{k}, we infer from \rcor{F.C.diaFT} immediately:
 
\bpropl{F.P.diaalg}
 Let \(k\in\mn{0}{\kappa}\) and let \(\seqska\in\Fggqka \) with \tfdf{} \(\seqdiaka\) and \tnFT{k} \(\seq{\su{j}^\FTa{k}}{j}{0}{\kappa-k}\).
 Then \(\seq{\ba^k\dia{k+j}}{j}{0}{\kappa-k}\) coincides with the \tfdfa{\(\seq{\su{j}^\FTa{k}}{j}{0}{\kappa-k}\)}.
\eprop

\bcorl{F.C.abCDalg}
 Let \(k,\ell\in\mn{0}{\kappa}\) and let \(\seqska\in\Fggqka\) be \tabCDo{\ell}.
 Then \(\seq{\su{j}^\FTa{k}}{j}{0}{\kappa-k}\) is \tabCDo{\max\set{0,\ell-k}}.
\ecor
\bproof
 According to \rdefn{F.D.abCD} and \rprop{ab.C1101}, we have \(\dia{j}=\Oqq\) for all \(j\in\mn{\ell}{\kappa}\).
 By virtue of \rprop{F.P.diaalg} and \rdefn{F.D.abCD}, then the assertion follows.
\eproof

\bthml{F.T.SalgFP}
 Let \(\seqska\in\Fggqka \) with \tfpf{} \(\fpseqka\) and \tfdf{} \(\seqdiaka \).
 Let \(k\in\mn{0}{\kappa}\) and denote by \(\seq{\su{j}^\FTa{k}}{j}{0}{\kappa-k}\) the \tnFTv{k}{\(\seqska\)} and by \(\fgpseq{2(\kappa-k)}\) the \tfpfa{\(\seq{\su{j}^\FTa{k}}{j}{0}{\kappa-k}\)}.
 Then \(\fgpu{0}=\ba^{k-1}\dia{k}\) and \(\fgpu{j}=\ba^k\fpu{2k+j}\) for all \(j\in\mn{1}{2(\kappa-k)}\).
\ethm
\bproof
 In view of \eqref{F.G.f012} and \eqref{F.G.d01}, the assertions are obviously valid for \(k=0\).
 Now assume \(\kappa\geq1\) and \(k\in\mn{1}{\kappa}\).
 Denote by \(\seqt{\kappa-(k-1)}\) the \tnFTv{(k-1)}{\(\seqska\)} and by \(\seq{\mathfrak{l}_j}{j}{0}{\kappa-(k-1)}\) the \tfdfa{\(\seqt{\kappa-(k-1)}\)}.
 According to \rthm{ab.P1030}, we have \(\seq{\su{j}^\FTa{m}}{j}{0}{\kappa-m}\in\Fggqu{\kappa-m}\) for all \(m\in\mn{0}{\kappa}\).
 In particular, \(\seqt{\kappa-(k-1)}\in\Fggqu{\kappa-(k-1)}\).
 From \rprop{F.R.Fgg<D} then \(\seqt{\kappa-(k-1)}\in\Dqqu{\kappa-(k-1)}\) follows.
 Hence, \rlem{ab.P1728a} yields \(t_0^\FT=\mathfrak{l}_1\).
 Using \rprop{F.P.diaalg}, we conclude \(\mathfrak{l}_1=\ba^{k-1}\dia{k}\).
 In view of \(t_0^\FT=\su{0}^\FTa{k}=\fgpu{0}\), we obtain thus \(\fgpu{0}=\ba^{k-1}\dia{k}\).
 By virtue of the recursive construction of the \tnFT{k}, the repeated application of \rcor{F.C.Salg1} provides \(\fgpu{j}=\ba^k\fpu{2k+j}\) for all \(j\in\mn{1}{2(\kappa-k)}\).
\eproof

\bthml{F.T.SalgIP}
 Let \(\seqska\in\Fggqu{\kappa}\) with \tfdf{} \(\seqdiaka\) and \tfcf{} \(\seqciaka \).
 Let \(k\in\mn{0}{\kappa}\) and denote by \(\seq{\su{j}^\FTa{k}}{j}{0}{\kappa-k}\) the \tnFTv{k}{\(\seqska\)}.
 Then \(\seq{\su{j}^\FTa{k}}{j}{0}{\kappa-k}\) belongs to \(\Fggqu{\kappa-k}\) and the \tfcf{} \(\seq{\mathfrak{p}_{j}}{j}{0}{\kappa-k}\) of \(\seq{\su{j}^\FTa{k}}{j}{0}{\kappa-k}\) fulfills \(\mathfrak{p}_{0}=\ba^{k-1}\dia{k}\) and \(\mathfrak{p}_{j}=\cia{k+j}\) for all \(j\in\mn{1}{\kappa-k}\).
\ethm
\bproof
 \rthm{ab.P1030} yields \(\seq{\su{j}^\FTa{k}}{j}{0}{\kappa-k}\in\Fggqu{\kappa-k}\).
 Denote by \(\fgpseq{2(\kappa-k)}\) the \tfpfa{\(\seq{\su{j}^\FTa{k}}{j}{0}{\kappa-k}\)}.
 In view of \rthm{F.T.SalgFP}, we have \(\mathfrak{p}_{0}=\fgpu{0}=\ba^{k-1}\dia{k}\).
 From \rprop{F.P.diaalg}, \rthm{F.T.SalgFP}, and \rremss{A.R.l*A}{A.R.A-1}, we get, for all \(j\in\mn{1}{\kappa-k}\), furthermore
\[\begin{split}
 \mathfrak{p}_{j}
 &=\ek*{\rk{\ba^k\dia{(k+j)-1}}^\varsqrt}^\mpi\fgpu{2j}\ek*{\rk{\ba^k\dia{(k+j)-1}}^\varsqrt}^\mpi
 =\rk{\ba^{k/2}\dia{(k+j)-1}^\varsqrt}^\mpi\rk{\ba^k\fpu{2k+2j}}\rk{\ba^{k/2}\dia{(k+j)-1}^\varsqrt}^\mpi\\
 &=\ek*{\ba^{-k/2}\rk{\dia{(k+j)-1}^\varsqrt}^\mpi}\rk{\ba^k\fpu{2(k+j)}}\ek*{\ba^{-k/2}\rk{\dia{(k+j)-1}^\varsqrt}^\mpi}
 =\cia{k+j}.\qedhere
\end{split}\]
\eproof

\bcorl{T1226}
 If \(\seqska\in\Fgqka\), then \(\seq{\su{j}^\FTa{k}}{j}{0}{\kappa-k}\in\Fgqu{\kappa-k}\) for all \(k\in\mn{0}{\kappa}\).
\ecor
\bproof
 In view of \(\ba>0\) and \rprop{ab.C0929}, we have \(\ba^{k-1}\dia{k}\in\Cgq\) for all \(k\in\mn{0}{\kappa}\).
 The combination of \rthm{F.T.SalgIP} and~\zitaa{MR3979701}{\cprop{6.36}{32}} completes the proof.
\eproof

\bcorl{F.C.abZalg}
 Assume \(\kappa\geq1\), let \(k\in\mn{0}{\kappa-1}\), let \(\ell\in\mn{1}{\kappa}\), and let \(\seqska\in\Fggqka\) be \tabZo{\ell}.
 Then \(\seq{\su{j}^\FTa{k}}{j}{0}{\kappa-k}\) is \tabZo{\max\set{1,\ell-k}}.
\ecor
\bproof
 According to \rprop{F.P.abZe} and \rnota{F.N.Pd}, we have \(\cia{j}=\frac{1}{2}\Pd{\ell-1}\) for all \(j\in\mn{\ell}{\kappa}\).
 Denote by \(\seq{\mathfrak{l}_j}{j}{0}{\kappa-k}\) the \tfdfa{\(\seq{\su{j}^\FTa{k}}{j}{0}{\kappa-k}\)} and by \(\seq{\mathfrak{p}_{j}}{j}{0}{\kappa-k}\) the \tfcfa{\(\seq{\su{j}^\FTa{k}}{j}{0}{\kappa-k}\)}.
 \rprop{F.P.diaalg} yields \(\mathfrak{l}_j=\ba^k\dia{k+j}\) for all \(j\in\mn{0}{\kappa-k}\).
 In view of \(\ba>0\), hence \(\OPu{\ran{\mathfrak{l}_j}}=\OPu{\ran{\dia{k+j}}}=\Pd{k+j}\) for all \(j\in\mn{0}{\kappa-k}\).
 From \rthm{F.T.SalgIP} we infer \(\mathfrak{p}_j=\cia{k+j}\) for all \(j\in\mn{1}{\kappa-k}\).
 For all \(m\in\N\) with \(\ell-k\leq m\leq\kappa-k\), we have then
\(
 \mathfrak{p}_m
 =\cia{k+m}
 =\frac{1}{2}\Pd{k+m-1}
 =\frac{1}{2}\OPu{\ran{\mathfrak{l}_{m-1}}}
\).
 Since \(\seq{\su{j}^\FTa{k}}{j}{0}{\kappa-k}\) belongs to \(\Fggqu{\kappa-k}\), in view of \rthm{ab.P1030}, the application of \rprop{F.P.abZe} completes the proof.
\eproof

\section{A Schur type transformation for matrix measures on \(\ab\)}\label{F.CM}
 As a main result of this last section we are going to characterize for \tnnH{} measures on the compact interval \(\ab\) of the real axis a centrality property associated to \rdefn{F.D.abZ}.
 To that end we realize that centrality in this sense is equivalent to obtaining essentially a matricial version of the arcsine distribution after transforming the given \tnnH{} measure several times in accordance with the algorithm considered in \rsec{S.F.kFT} applied to the corresponding moment sequence.
 Observe that the following construction is well defined due to \rprop{I.P.ab8} and \rthm{ab.P1030}:

\bdefnl{F.D.SN-M}
 Let \(\sigma\in\MggqF\) with \tfpmf{} \(\seqsinf\) and let \(k\in\NO\).
 Denote by \(\seq{\su{j}^\FTa{k}}{j}{0}{\infi}\) the \tnFTv{k}{\(\seqsinf\)} and by \(\sigma^\FTa{k}\) the uniquely determined element in \(\MggqFag{s^\FTa{k}}{\infi}\).
 Then we call \(\sigma^\FTa{k}\) the \emph{\tnFmTv{k}{\(\sigma\)}}.
\edefn 

 Against to the background of \rprop{I.P.ab8Fgg} and \rthm{F.T.Fggcia}, we generalized in~\zitaa{MR3979701}{\csec{8}} several results from the scalar theory of canonical moments to the matrix case:

\bdefnl{F.D.meacia}
 Let \(\sigma\in\MggqF\) with \tfpmf{} \(\seqmpm{\sigma}\).
 Denote by \(\seqmcm{\sigma}\) the \tfcfa{\(\seqmpm{\sigma}\)} and by \(\seqmdm{\sigma}\) the \tfdfa{\(\seqmpm{\sigma}\)}.
 Then we call \(\seqmcm{\sigma}\) the \emph{\tfmcmfa{\(\sigma\)}} and \(\seqmdm{\sigma}\) the \emph{\tfmdmfa{\(\sigma\)}}.
\edefn

\bthmnl{\zitaa{MR3979701}{\cthm{8.2}{36}}}{F.T.Mabcia}
 The mapping \(\Pi_\ab\colon\MggqF\to\es{q}{\infi}{\bam}\) given by \(\sigma\mapsto\seqmcm{\sigma}\) is well defined and bijective.
\ethm

 Regarding \rprop{I.P.ab8Fgg}, we have:
 
\bpropl{F.R.mdealg}
 Let \(k\in\NO\) and let \(\sigma\in\MggqF\) with \tnFmT{k} \(\mu\).
 Then \(\mcm{\mu}{0}=\ba^{k-1}\mdm{\sigma}{k}\) and \(\mcm{\mu}{j}=\mcm{\sigma}{k+j}\) for all \(j\in\N\).
 Furthermore, \(\mdm{\mu}{j}=\ba^k\mdm{\sigma}{k+j}\) for all \(j\in\NO\).
 In particular, \(\mu\rk{\ab}=\ba^{k-1}\mdm{\sigma}{k}\).
\eprop
\bproof
 Use \rthm{F.T.SalgIP} and \rprop{F.P.diaalg} together with \eqref{I.G.mom} and \eqref{F.G.d01}.
\eproof

 Let \(\Omega\in\BsAR\setminus\set{\emptyset}\).
 Then let \(\MggqmolO\) be the set of all \(\sigma\in\MggqO\) for which there exists a finite subset \(B\) of \(\Omega\) satisfying \(\sigma\rk{\Omega\setminus B}=\Oqq\).
 The elements of \(\MggqmolO\) are said to be \emph{molecular}.
 Obviously, \(\MggqmolO\) is the set of all \(\sigma\in\MggqO\) for which there exist an \(m\in\N\) and sequences \(\seq{\xi_\ell}{\ell}{1}{m}\) and \(\seq{A_\ell}{\ell}{1}{m}\) from \(\Omega\) and \(\Cggq\), \tresp{}, such that \(\sigma=\sum_{\ell=1}^m\Kronu{\xi_\ell}A_\ell\), where \(\Kronu{\xi_\ell}\) is the Dirac measure on \(\rk{\Omega,\BsAu{\Omega}}\) with unit mass at \(\xi_\ell\).
 In particular, \(\MggqmolO\subseteq\Mggoua{q}{\infi}{\Omega}\).
 
 It was shown in~\zitaa{MR3979701}{\cprop{8.4}{37}} that \(\sigma\in\MggqF\) is molecular if and only if \(\mdm{\sigma}{k}=\Oqq\) for some \(k\in\N\), which, in view of \rdefn{F.D.abCD}, is equivalent to \(\seqmpm{\sigma}\) being \tabCDo{k}.
 This leads to a characterization in terms of the \tnFmT{k}s:

\bpropl{F.P.algmol}
 Let \(\sigma\in\MggqF\).
 Then \(\sigma\) is molecular if and only if \(\sigma^\FTa{k}\) coincides with the \tqqa{zero} measure on \(\rk{\ab,\BsAF}\) for some \(k\in\NO\).
\eprop
\bproof
 First observe that \rprop{F.R.mdealg} yields \(\sigma^\FTa{\ell}\rk{\ab}=\ba^{\ell-1}\mdm{\sigma}{\ell}\) for all \(\ell\in\NO\).
 Taking into account \(\ba>0\),  for all \(\ell\in\NO\) we have then \(\sigma^\FTa{\ell}\rk{\ab}=\Oqq\) if and only if \(\mdm{\sigma}{\ell}=\Oqq\).
 
 If \(\sigma\) is molecular, then~\zitaa{MR3979701}{\cprop{8.4}{37}} yields \(\mdm{\sigma}{k}=\Oqq\) for some \(k\in\N\), implying \(\sigma^\FTa{k}\rk{\ab}=\Oqq\), \tie{}\ \(\sigma^\FTa{k}\) is the \tqqa{zero} measure on \(\rk{\ab,\BsAF}\).
 
 Conversely, assume that there is some \(k\in\NO\) such that \(\sigma^\FTa{k}\) is the \tqqa{zero} measure on \(\rk{\ab,\BsAF}\).
 If \(k=0\), then \(\sigma\) is molecular, since \(\sigma=\sigma^\FTa{0}\).
 Now assume \(k\geq1\).
 Then \(\sigma^\FTa{k}\rk{\ab}=\Oqq\), implying \(\mdm{\sigma}{k}=\Oqq\).
 From~\zitaa{MR3979701}{\cprop{8.4}{37}} we can conclude then that  \(\sigma\) is molecular. 
\eproof

 The notion of centrality in the sense of \rdefn{F.D.abZ} carries over to \tnnH{} measures on \(\rk{\ab,\BsAF}\) in an obvious way:
 
\bdefnl{F.D.Zmo}
 Let \(\sigma\in\MggqF\) with \tfpmf{} \(\seqmpm{\sigma}\) and let \(k\in\N\).
 We call \(\sigma\) \emph{\tZmo{k}} if \(\seqmpm{\sigma}\) is \tabZo{k}.
\edefn

 We are going to characterize \tZm{}ity in terms of the \tnFmT{k}s.
 First we supplement~\zitaa{MR3979701}{\clem{8.11}{40}}:

\bleml{F.L.sM}
 Let \(M\in\Cggq\) and let \(\nu\in\Mggoa{1}{\ab}\).
 Then \(\sigma\colon\BsAF\to\Cggq\) defined by \(\sigma\rk{B}\defeq\ek{\nu\rk{B}}M\) belongs to \(\MggqF\) and has \tfpmf{} \(\seq{\mpm{\nu}{j}M}{j}{0}{\infi}\) and \tfmdmf{} \(\seq{\mdm{\nu}{j}M}{j}{0}{\infi}\).
 Furthermore, \(\mcm{\sigma}{0}=\mcm{\nu}{0}M\) and \(\mcm{\sigma}{j}=\mcm{\nu}{j}\OPu{\ran{M}}\) for all \(j\in\N\).
\elem
\bproof
 This follows from \eqref{I.G.mom}, \rrem{F.R.dxB}, and the combination of \rprop{I.P.ab8Fgg} with \rlem{F.L.exB}.
\eproof

 Let \(\Leb\colon\BsAR\to[0,\infp)\) be the Lebesgue measure.
 For each \(\Omega\in\BsAR\setminus\set{\emptyset}\), we denote by \(\Lebu{\Omega}\) the restriction of \(\Leb\) to \(\BsAO\).
 Furthermore, for all \(a,b\in\R\) with \(a<b\), let \(f_{[a,b]}\colon[a,b]\to[0,\infp)\) be defined by \(f_{[a,b]}\rk{x}\defeq0\) if \(x\in\set{a,b}\) and by \(f_{[a,b]}\rk{x}\defeq\ek{\pi\sqrt{\rk{x-a}\rk{b-x}}}^\inv\) if \(x\in(a,b)\). 

\bexaml{F.E.arcsin}
 Let \(\nu\colon\BsAu{[0,1]}\to[0,1]\) be the \emph{arcsine distribution} given by \(\nu\rk{A}\defeq\int_Af_{[0,1]}\dif\Lebu{[0,1]}\). 
 It is a probability distribution having (classical) canonical moments \(p_k=1/2\) for all \(k\in\N\) (\tcf{}~\zitaa{MR1468473}{\cexa{1.3.6}{15}}).
 By virtue of \eqref{F.G.epq}, then its \tfmcmf{} in the sense of \rdefn{F.D.meacia} fulfills \(\mcm{\nu}{0}=1\) and \(\mcm{\nu}{j}=1/2\) for all \(j\in\N\).
 Let \(T\colon[0,1]\to\ab\) be defined by \(T(x)\defeq\ba x+\ug\) and denote by \(\tau\) the image measure of \(\nu\) under \(T\).
 Then \(\tau\rk{B}=\int_B f_{\ab}\dif\Lebu{\ab}\) for all \(B\in\BsAF\).
 In view of \(\ba>0\), furthermore the \tfmcmfa{\(\tau\)} coincides, according to~\zitaa{MR3979701}{\cprop{8.12}{40}}, with \(\seqmcm{\nu}\). 
 Let \(M\in\Cggq\).
 Because of \rlem{F.L.sM}, then \(\mu\colon\BsAF\to\Cggq\) defined by \(\mu\rk{B}\defeq\ek{\tau\rk{B}}M\) belongs to \(\MggqF\) and  its \tfmcmf{} fulfills \(\mcm{\mu}{0}=\mcm{\tau}{0}M=M\) and \(\mcm{\mu}{j}=\mcm{\tau}{j}\OPu{\ran{M}}=\frac{1}{2}\OPu{\ran{M}}\) for all \(j\in\N\).
 (This sequence already occurred in \rexam{F.E.Z1}.)
\eexam

\bthml{F.P.algZmo}
 Let \(\sigma\in\MggqF\) and let \(k\in\N\).
 Then \(\sigma\) is \tZmo{k} if and only if \(\sigma^\FTa{k-1}\) fulfills
\begin{align}\label{F.P.algZmo.A}
 \sigma^\FTa{k-1}\rk{B}&
 =\ba^{k-2}\rk*{\int_Bf_{\ab}\dif\Lebu{\ab}}\mdm{\sigma}{k-1}&\text{for all }B&\in\BsAF.
\end{align}
\ethm
\bproof
 Setting \(\mu\defeq\sigma^\FTa{k-1}\) and \(M\defeq\mu\rk{\ab}\), \rprop{F.R.mdealg} yields \(M=\ba^{k-2}\mdm{\sigma}{k-1}\).
 Denote by \(\seqsinf\) the \tfpmfa{\(\sigma\)}.
 Then the \tnFT{{\rk{k-1}}} \(\seq{\su{j}^\FTa{k-1}}{j}{0}{\infi}\) of \(\seqsinf\) coincides, by definition, with the \tfpmfa{\(\mu\)}.
 In particular, \(\su{0}^\FTa{k-1}=M\) and the \tfmcmf{} \(\seqmcm{\mu}\) associated with \(\mu\) is exactly the \tfcfa{\(\seq{\su{j}^\FTa{k-1}}{j}{0}{\infi}\)}.
 According to \rprop{I.P.ab8Fgg}, the sequences \(\seqsinf\) and \(\seq{\su{j}^\FTa{k-1}}{j}{0}{\infi}\) both belong to \(\Fggqinf\).
 
 First suppose that \(\sigma\) is \tZmo{k}, \tie{}, \(\seqsinf\) is \tabZo{k}.
 In view of \rcor{F.C.abZalg}, hence \(\seq{\su{j}^\FTa{k-1}}{j}{0}{\infi}\) is \tabZo{1}.
 Consequently, from \rexam{F.E.Z1} we get \(\mcm{\mu}{0}=\su{0}^\FTa{k-1}=M\) and \(\mcm{\mu}{j}=\frac{1}{2}\OPu{\ran{M}}\) for all \(j\in\N\).
 Taking into account \rthm{F.T.Mabcia}, we can conclude from \rexam{F.E.arcsin} for all \(B\in\BsAF\) then \(\mu\rk{B}=\rk{\int_Bf_{\ab}\dif\Lebu{\ab}}M\), implying \eqref{F.P.algZmo.A} by \(M=\ba^{k-2}\mdm{\sigma}{k-1}\).
 
 Conversely suppose that \eqref{F.P.algZmo.A} holds true.
 Since \(M\) is \tnnH{} by construction, we can conclude from \rexam{F.E.arcsin} then \(\mcm{\mu}{0}=M\) and \(\mcm{\mu}{j}=\frac{1}{2}\OPu{\ran{M}}\) for all \(j\in\N\).
 \rprop{F.R.mdealg} yields \(\mcm{\mu}{j}=\mcm{\sigma}{k-1+j}\) for all \(j\in\N\).
 Because of \(\ba>0\), for all \(\ell\in\minf{k}\), we obtain consequently
\(
 \mcm{\sigma}{\ell}
 =\mcm{\mu}{\ell-k+1}
 =\frac{1}{2}\OPu{\ran{M}}
 =\frac{1}{2}\OPu{\ran{\mdm{\sigma}{k-1}}}
\).
 According to \rprop{F.P.abZe}, then \(\seqsinf\) is \tabZo{k}, \tie{}, \(\sigma\) is \tZmo{k}.
\eproof

\appendix
\section{Some facts from matrix theory}\label{S.A}
 This appendix contains a summary of results from matrix theory, which are used in this paper.
 What concerns results on the Moore--Penrose inverse \(A^\mpi\) of a complex matrix \(A\), we refer to~\zitaa{MR1152328}{\csec{1}}.
 Results on the Kronecker product \(A\kp B\) of complex matrices \(A\) and \(B\) can be found, \teg{}, in~\zitaa{MR1288752}{\csec{4.2}}.

\breml{A.R.AB=BA}
 Let \(m,n\in\N\) and let \(A_1,A_2,\dotsc,A_m\) and \(B_1,B_2,\dotsc,B_n\) be complex \tpqa{matrices}.
 If \(M\in\Cqp\) is such that \(A_jMB_k=B_kMA_j\) holds true for all \(j\in\mn{1}{m}\) and all \(k\in\mn{1}{n}\), then \(\rk{\sum_{j=1}^m\eta_jA_j}M\rk{\sum_{k=1}^n\theta_kB_k}=\rk{\sum_{k=1}^n\theta_kB_k}M\rk{\sum_{j=1}^m\eta_jA_j}\) for all complex numbers \(\eta_1,\eta_2,\dotsc,\eta_m\) and \(\theta_1,\theta_2,\dotsc,\theta_n\).
\erem

\breml{A.R.rs+}
 If \(n\in\N\) and \(A_1,A_2,\dotsc,A_n\in\Cpq\), then \(\ran{\sum_{j=1}^n\eta_jA_j}\subseteq\sum_{j=1}^n\ran{A_j}\) and \(\bigcap_{j=1}^n\nul{A_j}\subseteq\nul{\sum_{j=1}^n\eta_jA_j}\) for all \(\eta_1,\eta_2,\dotsc,\eta_n\in\C\).
\erem

\breml{A.R.wxB}
 Let \(\eta\in\C\) and let \(B\in\Cpq\).
 Denote by \(A\) the complex \taaa{1}{1}{matrix} with entry \(\eta\).
 Then \(A\kp B=\eta B\).
\erem

\breml{A.R.AxB*CxD}
 If \(A\in\Coo{m}{n}\), \(B\in\Cpq\), \(C\in\Coo{n}{s}\), and \(D\in\Coo{q}{v}\), then \(\rk{A\kp B}\rk{C\kp D}=\rk{AC}\kp\rk{BD}\).
\erem

\breml{A.R.A++*}
 If \(A\in\Cpq\), then \(\rk{A^\mpi}^\mpi=A\) and \(\rk{A^\ad}^\mpi=\rk{A^\mpi}^\ad\).
\erem

\breml{A.R.l*A}
 If \(\eta\in\C\) and \(A\in\Cpq\), then \(\rk{\eta A}^\mpi=\eta^\mpi A^\mpi\).
\erem

\breml{A.R.UA+V}
 Let \(U\in\Coo{u}{p}\) with \(U^\ad U=\Ip\), let \(V\in\Coo{q}{v}\) with \(VV^\ad=\Iq\), and let \(A\in\Cpq\).
 Then \(\rk{UAV}^\mpi=V^\ad A^\mpi U^\ad\).
\erem

\breml{A.R.diag+}
 If \(D\in\Coo{p}{r}\) and \(E\in\Coo{q}{s}\), then \(\rk{\zdiag{D}{E}}^\mpi=\zdiag{D^\mpi}{E^\mpi}\).
\erem

\bremnl{\tcf{}~\zitaa{MR1987382}{Ex.~3, \cpage{54}}}{A.R.AxB+}
 If \(A\in\Coo{m}{n}\) and \(B\in\Cpq\), then \(\rk{A\kp B}^\mpi=A^\mpi\kp B^\mpi\).
\erem

\breml{A.R.A-1}
 If \(A\in\Cqq\) fulfills \(\det A\neq0\), then \(A^\mpi=A^\inv\).
\erem

\breml{ab.R1052}
 If \(A\in\Cpq\), then \(\OPu{\ran{A}}=AA^\mpi\) and \(\OPu{\ran{A^\ad}}=A^\mpi A\).
\erem

\breml{R.AA+B}
 If \(A\in\Cpq\), then:
\benui
 \il{R.AA+B.a} Let \(B\in\Coo{p}{m}\).
 Then \(\ran{B}\subseteq\ran{A}\) if and only if \(AA^\mpi B=B\).
 \item Let \(C\in\Coo{n}{q}\).
 Then \(\nul{A}\subseteq\nul{C}\) if and only if \(CA^\mpi A=C\).
\eenui
\erem

\breml{A.R.kK}
 The set \(\Cggq\) is a convex cone in the \(\R\)\nobreakdash-vector space \(\CHq\).
\erem

\breml{A.R.XAX}
 If \(A\in\Cggq\) and \(X\in\Cqp\), then \(X^\ad AX\in\Cggp\).
\erem

\bremnl{\tcf{}~\zitaa{MR1152328}{\clem{1.1.9(b)}{18}}}{L.AEP}
 If \(M=\tmat{A&B\\C&D}\) is the \tbr{} of a \tnnH{} \taaa{(p+q)}{(p+q)}{matrix} \(M\) with \tppa{block} \(A\), then \(C=B^\ad\) and the matrices \(A\) and \(M\schca A\defeq D-CA^\mpi B\) are both \tnnH{}.
\erem

\bremnl{\tcf{}~\zitaa{MR1288752}{\ccor{4.2.13}{245}}}{A.R.AxB>=0}
 If \(A\in\Cggo{m}\) and \(B\in\Cggq\), then \(A\kp B\in\Cggo{mq}\).
\erem

\breml{A.R.rA<rB}%
 If \(A,B\in\CHq\) fulfill \(\eta B\lleq A\lleq\theta B\) for some \(\eta,\theta\in\R\), then \(\nul{B}\subseteq\nul{A}\) and \(\ran{A}\subseteq\ran{B}\).
\erem

 The following result on Schur complements \eqref{E/A} is a slight modification of~\zitaa{MR2160825}{\cthm{1.3}{21}}:
 
\bleml{ab.R1038}
 Let \(M=\tmat{A & B\\ C & D}\) be the \tbr{} of a complex \taaa{(p+q)}{(r+s)}{matrix} \(M\) with \taaa{p}{r}{Block} \(A\), let \(X\in\Coo{m}{p}\), let \(Y\in\Coo{n}{p}\), let \(Z\in\Coo{n}{q}\), let \(U\in\Coo{r}{u}\), let \(V\in\Coo{r}{v}\), and let \(W\in\Coo{s}{v}\).
 If \(X^\ad X=\Ip\) and \(UU^\ad=\Iu{r}\), then \(N\defeq\smat{X&\Ouu{m}{q}\\Y&Z}M\smat{U&V\\\Ouu{s}{u}&W}\)
 admits the \tbr{}
\[
 N
 =
 \begin{pmat}[{|}]
  XAU&X(AV+BW)\cr\-
  (YA+ZC)U&YAV+YBW+ZCV+ZDW\cr
 \end{pmat}
\]
 and \(N\schca (XAU)=Y(\Ip-AA^\mpi)BW+ZC(\Iu{r}-A^\mpi A)V+Z(M\schca A)W\).
\elem
\bproof
 In view of \rrem{A.R.UA+V}, we have \(U\rk{XAU}^\mpi X=UU^\ad A^\mpi X^\ad X=A^\mpi\).
 Taking into account \eqref{mpi}, we hence obtain 
\[
 (YA+ZC)U\rk{XAU}^\mpi X(AV+BW)
 =YAV+YAA^\mpi BW+ZCA^\mpi AV+ZCA^\mpi BW.
\]
 By virtue of \eqref{E/A}, consequently,
\[\begin{split}
  N\schca (XAU)
  &=YBW+ZCV+ZDW-YAA^\mpi BW-ZCA^\mpi AV-ZCA^\mpi BW\\
  &=Y\rk{\Ip-AA^\mpi}BW+ZC\rk{\Iu{r}-A^\mpi A}V+Z(M\schca A)W.\qedhere
\end{split}\]
\eproof

 A square matrix \(A\) is called \emph{unipotent}, if its difference \(A-\EM\) with the identity matrix \(\EM\) is nilpotent, \tie{}, if \(\rk{A-\EM}^k=\NM\) for some \(k\in\N\).
 Obviously, the following block triangular matrices have this property:
 
\bnotal{ab.N1308}
 For each \(n\in\NO\) denote by \(\nudpu{n}\) (\tresp{}, \(\nodpu{n}\)) the set of all lower (\tresp{}, upper) \tppa{block} triangular matrices belonging to \(\Coo{(n+1)p}{ (n+1)p}\) with matrices \(\Ip\) on its block main diagonal.
\enota

\breml{R1547}
 For each \(n\in\NO\) the sets \(\nudpu{n}\) and \(\nodpu{n}\) are both subgroups of the general linear group of invertible complex \taaa{(n+1)p}{(n+1)p}{matrices}.
\erem

\breml{R1553}
 Let \(n\in\NO\) and let \(M\in\nudpu{n}\cup\nodpu{n}\).
 Then \(\det M=1\).
\erem

\breml{R1029}
 Let \(\ell,m\in\N\), let \(M=\tmat{A&B\\C&D}\) be the \tbr{} of a complex \taaa{(n+1)p}{(n+1)p}{matrix} \(M\) with \taaa{\ell p}{\ell p}{Block} \(A\), and let \(n\defeq\ell+m-1\).
 Then:
\benui
 \il{R1029.a} \(M\in\nudpu{n}\) if and only if \(A\in\nudpu{\ell-1}\), \(D\in\nudpu{m-1}\), and \(B=\Ouu{\ell p}{mp}\).
 \il{R1029.b} \(M\in\nodpu{n}\) if and only if \(A\in\nodpu{\ell-1}\), \(D\in\nodpu{m-1}\), and \(C=\Ouu{mp}{\ell p}\).
\eenui
\erem

\bremnl{\tcf{}~\zitaa{MR2805417}{\clem{A.3}{509}}}{091.YM}
 Let \(n\in\NO\), let \(D_0,D_1,\dotsc,D_n\) and \(E_0,E_1,\dotsc,E_n\) be complex \tpqa{matrices}, let \(\mathbb{L},\mathbb{M}\in \nudpu{n}\), and let \(\mathbb{U},\mathbb{V}\in\nodqu{n}\).
 If
\(
 \mathbb{L}\cdot\diag\rk{D_0,D_1,\dotsc,D_n}\cdot\mathbb{U}
 =\mathbb{M}\cdot\diag\rk{E_0,E_1,\dotsc,E_n}\cdot\mathbb{V},
\)
 then \(D_j=E_j\) for all \(j\in\mn{0}{n}\).
\erem

 Using the classes \(\nudpu{n}\) and \(\nodqu{n}\), we can introduce an equivalence relation for complex \taaa{(n+1)p}{(n+1)q}{matrices}:

\bnotal{A.N.ldusim}
 If \(n\in\NO\) and \(\mathbb{A},\mathbb{B}\) are two complex \taaa{(n+1)p}{(n+1)q}{matrices}, then we write \(\mathbb{A}\ldua{n}{p}{q}\mathbb{B}\) if there exist matrices \(\mathbb{L}\in \nudpu{n}\) and \(\mathbb{U}\in\nodqu{n}\) such that \(\mathbb{B}=\mathbb{L}\mathbb{A}\mathbb{U}\).
 If the corresponding (block) sizes are clear from the context, we will omit the indices and write \(\mathbb{A}\ldu\mathbb{B}\).
\enota

\breml{A.R.ldueq}
 \rrem{R1547} shows that, for a fixed number \(n\in\NO\), the relation \(\ldua{n}{p}{q}\) is an equivalence relation on the set of complex \taaa{(n+1)p}{(n+1)q}{matrices}.
\erem

 In view of \rrem{R1029}, we have furthermore:
 
\breml{A.R.ldudiag}
 Let \(\ell,m\in\NO\), let \(\mathbb{A},\mathbb{B}\) be two complex \taaa{(\ell+1)p}{(\ell+1)q}{matrices}, and let \(\mathbb{X},\mathbb{Y}\) be two complex \taaa{(m+1)p}{(m+1)q}{matrices}.
 If \(\mathbb{A}\ldua{\ell}{p}{q}\mathbb{B}\) and \(\mathbb{X}\ldua{m}{p}{q}\mathbb{Y}\), then \(\rk{\zdiag{\mathbb{A}}{\mathbb{X}}}\ldua{n}{p}{q}\rk{\zdiag{\mathbb{B}}{\mathbb{Y}}}\), where \(n\defeq\ell+m+1\).
\erem

\nocite{MR3745463,MR3377989,MR3704346}
\bibliography{178arxiv}
\bibliographystyle{bababbrv}

\vfill\noindent
\begin{minipage}{0.5\textwidth}
 Universit\"at Leipzig\\
 Fakult\"at f\"ur Mathematik und Informatik\\
 PF~10~09~20\\
 D-04009~Leipzig
\end{minipage}
\begin{minipage}{0.49\textwidth}
 \begin{flushright}
  \texttt{
   fritzsche@math.uni-leipzig.de\\
   kirstein@math.uni-leipzig.de
  } 
 \end{flushright}
\end{minipage}

\bigskip\noindent
\begin{minipage}{0.5\textwidth}
 Universit\"at Leipzig\\
 PF~10~09~20\\
 D-04009~Leipzig
\end{minipage}
\begin{minipage}{0.49\textwidth}
 \begin{flushright}
  \texttt{
   conrad.maedler@uni-leipzig.de
  }
 \end{flushright}
\end{minipage}

\end{document}